\documentclass[final,leqno,onefignum,onetabnum]{siamltex1213}
\usepackage[utf8]{inputenc}

\usepackage[english]{babel}
\usepackage[notcite,notref]{showkeys}
\usepackage{fullpage,epsfig}
\usepackage{enumerate}
\usepackage{amsmath}
\usepackage{stmaryrd}
\usepackage{amssymb}
\usepackage{latexsym}
\usepackage{accents}
\usepackage{mathrsfs}
\usepackage{subfigure}
\usepackage{color}
\usepackage[colorinlistoftodos]{todonotes}

\numberwithin{equation}{section}

\newtheorem{example}{Example}[section]
\newtheorem{remark}[example]{Remark}
\newtheorem{problem}[theorem]{Open problem}

\definecolor{Green}{RGB}{0,0, 0}

\newcommand{\be}{\begin{eqnarray}}

\newcommand{\ee}{\end{eqnarray}}
\renewcommand{\O}{\Omega}
\newcommand{\R}{{\mathbb R}}
\newcommand{\N}{{\mathbb N}}
\newcommand{\cN}{{\mathcal N}}
\newcommand{\md}{{\rm d}}
\newcommand{\ben}{\begin{eqnarray*}}
\newcommand{\een}{\end{eqnarray*}}
\newcommand{\ra}{\right\rangle}

\newcommand{\la}{\left\langle}

\newcommand{\wto}{{\rightharpoonup}}
\newcommand{\cof}{{\rm Cof }}
\newcommand{\tr}{{\rm tr}}
\newcommand{\wstar}{{\stackrel{*}{\rightharpoonup}}\,}

\renewcommand{\ker}{\text{ker}\,\mathcal{A}}
\newcommand{\FIRST}{\color{black}}
\newcommand{\SECOND}{\color{black}}
\newcommand{\WE}{\color{black}}
\newcommand{\EOR}{\color{black}}

\begin{document}
\begin{center}
{\LARGE\bf Weak lower semicontinuity of integral functionals and applications } 

\medskip

{\large Barbora Bene\v{s}ov\'{a}}\footnote{Institute of Mathematics, University of W\"{u}rzburg, Emil-
Fischer-Stra\ss e 40, D-97074 W\"{u}rzburg, Germany}  {\large \& Martin Kru\v{z}\'{\i}k}\footnote{Institute of Information Theory and Automation of the CAS, Pod vod\'{a}renskou
v\v{e}\v{z}\'{\i}~4, CZ-182~08~Praha~8, Czech Republic and
Faculty of Civil Engineering, Czech Technical
University, Th\'{a}kurova 7, CZ-166~ 29~Praha~6, Czech Republic}

\end{center}

\medskip

{\small
{\hfill\begin{minipage}[c]{18em}
``Nothing takes place in the world whose meaning is not that of some maximum or minimum.''\\
\hspace*{.7em} \sc {Leonhard Paul Euler} (1707--1783)
\end{minipage}}
}

\medskip

\begin{center}
\end{center}

\smallskip

\begin{abstract}
Minimization  is a reoccurring theme in many mathematical disciplines ranging from pure to applied ones. Of particular importance  is the minimization of integral functionals that  is studied within the calculus of variations. Proofs of 
the existence of minimizers usually rely on a fine property of the involved functional called  weak lower semicontinuity. While  early  studies of lower semicontinuity  go back to the beginning of the 20th century  the  milestones of the modern theory  were set by
C.B.~Morrey Jr. \cite{morrey-orig} in 1952 and  N.G.~Meyers  \cite{meyers} in 1965. 
We recapitulate the   development on this topic from then  on. Special attention is paid to signed integrands and to  applications in continuum mechanics of solids. In particular, we review 
the concept of  polyconvexity and special properties of (sub)determinants with respect to weak lower semicontinuity.
Besides, we emphasize some recent progress in lower semicontinuity of functionals along sequences satisfying differential and algebraic constraints which have applications in elasticity to ensure injectivity and orientation-preservation of  deformations. Finally, we outline generalization of these results to more general first-order partial differential operators and make some suggestions for further reading.

\end{abstract}
\smallskip



\section{Introduction}
{\color{Green} Many tasks in the world surrounding us can be mathematically formulated as minimization or maximization problems. For example, in physics we minimize  the energy, in economy one tries to minimize the cost and maximize the profit, entrepreneurs may try to minimize the investment risk. In addition, minimization problems appear in many  more specific tasks: in a fitting procedure, or more generally inverse problems, one tries to minimize the deviation of the model prediction from the experimental observation or training of a neuronal network is based on minimizing a suitable cost function. 

In a very general manner, we may express these problems as   
\begin{equation}
\text{minimize $I$ over $\mathcal{Y}$}\ ,
\label{var-problem-intro}
\end{equation}
where $\mathcal{Y}$ is a set  over which the minimum is sought and $I:\mathcal{Y}\to\R$ is a functional whose  meaning  may be the energy, cost, risk, or gain, for instance. From the mathematical point of view, two questions are immediate when inspecting problem \eqref{var-problem-intro}:} {\color{black} firstly  whether \eqref{var-problem-intro} \emph{is solvable}, that is if $I$ possesses minimizers on $\mathcal{Y}$, and secondly \emph{ how to find a solution (i.e. a minimizer)} to \eqref{var-problem-intro}.}

{\color{Green} \emph{Calculus of variations} is devoted to solving \eqref{var-problem-intro}  when  $\mathcal{Y}$ is (a subset) of an infinite-dimensional vector space.} Its starting point may have been a question of Johann Bernoulli on which curve a mass point will descent the fastest in a gravitational field; the so-called {\it brachistochrone problem}. {\color{Green} In the most typical situation (that covers the brachistochone problem in particular), $I$ in \eqref{var-problem-intro} is an integral functional depending on functions $u: \Omega \to \R^m$ with $\Omega \subset \R^n$  and their derivatives.} {\color{black} In the easiest case, in which $n=m=1$, $\Omega = [a,b]$, and $f:\O\times\R\times\R\to\R$ is a suitable integrand,  the functional reads  
\begin{equation}
I(u): = \int_a^b f(x,u(x),u'(x))\,\md x \quad \text{ with $u(a)=u_a$ and $u(b)=u_b$,}
\label{I-1D}
\end{equation} 
where $u_a$ and $u_b$  are given boundary data. The task is to either solve \eqref{var-problem-intro} or  at least to prove existence of minimizers.} 

Foundations of the calculus of variations were laid down in the 18th century by L.P.~Euler and J.L.~Lagrange who also realized its important connections to physics and to mechanics. {\color{black} These early works quite naturally concentrated on the question on how to find (candidates for) solutions of \eqref{var-problem-intro}. The classical method to do so, is to consider so-called \emph{variations}. Indeed, if $u_0$ is a minimizer of $I$ then 
\begin{equation}
I(u_0) \leq I(u_0 + \varepsilon \varphi) \qquad \text{for all $\varphi \in C^\infty_0([a,b])$,}
\label{variatons}
\end{equation}
where $\varepsilon \varphi$ is called a variation of the minimizer. Now, assume that $f$ is twice continuously differentiable and $u_0 \in C^2([a,b])$, then by the classical calculus \eqref{variatons} implies that $\frac{\md}{\md \varepsilon} I(u_0(x) + \varepsilon \varphi(x)){\big\vert}_{\varepsilon = 0} $ vanishes for all $\varphi \in C^\infty_0([a,b])$. This is equivalent to solving
\begin{equation}
\frac{\partial f}{\partial r}(x,u_0,u_0')- \frac{\md }{\md x}\frac{\partial f}{\partial s}(x,u_0,u_0') =0 \quad \text{on $[a,b]$,}
\label{E-L-eq}
\end{equation}
where $\frac{\partial f}{\partial r}$ and $\frac{\partial f}{\partial s}$ denote the partial derivative of $f$ with respect to the second and third variable, respectively. Equation \eqref{E-L-eq} is referred to as the \emph{Euler-Lagrange equation} and solving it is the classical path to finding solutions of \eqref{var-problem-intro}. Of course, any critical point of $I$ (and not only the minimizer) is a solution to \eqref{E-L-eq} but solving  \eqref{E-L-eq} is still an efficient approach to \eqref{var-problem-intro} at least in a situation in which all critical points are minimizers, for example if $f$ is convex. For more details,  see for example the book by Bolza \cite{bolza}.}

Nevertheless, solving the Euler-Lagrange equation naturally relies on smoothness properties of $f$ which might not be available. {\color{Green} Therefore, it is often advantageous to address existence of solutions to \eqref{var-problem-intro} in a non-constructive way by using suitable compactness properties of $\mathcal{Y}$ and continuity properties of $I$. For example, if $\mathcal{Y}$ is a bounded closed interval of reals and $I: \mathcal{Y} \to \R$ is a function then \eqref{var-problem-intro} has a solution whenever $I$ is continuous.  This observation }goes back to Bernard Bolzano who proved it in his work ``Function Theory'' in  1830 {\color{Green} and} is called the {\it Extreme Value Theorem}. Later on, it was independently shown by Karl Weierstrass around 1860. The main ingredient of the proof, namely the fact that one can extract a convergent subsequence from a closed bounded interval of reals, is nowadays known as the Bolzano-Weierstrass theorem. 

{\color{black}The results of Bolzano and Weierstrass easily extend to the situation when $\mathcal{Y}$ is a bounded and closed set of a \emph{finite-dimensional} vector space. However, they cannot be generalized to the situation in which, for example, $\mathcal{Y}$ is a ball in an infinite dimensional vector space since the Bolzano-Weierstrass theorem is false in this case. In fact, being able to extract a convergent subsequence from a sequence of elements in the unit ball of a normed vector space $\mathcal{X}$ is \emph{equivalent} to $\mathcal{X}$ being finite dimensional. This is a classical result attributed to F.~Riesz.} 

{\color{Green} Thus, the only hope to transfer a variant of the Bolzano-Weierstrass theorem to infinite dimensional spaces is to seek compactness in a ``weaker" topology than the one induced by the norm.} This possibility has been opened by Riesz and Hilbert who used the \emph{weak topology} on Hilbert spaces from the beginning of the 20th century and Stefan Banach who defined it on other normed spaces around 1929 \cite{pietsch, wehausen}.

\medskip

{\color{Green}
\begin{definition}
Let $\mathcal{X}$ be a Banach space and $\mathcal{X}'$ its dual. We say that a sequence $\{u_k\}_{k \in \N} \subset \mathcal{X}$ converges \emph{weakly} in $\mathcal{X}$ to $u \in \mathcal{X}$ if 
$$
\psi(u_k) \to \psi(u) \quad \text{for all $\psi \in \mathcal{X}'$ and we write that} \quad u_k \rightharpoonup u.
$$
Similarly a sequence $\{v_k\}_{k \in \N} \subset \mathcal{X}'$ converges \emph{weakly*} in $\mathcal{X}'$ to $v \in \mathcal{X}'$ if 
$$
v_k(\varphi) \to v(\varphi) \quad \text{for all $\varphi \in \mathcal{X}$ and we write that} \quad v_k \wstar v.
$$
\end{definition}

\medskip

A crucial property of the weak topology is that it allows for a generalization of the Bolzano-Weierstrass theorem to infinite dimensional vector spaces. Indeed, take $\mathcal{X}'$ the dual to a Banach space $\mathcal{X}$. Then bounded subsets of $\mathcal{X}'$ are \emph{precompact in the weak* topology} by the Banach-Alaoglu theorem, even though they are generically not compact if $\mathcal{X}'$ is infinitely dimensional. As an immediate consequence, we have that bounded subsets of a reflexive Banach  space $\mathcal{X}$ are precompact in the weak topology. 

Having the weak topology at hand, a generalization of the Bolzano extreme value theorem becomes possible and is today known as the \emph{direct method of the calculus of variations}.} This algorithm was proposed by David Hilbert around 1900, to show (in a non-constructive way) the existence of a solution to the minimization problem \eqref{var-problem-intro}.
 It consists of three steps: 
 
 \smallskip
 
 \begin{enumerate}
 \item We find a minimizing sequence along which $I$ converges to its infimum on $\mathcal{Y}$. 
 \item We show that a subsequence of the  minimizing sequence converges to an element of $\mathcal{Y}$ in some topology $\tau$. 
 \item We prove that this limit element 
is a minimizer. 
\end{enumerate}

\smallskip
{\color{Green} The first step of the direct  method is easily handled if the infimum of $I$ is finite.
For the second step,  the appropriate choice of the topology $\tau$ is crucial.  
 In the most typical situation, the set $\mathcal{Y}$ is a subset of a Banach space or its dual  and} $\tau$ is either the weak or the weak$^*$ topology. {\color{Green} In this case, if $\mathcal{Y}$ is bounded, existence of a converging subsequence of the minimizing sequence is immediate from the Banach-Alaoglu theorem. If $\mathcal{Y}$ is not bounded, the usual remedy is to realize that the minimizer can only lie in a bounded subset of $\mathcal{Y}$ due to \emph{coercivity} of $I$.} {\color{black} Coercivity refers to the property of $I$ that it takes arbitrarily large values  if the norm of its argument is sufficiently large. More precisely, we say that $I$ is coercive if 
\begin{align}\lim_{\|u\|\to\infty} I(u)=\infty\ .\end{align}  This allows us to say that 
all minimizers of $I$ are contained in some closed ball centered at the origin.  }

{\color{Green} The third step of the direct needs to rely on suitable semicontinuity properties of $I$; a sufficient and widely used condition is the} \emph{(sequential) lower semicontinuity} of $I$ with respect to the {\color{Green} weak/weak*} topology:

\medskip

\begin{definition}
\label{def-wlsc}
Let $\mathcal{Y}$ be a subset of a  Banach space. We say that the functional $I: \mathcal{Y} \to \mathbb{R}$ is (sequentially) weakly/weakly* lower-semicontinuous on $\mathcal{Y}$ if for every sequence $\{u_k\}_{k \in \mathbb{N}} \subset \mathcal{Y}$ converging weakly/weakly* to $u\in\mathcal{Y}$, we have that
$$
I(u) \leq \liminf_{k \to \infty} I(u_k).
$$
\end{definition}
\medskip

 {\color{Green}

If $I$ is not weak/weak* lower semicontinuous solutions to \eqref{var-problem-intro} need not to exist. However, weak lower semicontinuity of $I$ is not a necessary condition for the existence of minimizers. These facts are demonstrated by the following example.
\medskip

\begin{example}
\label{ex-zigzag}
Consider the following special case of \eqref{I-1D}:
\begin{align}
I(u) = \int_0^1 \big(1-(u'(x))^2\big)^2 + (u(x))^2\,\md x  \end{align}
with 
$$\mathcal{Y}:=\{u\in W^{1,\infty}(0,1);\ -1\le u'\le 1,\ u(0)=u(1)=0\}\ .$$

We can see, for example by the Lebesgue dominated convergence theorem, that $I$ is continuous on $W^{1,\infty}(0,1)$ but it is not weakly lower semicontinuous. To show this, define  
$$
u(x) = \begin{cases} x & \text{if $0\le x\le 1/2$}\\ 
-x+1 & \text{if $1/2\le x\le 1$}\\
\end{cases} 
$$
and extend it periodically to the whole $\R$. Let $u_k(x):=k^{-1}u(kx)$ for all $k\in\N$ and all $x\in\R$. Notice that 
$\{u_k\}_{k\in\N}\subset\mathcal{Y}$.

The sequence of ``zig-zag" functions $\{u_k\}_{k\in \N}$  converges weakly* to zero in $W^{1,\infty}(0,1)$. It is not hard to see that $I(u_k) \to 0$ for $k\to\infty$  but 
$$
1  = I(0) > \lim_{k \to \infty } I(u_k) = 0;
$$
so that $I$ is not weakly* lower semicontinuous on $W^{1,\infty}(0,1)$ and, in fact, no minimizer exists in this case.

Indeed, $0=\inf_{\mathcal{Y}} I\ne\min_{\mathcal{Y}} I$ because $I\ge 0$ and $I(u_k)\to 0$, so that $0=\inf_{\mathcal{Y}} I$. However, $I(u)>0$ for every $u\in\mathcal{Y}$, for otherwise it would mean that we could find a Lipschitz function whose derivative is $\pm 1$ almost everywhere on $(0,1)$ but the function value is identically zero.

If we, however, consider a slight modification of $\mathcal{Y}$ by changing the boundary condition at $x=1$, and consider 
$$ \mathcal{Y}_1:=\{u\in W^{1,\infty}(0,1);\ -1\le u'\le 1,\ u(0)=0,\ u(1)=1\}\ $$
then $\min_{\mathcal{Y}_1} I=1/3$ and the unique minimizer is $u(x)=x$ for $x\in(0,1)$.  

Firstly, this shows that weak/weak* lower semicontinuity of $I$ is not necessary for the existence of a minimizer, and, secondly, it stresses the influence of  boundary conditions on the solvability of \eqref{var-problem-intro}. This phenomenon is even more pronounced in higher dimensions.

\end{example}
}
\medskip

Although the study of weak lower semicontinuity is motivated by understanding minimization problems, it has become an independent subject in mathematical literature that has been studied for its own right. 
{\color{Green} In the case of integral functionals as in \eqref{I-1D}, further properties of the integrand besides continuity are needed to assure weak/weak* lower-semicontinuity: the right additional property is always some type of \emph{convexity} of $f$. Indeed, notice that $I$ in Example \ref{ex-zigzag} is not convex.} 

{\color{Green} The importance of convexity for weak/weak* lower semicontinuity for integral functionals has been discovered by Tonelli in 1920 \cite{tonelli}, who pioneered the study of lower semicontinuity of an integral functional rather than studying the associated Euler-Lagrange equation. Tonelli considered $f:\O\times\R\times\R$ in  \eqref{I-1D} that is  {\color{black} twice continuously differentiable} and showed that $I$ is  lower semicontinuous subject to a ``convergence of curves''\footnote{Notice that the notion of weak topology was invented later than Tonelli's studies.}} if and only if $f$ is  convex in its last variable, i.e., in the derivative $u'$. Later, several authors generalized this result to  functions in $W^{1,\infty}(\Omega; \R)$ with $\Omega \subset \R^n$ and $n > 1$; see for example Serrin \cite{serrin}, where differentiability properties of $f$ were removed from assumptions and $f$ {\color{Green} was only assumed to be continuous}, and Marcellini and Sbordone \cite{marcellini-sbordone} who {\color{Green} allowed for} Carath\'{e}odory integrands\footnote{ i.e. $f(x,\cdot,\cdot)$ is continuous for almost all $x\in\O$ and $f(\cdot, s,A)$ is measurable for all $(s,A)\in\R\times\R$}. {\color{Green} Similarly as in this one-dimensional situation}, relaxing smoothness/continuity assumptions of $f$ will be a re-occurring theme throughout this review in which we focus on the higher-dimensional case.  

{\color{black} Let us now allow the function $u$ to be vector-valued, i.e., $u \in W^{1,\infty}(\Omega; \R^m)$ with $\Omega \subset \R^n$ and $n>1$ as well as  $m>1$ and consider an integral functional of the form
\begin{align}\label{functional0}
I(u):=\int_\O f(x,u(x),\nabla u(x))\,\md x\ .
\end{align}}
In this case, the convexity hypothesis turns out to be sufficient for weak/weak* lower semicontinuity but unnecessary. A suitable condition, termed  quasiconvexity, was introduced by Morrey \cite{morrey-orig}. 

\medskip

\begin{definition}
\label{def-quasiconvexity}
 Let $\O\subset\R^n$ be a bounded Lipschitz  domain with the Lebesgue measure $\mathcal{L}^n(\O)$. A function $f:\R^{m\times n}\to\R$ is quasiconvex {\color{Green} at $A\in\R^{m\times n}$ if for} every $\varphi\in W^{1,\infty}_0(\O;\R^m)$
\be\label{quasiconvexity}
f(A)\mathcal{L}^n(\O)\le \int_\O f(A+\nabla \varphi(x))\,\md x\ .\ee
{\color{Green} The function $f$ is termed quasiconvex if it is quasiconvex in all $A\in\R^{m\times n}$. }
\end{definition}

\medskip

{\color{black} Quasiconvexity is implied by convexity and can be understood as, roughly speaking, convexity over gradients. Indeed take a convex function $f:\R^{m\times n}\to\R$. Then for some arbitrary $A \in \R^{m\times n}$ fixed  and every $B \in \R^{m\times n}$, we know that $f(A + B) \geq  f(A) + g(A) {\cdot} B$; i.e., we can find an affine function that touches $f$ at $A$ and its values are not greater than $f$ (in fact, this can be found by taking $g(\cdot)$ is one element of the subdifferential of $f$).  Let us now take some arbitrary $\varphi\in W^{1,\infty}_0(\O;\R^m)$ and plug in $\nabla \varphi (x)$ in the position of $B$ and take an average of the inequality over $\Omega$ to get that
$$
\frac{1}{\mathcal{L}^n(\O)} \int_\Omega f(A+\nabla \varphi(x))\,\md x \geq \frac{1}{\mathcal{L}^n(\O)} \int_\Omega f(A)\,\md x + \frac{1}{\mathcal{L}^n(\O)} \int_\Omega  g(A) {\cdot} \nabla \varphi(x)\,\md x\ ,
$$
where the last integral vanishes due to integration by parts because $\varphi=0$ on $\partial\O$ so that we truly obtain \eqref{quasiconvexity}. We also note that quasiconvex fuctions are continuous \cite{dacorogna}.
}

Morrey showed, under strong regularity assumptions on $f$, that $I$ from \eqref{functional0} is weakly lower semicontinuous in $W^{1,\infty}(\O;\R^m)$ if and only if $f$ is quasiconvex in the last variable (i.e. in the gradient). 
{\color{black} To see how quasiconvexity is used in the proof of lower-semicontinuity, let us consider the following simplified example:

\medskip

\begin{example}
\label{ex-intro-quasi}
Assume that $\{u_k\}_{k \in \N} \subset W^{1,\infty}(\O;\R^m)$ is such that $u_k \wstar u$ with $u(x)=Ax$ for  some fixed matrix  $A\in\R^{m\times n}$. We show how in this case weak* lower semicontinuity on $W^{1,\infty}(\O;\R^m)$ is obtained for 
$$
\tilde{I}(u) := \int_\Omega f(\nabla u(x)) \md x\ ,
$$
for $f:\R^{m\times n}\to\R$ quasiconvex. 
To this end, let us take a smooth cut-off function $\eta_\ell:\Omega \to \R$ such that $\eta_\ell = 1$ on $\Omega_\ell$ and $\eta_\ell = 0$ on $\partial \Omega$, where $\Omega_\ell \subset \Omega$ is a Lipschitz domain satisfying that $\mathcal{L}^n(\Omega \setminus \Omega_\ell) \leq \frac{1}{\ell}$. We may find $\eta_\ell$ in such a way that $|\nabla \eta_\ell| \leq C \ell$ uniformly on $\Omega$, where $C$ is a constant that depends just on $\Omega$. Let us now define 
$$
u_{k, \ell}(x) = \eta_\ell u_k + (1-\eta_\ell) Ax \quad \text{so that} \quad \nabla u_{k, \ell}(x) = \eta_\ell \nabla u_k + (1-\eta_\ell) A + (u_k-Ax)\otimes\nabla \eta_\ell ;
$$
notice that $u_{k, \ell}$ coincides with $u_k$ on $\Omega_\ell$. Now, since $u_k \to u$ strongly in $L^\infty(\O;\R^m)$, we may choose a subsequence of $k$'s, labeled $k(\ell)$,  such that $(u_{k(\ell)}-Ax)\otimes\nabla \eta_\ell$ stays uniformly bounded (whence $u_{k(\ell), \ell}$ is bounded in $W^{1,\infty}(\O;\R^m)$). Due to the fact that  $u_{k(\ell), \ell}(x) = Ax$ on $\partial \Omega$ we get from \eqref{quasiconvexity} that
\begin{equation}
f(A)\mathcal{L}^n(\O) \leq \int_\Omega f(\nabla u_{k(\ell), \ell}(x)) \md x = \int_{\Omega} f(\nabla u_{k(\ell)}(x)) \md x + \int_{\Omega\setminus\Omega_\ell} f(\nabla u_{k(\ell), \ell}(x))-f(\nabla u_{k(\ell)}(x))  \md x.
\label{example-intro-quasi}
\end{equation}
As $f$ is continuous and $\{\nabla u_{k(\ell), \ell}\}_{\ell\in\N}$ is uniformly bounded on $\Omega$, so is $f(\nabla u_{k(\ell),\ell})-f(\nabla u_{k(\ell)})$ and thus the last integral in \eqref{example-intro-quasi} vanishes as $\ell \to \infty$. Therefore, taking the limit $\ell \to \infty$  yields the claim.
\end{example}

\medskip

}

{\color{Green} The results of Morrey} were generalized more than fifty years ago, in 1965, by Norman G.~Meyers  in his seminal paper \cite{meyers}. Taking $k\in\N$ and $1\le p\le +\infty$ he investigated  the $W^{k,p}$-weak (weak* if $p=+\infty$) lower semicontinuity of integral functionals of the form 
\begin{align}\label{functional1}
I(u):=\int_\O f(x,u(x),\nabla u(x),\ldots, \nabla^{k}u(x))\,\md x\ ,
\end{align}
where $\O\subset\R^n$ is a bounded domain and $u:\O\to\R^m$ is  a mapping possessing (weak) derivatives up to the order $k\in\N$. The function $f$ was supposed to be continuous in all  its arguments. Since now higher gradients {\color{Green} (and not only the first ones)}  are considered, the definition of quasiconvexity also needs to be generalized accordingly;   {\color{Green} see Section \ref{result-meyers} for details.}

{\color{Green} Moreover}, more generally than in Morrey's work, the function $f$ is not necessarily  bounded from bellow in \cite{meyers}. From this, additional difficulties arise and, in fact, quasiconvexity is no longer a sufficient condition for weak lower semicontinuity (cf. Section \ref{result-meyers}). In addition, the regularity assumptions on the integrand in \eqref{functional1} were weakened in Meyers' work.

The motivation for studying functionals of the type \eqref{functional1} is twofold: from the point of view of applications in continuum mechanics it is reasonable to let $f$ depend also on higher-order gradients since their appearance in the energy usually models interfacial energies or multipolar elastic materials \cite{Green-rivlin}. Another reason might be to consider deformation-gradient dependent surface loads \cite{ball-curie-olver}. On the other hand, not assuming a constant lower bound on $f$ is important to consider for mathematical completeness. Additionally, integrands  of the type $f(A):=\det A$, which are unbounded from below, are of crucial importance in continuum mechanics.  

Meyers' main results are necessary and sufficient conditions on $f$ so that $I$ is weakly lower semicontinuous on $W^{k,p}(\O;\R^m)$. We review these results in Section \ref{result-meyers}. He first discusses the problem on {\color{black} $W^{k,\infty}(\O;\R^m)$}, where quasiconvexity in the highest-order gradient (cf. Theorem \ref{meyers1}) turns out to be a necessary and sufficient  condition for weak*-lower semicontinuity. {\color{black} Lower semicontinuity on $W^{k,p}(\O;\R^m)$ with $1< p<+\infty$ is}, however,  much more subtle, and an additional condition (cf.~Theorem \ref{meyers4} and Section \ref{sec-BdCon}) is needed.

Since the appearance of Meyers' work, significant progress has been achieved with respect to the characterization of weak lower semicontinuity of functionals of the type \eqref{functional1}. In particular, for $k=1$ in \eqref{functional1} the additional condition for sequential weak lower semicontinuity was characterized more explicitly and results relaxing Meyers' continuity assumptions were obtained for functionals bounded from below; cf. Section \ref{result-meyers}. 

Moreover, it has been identified for which functions $f$ the functional $I$ in \eqref{functional1} is even weakly continuous (see Section \ref{sec-NullLag}) -- these functions are the so-called null Lagrangians -- and this knowledge led to the notion of polyconvexity (see Section \ref{sec-polyconvexity}) that is sufficient for weak lower semicontinuity and of particular importance in mathematical elasticity. In fact, quasiconvexity, which is, for a large class of integrands, the necessary and sufficient condition for weak lower semicontinuity  is not well-suited for elasticity. We explain this issue in Section \ref{sec-quasicon-elas} and review some recent progress in this field. Null Lagrangians have also been identified for functionals defined on the boundary (see Section \ref{sec-NullBd}). Finally, we review weak lower semi-continuity results for functionals depending on maps that satisfy general differential constraints in Section \ref{sec-A-quasi} and we conclude with some suggestions for further reading in Section~\ref{reading}.  

\section{Notation}
\label{sec-notation}
In this section, we summarize the notation that shall be used throughout the paper. It largely coincides with the one  in \cite{ball-curie-olver}. In what follows, $\O\subset\R^n$ is a bounded  domain  the boundary of which is  Lipschitz or smoother. This domain is mapped to a set in $\R^m$ by means of a mapping $u:\O\to\R^m$.  

Let $\N$ be the set of natural numbers and $\N_0:=\N\cup\{0\}$.
If $J:=(j_1,\ldots, j_n)\in\N_0^n$ and $K:=(k_1,\ldots, k_n)\in\N_0^n$ are two multiindices we define $J\pm K:=(j_1\pm k_1,\ldots, j_n\pm k_n)$, further $|J|=\sum_{i=1}^nj_i$,
$J!:=\Pi_{i=1}^nj_i!$, and we say that $J\le K$ if $j_i\le k_i$ for all $i$. Then we also define ${J\choose K}:=\frac{J!}{K!(K-J)!}$, $\partial u^j_{K}:=\frac{\partial^{k_1}\ldots\partial^{k_n}}{\partial x_1^{k_1}\ldots\partial x_n^{k_n}} u^j$, 
$x_K=x^K:=x_1^{k_1}\ldots x^{k_n}_n$, and $(-D)^K:=\frac{(-\partial)^{k_1}\ldots (-\partial)^{k_n}}{\partial x_1^{k_1}\ldots\partial x_n^{k_n}}$. 

We will  work with the space of {\color{black} higher-order matrices} $X=X(n,m,k)$ with the dimension 
$m{n+k-1\choose k}$. This is the space of {\color{black} (higher-order) matrices} $M=(M^i_K)$ for $1\le i \le m$  and $|K|=k$. Similarly, $Y=Y(n,m,k)$ is a space of {\color{black} (higher-order) matrices} $M=(M^i_K)$ for $1\le i \le m$  and $|K|\le k$. Its dimension is $m{n+k\choose k}$. We denote the elements of $X(n,m,k)$ by $A^k$ while $A^{[k]} = (A, A^2, \ldots, A^k)$ is an element of $Y(n,m,k)$. We use an analogous notation also for gradients; thus, if $x\in\O$, then $\nabla^k u(x)\in X(n,m,k)$ while $\nabla^{[k]} u(x)\in Y(n,m,k)$.

Integrands which define the integral functional will be denoted by $f$. They will depend on $x$, $u$, and (higher-order) gradients of $u$, in general. Occasionally, we will work with integrands independent of $u$ or $x$, however, this will be clear from the context and will not cause any ambiguity.   We denote by $B(x_0, r)$ the ball of origin $x_0$ with the radius $r$ while $D_\varrho(x_0,r)$ is the half-ball with $\varrho$ being the normal of the planar component of its boundary; i.e. 
$$
D_\varrho(x_0,r):= \{x\in B(x_0,r):\, (x-x_0)\cdot\varrho<0\},
$$
and we write $D_\varrho:= D_\varrho(0,1)$.

{\color{black} For this review, we will assume that the reader is familiar with functional analysis and measure theory; in particular theory of Lebesgue and Sobolev spaces and refer for example to the books  by Rudin \cite{rudin} and  Leoni \cite{leoni}  for an introduction.} We  shall use the standard notation for the Lebesgue spaces $L^p(\O; \R^m)$ and Sobolev spaces $W^{k,p}(\O; \R^m)$. Moreover, $\mathrm{BV}(\O; \R^m)$ is the space of functions of a bounded variation. If $m=1$, we may omit the target space. If $\O$ is a bounded open domain we denote $\mathcal{M}(\O)$ the space of Radon measures on $\O$ and  $\mathcal{L}^n(\O)$ stands for  the $n$-dimensional Lebesgue measure of $\Omega$; cf.~e.g.~Halmos \cite{halmos}. Further, $\mathcal{M}^1_+(\R^{m\times n})$ is the set of probability measures on $\R^{m\times n}$.   Moreover, $\mathcal{D}(\Omega)$ is the space of infinitely differentiable functions with compact support in $\O$  and its dual $\mathcal{D}'(\Omega)$ is the space of distributions. 

If  $n=m=3$ and  $F\in\R^{3\times 3}$ the cofactor matrix $\cof F\in\R^{3\times 3}$  is  a matrix whose entries are 
signed subdeterminants of $2\times 2$ submatrices of $F$. More precisely, $[\cof F]_{ij}:=(-1)^{i+j}\det F'_{ij}$ where $F'_{ij}$ for $i,j\in\{1,2,3\}$ is a submatrix of $F$ obtained by removing the $i$-th row and $j$-th column. If $F$ is invertible, we have  $\cof F=(\det F) F^{-\top}$. Rotation matrices with determinants equal one  are denoted SO$(n)$ while orthogonal matrices with  determinants $\pm 1$  are {\color{black} denoted} O$(n)$. Additional notation needed locally in the text will be explained on particular spots.

\medskip

\section{A review of Meyers' results}
\label{result-meyers}
{ \color{black}
Within this section we review the results of Meyers' seminal paper \cite{meyers} and give generalizations of his results that were proved since the appearance of his work. As highlighted above, Meyers generalized Morrey's results \cite{morrey-orig} particularly within two respects: First he considers integral functionals of the type 
\begin{equation}\tag{\ref{functional1}}
I(u):=\int_\O f(x,u(x),\nabla u(x),\ldots, \nabla^{k}u(x))\,\md x\ ,
\end{equation}
i.e. those that depend also on higher gradients and second he allows for $f$ unbounded from below. Now if \eqref{functional1} depends on higher gradients, also the definition of quasiconvexity needs to be generalized accordingly: 

\medskip

\begin{definition}
 Let $\O\subset\R^n$ be a bounded Lipschitz domain. We say that  a function $f:X(n,m,k)\to\R$ is $k$-quasiconvex\footnote{In the original paper \cite{meyers}, quasiconvexity with respect to the $k$-th gradient is also referred to as quasiconvexity.} if
for every $A\in X(n,m,k)$ and any $\varphi\in W^{k,\infty}_0(\O;\R^m)$
\be\label{quasiconvexity1}
f(A)\mathcal{L}^n(\O)\le \int_\O f(A+\nabla^k \varphi(x))\,\md x\ .\ee 
\end{definition}

Thus, more precisely, $k$-quasiconvexity of $f$ (i.e. quasiconvexity with respect to the $k$-th gradient) means that $A^k\mapsto f(x,A^{[k-1]},A^k)$ is quasiconvex for all fixed $(x,A^{[k-1]})\in\O\times Y(m,n,k-1)$; here, cf. Section \ref{sec-notation} for notation.
\medskip

\begin{remark}
In fact, it was shown in \cite{dalmasoetal} that if $k=2$  and if $f$ satisfies  a (slightly) stronger version of $2$-quasiconvexity then $2$-quasiconvexity coincides with $1$-quasiconvexity. See \cite{cagnetti} for an analogous result with general $k$.

\end{remark}

\medskip
}

{\color{black} With this definition at hand,} Meyers proves an analogous result to the one found in  the original work of Morrey for $k=1$ \cite{morrey-orig}:

\medskip

\begin{theorem}[from \cite{meyers}]\label{meyers1}
Let $\O$ be a bounded domain and {\color{black}$f$} a continuous function. Then $I$ from \eqref{functional1} is weakly$^*$ lower semicontinuous on $W^{k,\infty}(\O;\R^m)$  if and only if it is $k$-quasiconvex {\color{black} in the last variable}.
\end{theorem}

\medskip

Nevertheless, when it comes to the case of $W^{k,p}(\O;\R^m)$ with $1<p<+\infty$   the situation is substantially more involved; in particular, because the considered integrands are not bounded from below. {\color{black} In fact, as can be seen from the definition of the class $\mathcal{F}_p(\O)$, Meyers studies weak lower semicontinuity of \eqref{functional1} on a fairly  general class of integrands including those with critical negative growth.}

\medskip

\begin{definition}[Class $\mathcal{F}_p(\Omega)$]
\label{class-F_R}
Let $\O\subset\R^n$ be a bounded domain. A continuous integrand ${\color{black}f}:\O\times Y(n,m,k)\to\R$ is said to be in the class $\mathcal{F}_p(\O)$ for $1\le p<+\infty$ if ($C>0$ is a constant depending only on ${\color{black}f}$)
\begin{itemize}
\item[(i)] ${\color{black}f}(x,A^{[k]})\le C\big(1+|A^{[k]}|\big)^p$, 
\item[(ii)] $|{\color{black}f}(x, A^{[k]}+B^{[k]})-{\color{black}f}(x,A^{[k]})|\le C\big(1+|A^{[k]}|+|B^{[k]}|\big)^{p-\gamma}|B^{[k]}|^\gamma$, for some $0<\gamma\le 1$,
\item[(iii)] $|{\color{black}f}(x+y,A^{[k]})-{\color{black}f}(x,A^{[k]})|\le (1+|A^{[k]}|)^p\eta(|y|)$ with $\eta:[0;+\infty)\to[0;+\infty)$  continuous, increasing and vanishing at zero.
\end{itemize}
\end{definition}

\medskip

{\color{black} 
\begin{remark}[Class $\mathcal{F}_p(\Omega)$ for $k=1$]
\label{F_R-k1}
Let us, for clarity, repeat the conditions given in Definition \ref{class-F_R} for the case $k=1$. In this case, the notation is much simpler so that the important features of functions in the class $\mathcal{F}_p(\Omega)$ can be seen more easily. 

We say that $f:\Omega \times \R^m \times \R^{m \times n} \to \R$  is in the class $\mathcal{F}_p(\O)$ for $1\le p<+\infty$ if 
\begin{itemize}
\item[(i)] $f(x,s,A)\le C\big(1+|s|+|A|\big)^p$, 
\item[(ii)] $|f(x, s+r,A+B)-f(x,s,A)|\le C\big(1+|s|+|r|+|A|+|B|\big)^{p-\gamma}(|r|+|B|)^\gamma$, {\color{black} for some} $0<\gamma\le 1$,
\item[(iii)] $|f(x+y,s,A)-f(x,s,A)|\le (1+|s|+|A|)^p\eta(|y|)$ with $\eta:[0;+\infty)\to[0;+\infty)$  continuous, increasing and vanishing at zero.
\end{itemize}
Above, $C>0$ is a constant depending on $f$.
\end{remark}
}

\medskip 

When setting $A^{[k]}=0$ in (ii) {\color{Green} in Definition \ref{class-F_R}} {\color{black}(or alternatively $s=0$ and $A=0$ in (ii) of Remark \ref{F_R-k1})} we get that $|{\color{black}f}(x,B^{[k]})|\le C(1+|B^{[k]}|)^p$ and thus the class $\mathcal{F}_p(\O)$ contains also non-coercive integrands and, in particular, those which decay as $A\mapsto -|A|^p$. {\color{Green} Quasiconvexity is not sufficient to prove sequential weak lower semicontinuity for such integrands. We shall devote Section \ref{sec-BdCon} to a detailed discussion of this issue and state at this point Meyers' original theorem which copes with the non-coercivity of $f$ by introducing an additional condition (item (ii) in Theorem \ref{meyers4}).} 

\medskip

\begin{theorem}\label{meyers4}
Let $\O$ be a bounded domain and ${\color{black} f}\in\mathcal{F}_p(\O)$.  Then $I$ from \eqref{functional1} is weakly lower semicontinuous on $W^{k,p}(\O;\R^m)$ with $1\le p<\infty$ if and only if the following two conditions hold simultaneously:
\begin{itemize}
\item[(i)] ${\color{black} f}(x,A^{[k-1]},\cdot)$ is $k$-quasiconvex for all values of $(x,A^{[k-1]})$,
\item[(ii)] $\liminf_{j\to\infty} I(u_j,\O')\ge -\mu(\mathcal{L}^n(\O'))$ for every subdomain $\O'\subset\O$ and every sequence $\{u_j\}_{j\in\N}\subset W^{k,p}(\O;\R^m)$ such that $u_j= u$ on $\O\setminus\O'$ and $u_j\wto u$ in $W^{k,p}(\O;\R^m)$.  Here $\mu$ is an increasing continuous function with $\mu(0)=0$ which only depends on $u$ and on $\limsup_{j\to\infty} \|u_j\|_{W^{k,p}(\Omega; \R^m)}$.
\end{itemize}
\end{theorem} 

\medskip
Above, $
I(\cdot,\O')$ denotes the functional $I$ when the integration domain $\O$ is replaced by $\O'$. We immediately see that condition (ii) is satisfied if ${\color{black}f}$ has a lower bound; for example, if ${\color{black}f}\ge 0$. {\color{Green} This is a very common case, in which Theorem \ref{meyers4} can be sharpened; we refer to Section \ref{sec-BoundedFromBelow} where this situation is handled in detail.}

\medskip

{\color{black}
\begin{remark}[Theorem \ref{meyers4} for $k=1$] If $k=1$ and $f:\Omega \times \R^m \times \R^{m \times n} \to \R$ is in $\mathcal{F}_p(\Omega)$ in the sense of Remark \ref{F_R-k1} then Theorem \ref{meyers4} assures that the functional $I$ from \eqref{functional0} is lower semicontinuous in $W^{1,p}(\Omega;\R^m)$ if f is quasiconvex in the last variable and condition (ii) in Theorem \ref{meyers4} is fulfilled with $k=1$.
\end{remark}
}

\medskip

\subsection{Understanding condition (ii) in Theorem~\ref{meyers4}} 
\label{sec-BdCon}

Condition (ii) in Theorem \ref{meyers4} is rather implicit and thus hard to verify. Nevertheless, in this section, {\color{Green} we will motivate that} it should be linked to concentrations on the boundary of the domain. To our best knowledge, this link has been fully drawn only in the case $k=1$ and for integrands ${\color{black}f}(x,u,\nabla u):={\color{black}f}(x,\nabla u)$ in \eqref{functional0}. {\color{Green} Thus, we will limit our scope  to this particular case and only point to some possible extensions at the end of the section.} 

{\color{black} In essence, (ii) in Theorem~\ref{meyers4} needs to cope with the potential non-equiintegrability of the negative part of the integrand $f$.} {\color{Green} To explain this statement in more detail, let us start with the definition of equiintegrability:

\medskip

\begin{definition}
\label{def-equiint}
We say that a sequence of functions $\{\varphi_k\}_{k \in \N} \subset L^1(\Omega)$ is equiintegrable, if for every $\varepsilon  > 0$ there is  $\delta > 0$ such that for every $\omega \subset \Omega$ with $\mathcal{L}^n(\omega) \leq \delta$ it holds that
$$
\sup_{k \in \N}  \int_\omega |\varphi_k(x)| \md x \leq \varepsilon.
$$
\end{definition}

As $L^1(\Omega)$ is not reflexive, a bounded sequence in $L^1(\Omega)$ does not necessarily contain a   weakly convergent subsequence  in $L^1(\Omega)$ (though it will always contain  a subsequence  weakly*-convergent in measures) but it follows from the Dunford-Pettis criterion \cite{d-s,fomu} that this measure is an $L^1$ function and the convergence improves from weak* in measures to the weak one in $L^1$ if and only if  the sequence is  equiintegrable. Since the failure of equiintegrability is caused by concentrations of the sequence $\{\varphi_n\}_{n \in \N}$},  we say that a sequence bounded in $L^1(\Omega)$ is concentrating if it converges weak* in measures but not weakly in $L^1(\Omega)$.

Recall that two effects may cause a sequence $\{u_n\}_{n \in N}\subset W^{1,p}(\O;\R^m)$ to converge weakly but not strongly to some limit function $u$: \emph{oscillations} and \emph{concentrations}. {\color{Green} Here, concentrations are understood in that sense that $\vert u_n \vert^p$ is a concentrating sequence. }{\color{black} In fact, it can be seen by Vitali's convergence theorem that if $\vert u_n \vert^p$ is equiintegrable (i.e. concentrations are excluded) and $u_n \to u$ a.e. in $\Omega$ (i.e. oscillations are excluded) $\{u_n\}_{n \in N}$ actually converges strongly to $u$ in $W^{1,p}(\O;\R^m)$.}

{\color{Green} Concentrations and oscillations in a sequence $\{u_n\}_{n \in N}\subset W^{1,p}(\O;\R^m)$ can be separated from each other by the so-called} decomposition lemma due to Kristensen  \cite{kristensen} and  Fonseca, M\"{u}ller, and Pedregal  \cite{fmp}.  

\medskip

\begin{lemma}[Decomposition lemma]\label{fons}
Let $1 < p<+\infty$ and  $\O\subset\R^n$ be an open bounded set and let $\{ u_k\}_{k\in\N}\subset W^{1,p}(\O;\R^m)$ be bounded. Then there is a subsequence $\{u_j\}_{j\in\N}$ and a sequence $\{z_j\}_{j\in\N}\subset W^{1,p}(\O;\R^m)$ such that
\be\label{rk}
\lim_{j\to\infty} \mathcal{L}^n(\{x\in\O;\ z_j(x)\ne u_j(x)\text{ or }  \nabla z_j(x)\ne \nabla u_j(x)\})=0
\ee
and $\{|\nabla z_j|^p\}_{j\in\N}$ is relatively weakly compact in $L^1(\O)$.
\end{lemma}

\medskip

This lemma allows us  to find, for a general  sequence bounded in $W^{1,p}(\Omega; \R^m)$,  another one, called $\{z_j\}\subset W^{1,p}(\O;\R^m)$, whose gradients  are  $p$-equiintegrable, i.e., for which $\{|\nabla z_j|^p\}$ is relatively weakly compact in $L^1(\Omega)$ and so it is a {\it purely oscillating} sequence. Thus, we decompose $u_j=z_j+w_j$, and $\{|\nabla w_j|^p\}_{j\in\N}$ tends to zero in measure for $j\to\infty$; i.e., it is a {\it purely concentrating} sequence. Roughly speaking, this means that for every weakly converging sequence in $W^{1,p}(\Omega; \R^m)$, $p>1$, we can decompose the sequence of gradients   into a purely oscillating and a purely concentrating one. Note, however, that due to \eqref{rk}, this decomposition is very special.  Notice that Lemma \ref{fons} inherited its name exactly from this decomposition.

{\color{Green} Moreover, for \emph{quasiconvex integrands} in \eqref{functional0} the effect of concentrations and oscillations splits additively also for the (non-linear) functional $I$; i.e.~we}  get for (non-relabeled) subsequences
\begin{align}
&\lim_{j\to\infty} \int_\O {\color{black}f}(x,\nabla u_j(x))\,\md x=\lim_{j\to\infty} \int_\O {\color{black}f}(x,\nabla z_j(x))\,\md x+\lim_{j\to\infty} \int_\O {\color{black}f}(x,\nabla w_j(x))\,\md x,
\label{split}
\end{align}
{\color{Green} with $\{u_j\}, \{w_j\}$ and $\{z_j\}$ as introduced in Lemma \ref{fons} and the discussion thereafter. Relation \eqref{split} can be proved by exploiting the} so-called $p$-Lipschitz continuity of quasiconvex functions {\color{Green} by a straightforward technical calculation (see e.g. \cite{Kroe10d})}. The $p$-Lipschitz continuity asserts that if $f:\R^{m\times n}\to \R$ is quasiconvex and $|f|\le
C(1+|\cdot|^p)$ for  some $C>0$, and $1\le p<+\infty$ then there is a  constant $\alpha\ge 0$  such that for all
$A,B\in\R^{m\times n}$
\begin{align}\label{p-lipschitz}
|f(A)-f(B)|\le\alpha (1+|A|^{p-1}+|B|^{p-1})|A-B|\ .
\end{align}

\medskip

\begin{remark}
The $p$-Lipschitz continuity holds even if $f$ is only separately convex, i.e. convex along the Cartesian axes in $\R^{m\times n}$. Various variants of this statement  are proven e.g.~in  \cite{fusco,mar} and in \cite{dacorogna}{\color{Green}; an analogous result for $k$-quasiconvex functions also holds and can be found e.g.~in  \cite{guidorzi,zappale}}. It follows from \eqref{p-lipschitz} that quasiconvex functions satisfying the mentioned bound are locally Lipschitz.   
\end{remark}

\medskip

{\color{Green} Owing to the decomposition lemma and the split \eqref{split}, we may inspect lower semicontinuity of $I$ in \eqref{functional0} along a sequence $\{u_j\}_{j \in \N}$ separately for the oscillating  and the concentrating part. Roughly speaking, the oscillating part is handled by quasiconvexity itself while additional conditions are needed for the concentrating part. This statement is formalized via the following theorem:}

\medskip

\begin{theorem}[adapted from Ka\l amajska and Kru\v{z}\'{\i}k \cite{mkak}] \label{wlsc1}
Let ${\color{black}f}\in C(\overline{\Omega}\times \R^{m\times n})$, $|{\color{black}f}|\le
C(1+|\cdot|^p)$, $C>0$, $ f(x,\cdot)$ quasiconvex for  all $x\in\overline{\Omega}$, and  $1<p<+\infty$. Then the
functional 
\begin{align}\label{kala-kru}
I(w):=\int_\O {\color{black}f}(x,\nabla w(x))\,\md x\end{align}   is sequentially weakly
lower semicontinuous on $W^{1,p}(\O;\R^m)$ if and only if for every
bounded sequence $\{w_j\}\subset W^{1,p}(\O;\R^m)$ such that
$\nabla w_j\to 0$ in measure we have
$\liminf_{j\to\infty}I(w_j)\ge I(0)$.
\end{theorem}

\medskip 

{\color{Green} Thus, let us study weak lower semicontinuity of $I$ only along purely concentrating sequences; i.e. along a sequence $\{w_j\}_{j \in \N} \subset W^{1,p}(\Omega;\R^m)$ such that $\nabla w_j \wto 0$ and the $\mathcal{L}^n(x\in \mathrm{supp}\, w_j) \to 0$ as $j\to\infty$.  For simplicity, we set $f(\cdot,0) = 0$. Then, we can write
$$
\int_\Omega f(x,\nabla w_j) \md x = \int_\Omega f^+(x,\nabla w_j) \md x - \int_\Omega f^-(x,\nabla w_j) \md x \geq \int_\Omega f(x,0) \md x - \int_\Omega f^-(x,\nabla w_j) \md x,
$$
where $f^-$ and $f^+$ are the negative and the positive part of $f$, respectively.
So, we see that lower semicontinuity of $I$ along the sequence $\{w_j\}_{j \in \N}$ is obtained if $\int_\Omega f^-(x,\nabla w_j) \md x \to 0$. Recall from Definition \ref{def-equiint} that this is always the case once the sequence $ \{f^-(\cdot, \nabla w_j)\}_j$ is equiintegrable and we conclude that for quasiconvex integrands only the fact that $\{f^-(\cdot, \nabla w_j)\}_{j \in \N}$ \emph{is a concentrating sequence might harm weak lower semicontinuity}. Notice that equiintegrability of $\{f^-(\cdot, \nabla u_j)\}_{j \in \N}$ can, for example, be achieved if the negative part of $f$ is of sub-critical growth (cf. Theorem~\ref{wlsc} below).

However, not all concentrations of $\{\vert\nabla u_j\vert^p \}_{j \in \N}$ affect the  weak lower semicontinuity of $I$. In fact, we show in Remark \ref{concentr-inter} that concentrations inside the domain $\Omega$ are ineffectual for  weak lower semicontinuity of $I$  in \eqref{kala-kru} if $f(x,\cdot)$ is quasiconvex for all $x\in\overline{\Omega}$ and $f(\cdot, A)$ is continuous for all $A\in\R^{m\times n}$. Therefore, only concentrations \emph{at the boundary} need to be excluded by further requirements since 
}  along concentrating sequences of gradients, energy may be \emph{gained} and hence the lower semicontinuity might be destroyed. {\color{Green} The following examples show that such a situation does appear:}

\medskip

\begin{example}[following \cite{KriRi10a}, \cite{BaCheMaSa13a}]
\label{concentration-exa-basic}
Choose $\Omega=(0,1)$ {\color{black} and a smooth, non-negative function $\Phi: \R \to \R$ with compact support in $(0,1)$ and such that $\int_0^1 \Phi(y) \md y = 1$. Let us now define the sequence $\{u_n\}_{n \in \N} \subset W^{1,1}(0,1)$ through 
$$
u_n(x) = 1-\int_0^x n \Phi\left(nt\right) \md t \qquad \text{ so that } \qquad u_n'(x) = - n \Phi\left(nx\right).
$$ It can bee seen that $\{u_n\}_{n \in \N}$ is a concentrating sequence that converges to 0 pointwise and in measure on $(0,1)$. Further let us choose $f(x,r,s):=s$   in \eqref{I-1D}; i.e. $f$ is a linear function and so quasiconvex. Then the functional \eqref{I-1D} fulfills $I(u_n)=-1$ for all $n$, but $u_n' \wstar 0$ in measure and $I(0)=0>-1$.}

The example illustrates the above mentioned effect that a sequence concentrating on the boundary (such as $\{u_n\}_{n \in \N}$) may actually lead to an energy gain in the limit. {\color{black} However, the failure of weak lower semicontinuity is shown with respect to the weak topology in measure for the derivative and not the weak convergence in $W^{1,1}(0,1)$. The reason is that this allows us to take a linear, and thus a particularly easy, integrand in \eqref{functional0}, which is however of critical negative growth only in $W^{1,1}(0,1)$. But any sequence converging weakly in $W^{1,1}(0,1)$ is also equiintegrable so the concentration effect could not be seen. Let us point to Example \ref{det-example} below for appropriate nonlinear integrands that lead to the same effect in $W^{1,p}(\O;\R^m)$ with $p>1$.
} 

{\color{Green} Let us also mention that the above example allows an easy adaptation to $\mathrm{BV}(0,1)$ that avoids the mollification kernel $\Phi$}. Take a sequence $\{u_n\}_{n \in \N} \subset \mathrm{BV}(0,1)$ defined through $u_n:=\chi_{(0,\frac{1}{n})}$, i.e. the characteristic function of $(0,\frac{1}{n})$ in $(0,1)$, so that $Du_n=-\delta_{\frac{1}{n}}$. Then
$$I(u)= \int_{(0,1)}\, \mathrm{d} Du(x),$$
which is a $\mathrm{BV}$-equivalent of \eqref{I-1D} with $f(x,r,s):=s$, is not weakly* lower semicontinuous on $\mathrm{BV}(0,1)$ because $I(u_n)=-1$ for all $n$, but $u_n \wstar 0$ in $\mathrm{BV}(0,1)$ and $I(0)=0>-1$. 
\end{example}

\medskip

\begin{example}[See \cite{murat}]
\label{det-example}
Let $n=m=p=2$, $0<a<1$, $\Omega:=(0,a)^2$  and  for $x\in\O$ define 
$$u_j(x_1,x_2)=\frac{1}{\sqrt{j}}(1-|x_2|)^j\big(\sin (jx_1),\cos (jx_1)\big).$$
We see that $\{u_j\}_{j\in\N}$ converges weakly in $W^{1,{\color{black}n}}(\O;\R^2)$ as well as pointwise to zero. Moreover, we calculate for $j\to\infty$
\begin{eqnarray*}
\int_{0}^a\int_0^a\det\nabla u_j(x)\,\md x\to \frac{-a}{2}< 0\ .
\end{eqnarray*} 

Hence, we see that $I(u):=\int_\O \det\nabla u(x)\,\md x$ is not weakly lower semicontinuous in $W^{1,{\color{black}n}}(\O;\R^2)$. 
This example can be generalized to arbitrary dimensions $m=n \geq 2$. Indeed, take $u\in W_0^{1,n}(B(0,1);\R^n)$ and extend $u$ by zero to the whole $\R^n$. We get that $\int_{B(0,1)}\det\nabla u(x)\,\md x=0$ because of the zero Dirichlet boundary conditions on $\partial B(0,1)$.  Take $\varrho\in\R^n$, a unit vector,  such that $\int_{D_\varrho}\det\nabla u(x)\,\md x<0$; {\color{black} here recall from Section \ref{sec-notation} that $D_\varrho:=\{x\in\R^n;\ x\cdot\varrho<0\}$} Notice that this condition can be fulfilled, if we take $u$ suitably. 

Denote $u_j(x):=u(jx)$  for all $j\in\N$; then $u_j\wto 0$ in $W^{1,n}(B(0,1);\R^n)$  but  also $\int_{D_\varrho}\det\nabla u_j(x)\,\md x\to\int_{D_\varrho}\det\nabla u(x)\,\md x<0$ by our construction.  The same conclusion can be drawn if we take $\O\subset\R^n$ with arbitrarily smooth boundary and such that 
$0\in\partial\O$.  Let $\varrho$ be the outer unit normal to $\partial\O$ at zero. Then we have for the same sequence as before 
 \begin{align*}
&\lim_{j\to\infty} \int_{\O}\det\nabla u_j(x)\,\md x= \lim_{j\to\infty} \int_{B(0,1)\cap\O}\det\nabla u_j(x)\,\md x \\
&= \lim_{j\to\infty}\int_{B(0,1)\cap\O}j^n \det\nabla u(jx)\,\md x\nonumber 
= \int_{D_\varrho} \det\nabla u(y)\,\md y<0\ .
\end{align*}
\end{example}

\medskip

\begin{remark}
\label{concentr-inter}
In this remark, we indicate why quasiconvexity is capable of preventing concentrations {\it in the domain} $\O$ from breaking weak lower semicontinuity.
Indeed, let  $\zeta\in \mathcal{D}(\O)$, $0\le\zeta\le 1$ and take a quasiconvex function ${\color{black} f}:\R^{m\times n}\to\R$ such that $|{\color{black} f}(A)|\le C(1+|A|^p)$ for some $C>0$ and all $A\in\R^{m\times n}$ with $p>1$. {\color{Green} Moreover, let $\{w_j\}_{j \in \N}$ be a purely concentrating sequence. From} 
Definition~\ref{def-quasiconvexity} for $A:=0$, {\color{Green} we have that}
\begin{eqnarray*}
\mathcal{L}^n(\O){\color{black} f}(0) \le
 \int_{\O}{\color{black} f}(\nabla (\zeta(x)w_j(x)))\,\md x
\end{eqnarray*}
{\color{Green} and by using the chain rule, the $p$-Lipschitz property \eqref{p-lipschitz} and the fact that} $w_j\to 0$ strongly in $ L^p(\O;\R^n)$ and $\{\nabla w_j\}_{k\in\N}$ is bounded in $L^p(\O;\R^{m\times n})$ {\color{Green} we get that}
  \begin{align}\label{est1}  \mathcal{L}^n(\O)f(0)\le\liminf_{j\to\infty}\int_{\O}f(\zeta(x)\nabla w_j(x))\,\md x\ .\end{align}
   Let $|\nabla w_j|^p\stackrel{*}{\wto} \sigma$  in $\mathcal{M}(\overline{\Omega})$ for a (non-relabeled) subsequence.  {\color{Green} Given the assumption that all concentrations appear inside the domain $\Omega$}, we have that  $\sigma(\partial \O)=0$, whence we continue with the following estimate
  \begin{align}\label{est2}
   & \lim_{j\to\infty} \int_{\O}{\color{black} f}(\zeta(x)\nabla w_j(x))\,\md x \nonumber\\&\qquad\le
    \lim_{j\to\infty} \int_{\O}{\color{black} f}(\nabla w_j(x))+\alpha(1-\zeta(x))
(1+\zeta^{p-1}(x))|\nabla w_j(x)|^p+ \alpha (1-\zeta(x))|\nabla
w_j(x)|\,\md x\nonumber\\
&\qquad  =\lim_{j\to\infty} \int_{\O}{\color{black} f}(\nabla w_j(x))\,\md x +\alpha\int_{\O}(1-\zeta(x))
(1+\zeta^{p-1}(x))\sigma(\md x) \ 
  \end{align}
  {\color{Green} where we again used the $p$-Lipschitz property.}
   Now, {\color{black} we choose} a sequence $\{\zeta_j\}_{j\in\N}\subset \mathcal{D}(\O)$,
   satisfying $0\le \zeta_j\le 1$ 
   that pointwise tends to  the characteristic function of $\O$, $\chi_{\O}$, $\sigma$-a.e. Taking into account (\ref{est1}) and (\ref{est2}), we have by the Lebesgue's dominated convergence theorem
   $$  \mathcal{L}^n(\O) {\color{black} f}(0)\le\lim_{j\to\infty}\int_{\O}{\color{black} f}(\nabla w_j(x))\,\md  x\ .$$
	Hence, weak lower semicontinuity is preserved.  
\end{remark}

\medskip 

This reasoning {\color{Green} of Remark \ref{concentr-inter}}, however, clearly breaks if $\partial\O$ is not a $\sigma$-null set, hence concentrations at the boundary appear. Nevertheless, {\color{Green} even} not every boundary concentration is fatal for weak lower semicontinuity.  Arguing heuristically, concentrations at $\partial\O$ are influenced by interior concentrations coming from $\O$  and exterior ones coming from the complement.  If exterior concentrations can be excluded then the interior ones  cannot spoil weak lower semicontinuity. That is, roughly speaking,  why Dirichlet boundary conditions suffice to ensure (ii) in Theorem~\ref{meyers4} at least if $k=1$. \WE If periodic boundary conditions are applicable, then they will do, as well, because exterior and interior concentrations mutually compensate due to periodicity. \EOR

The next theorem {\color{Green} formalizes the discussion concerning equiintegrability of the negative part of $f$ and Dirichlet boundary conditions}. 
 
\medskip
	
\begin{theorem}[taken from Ka\l amajska and Kru\v{z}\'{\i}k  \cite{mkak}]\label{wlsc}
Let the assumptions of Theorem~\ref{wlsc1} hold.
Let further $\{u_j\}\subset W^{1,p}(\O;\R^m) $, $u_j\rightharpoonup u$ in $W^{1,p}(\O;\R^m)$ and at least one of the following conditions be satisfied:
\begin{enumerate}[(i)]
\item for {\color{black} every} subsequence of $\{u_j\}_{j \in \N}$ (not relabeled) such that $|\nabla u_j|^p\wstar \sigma$ in $\mathcal{M}(\overline{\Omega})$,  where $\sigma\in \mathcal{M}(\overline{\Omega})$ depends on the particular subsequence, it holds that  $\sigma(\partial\O)=0$,\\
\item $\lim_{|A|\to\infty}\frac{{\color{black} f}^-(x,A)}{1+|A|^p}=0$ for all $x\in\overline{\Omega}$ where ${\color{black} f}^{-}:= {\rm max} \{ 0,-{\color{black} f}\}$,\\
\item $u_j=u$ on $\partial\O$  for every $j\in\N$ and $\O$ is
Lipschitz.
\end{enumerate}
Then  $I(u)\le \liminf_{j\to\infty}I(u_j)$.

\end{theorem}

\medskip

Notice that \emph{(ii)} is satisfied for example, if ${\color{black} f}\ge 0$ or  if
 ${\color{black} f}^-\le C(1+|\cdot|^q)$ for some $1\le q<p$ in which case  $-C(1+|A|^q)\le
{\color{black} f}(x,A)\le C(1+|A|^p)$, $C>0$ and $x\in\O$. This result can be found  e.g.~in  \cite{dacorogna}.
 
{\color{Green} It follows from the discussion in this section that} condition (ii) in Theorem \ref{meyers4} is connected with concentrations on the boundary. {\color{Green} This must have been clear already to Meyers who} conjectured \cite[p.~146]{meyers} that it can be dropped if $\partial \Omega$ is ``smooth enough" or a ``smooth enough" function is prescribed on the boundary {\it as the datum}. The second part of the conjecture turned out to be true  in the following special cases: for $k=1$ in \eqref{functional1} (see \cite[Thm.~5]{meyers} and   ~Thm.~\ref{wlsc}) or if the integrand in \eqref{functional1} depends just on the highest gradient (see end of Section \ref{sec-A-quasi}). However, the general case is still open:
\medskip
\begin{problem}
\label{Meyers-conj}
Is the functional \eqref{functional1} weakly lower semicontinuous along sequences with fixed Dirichlet boundary data if $f$ is a general function in the class $\mathcal{F}_p(\Omega)$ that is $k$-quasiconvex?
\end{problem}

\medskip

The first part of the conjecture of Meyers turned out \emph{not} to hold {\color{Green} as is illustrated by Example  \ref{det-example} where weak lower semicontinuity breaks down independently of the smoothness of $\partial \Omega$. }

Let us return to the issue of making  condition (ii) in Theorem \ref{meyers4} more explicit. It has been identified in \cite{Kroe10d}  that a suitable \emph{growth} from below of the whole functional in \eqref{functional1} (which does not necessarily imply a lower bound on the integrand ${\color{black} f}$ itself) equivalently  replaces this condition. First, let us illustrate  that some form of boundedness from below is indeed necessary for weak lower semicontinuity.

\medskip

\begin{example}
\label{growthfrombellow}
Take $u \in W^{1,p}_0(B(0,1);\R^m)$ ($1 <  p < \infty$)  and extend it by zero to the whole of $\R^n$. 
Define for $x\in\R^n$ and $j\in\mathbb{N}$  $u_j(x)=j^{\frac{n-p}{p}}u(jx)$ and  consider a smooth  domain
 $\Omega\subset \R^n$ such that $0\in\partial\Omega$; denote by $\varrho$ the outer unit normal to $\partial\Omega$ at $0$. Notice that  $u_j \rightharpoonup 0$ in $W^{1,p}(\O;\R^m)$ and $\{|\nabla u_j|^p\}_{j\in\N}$ concentrates at zero.  Moreover, take a  function  ${\color{black} f}:\R^{m\times n}\to\R$ that is positively $p$-homogeneous, i.e., ${\color{black} f}(\alpha \xi)=\alpha^p {\color{black} f}(\xi) $ for all $\alpha\ge 0$. If 
$$
I(u)= \int_\Omega {\color{black} f}(\nabla u(x)) \, \md x
$$ 
is weakly lower semicontinuous on $W^{1,p}(\O;\R^m)$  then 
\begin{eqnarray}\label{eq:intro}
0&=&I(0)\le \liminf_{j\to\infty} \int_\Omega {\color{black} f}(\nabla u_j(x))\,\md x= \liminf_{j\to\infty} \int_{B(0,1/j)\cap\Omega}{\color{black} f}(\nabla u_j(x))\,\md x\\ &=& \liminf_{j\to\infty}\int_{B(0,1/j)\cap\Omega}j^n {\color{black} f}(\nabla u(jx))\,\md x= \int_{D_\varrho} {\color{black} f}(\nabla u(y))\,\md y\ .\nonumber
\end{eqnarray}
Thus, we see that 
\begin{equation}
0\le \int_{D_\varrho} {\color{black} f}(\nabla u(y))\,\md y
\label{example}
\end{equation}
for all $u \in W^{1,p}_0(B(0,1);\R^m)$  forms a necessary condition for weak lower semicontinuity of $I$ whenever ${\color{black} f}$ is positively $p$-homogeneous.  
\end{example}

\medskip

For functions that are not $p$-homogeneous, S. Kr\"{o}mer \cite{Kroe10d}  generalized \eqref{example} as follows.

\medskip 

\begin{definition}[following \cite{Kroe10d}\footnote{In \cite{Kroe10d} this condition is actually not referred to as $p$-quasi-subcritical growth from below but is introduced in Theorem 1.6 (ii).}]
Assume that $\O\subset\R^n$ has a smooth boundary and let $\varrho(x)$ be the unit outer normal to $\partial\O$ at $x$.
We say that a function ${\color{black} f}: \Omega \times \R^{m \times n} \to \R$ is of \emph{$p$-quasi-subcritical growth from below} if for every $x\in\partial\O$ and for every $\varepsilon>0$, there exists $C_\varepsilon\geq 0$ such that
\begin{align}\label{r-qsb}
	&\int_{D_{\varrho(x)}(x,1)} {\color{black} f}(x,\nabla u(z))\md z\geq -\varepsilon \int_{D_{\varrho(x)}(x,1)} |\nabla u(z)|^p\md z-C_\varepsilon~~
	\text{for all $u \in W^{1,p}_0(B(0,1);\R^m)$}.
\end{align}
\end{definition}

\medskip

It has been proved in \cite{Kroe10d} that the $p$-quasi-subcritical growth from below of the function ${\color{black} f}:= {\color{black} f}(x,\nabla u)$ equivalently  replaces (ii) in Theorem~\ref{meyers4}.

Notice that \eqref{r-qsb} is expressed only in terms of ${\color{black} f}$ and that it is local in $x$. Moreover, it shows again that, at least in the case when ${\color{black} f}$ does depend only on the first gradient of $u$ but not on $u$ itself, only \emph{concentrations at the boundary} may interfere with weak lower semicontinuity of functionals involving quasiconvex functions. 

\medskip

\begin{remark}
Let us realize that  \eqref{r-qsb} implies \eqref{example} if ${\color{black} f}$ is positively $p$-homogeneous and independent of $x$. To this end, we use, for $t \geq 0$, $u = t \tilde{u}$ in \eqref{r-qsb} to see that
$$
0 \leq \frac{1}{t^p}\left(\int_{D_\varrho(x_0,1)} {\color{black} f}(t\nabla \tilde{u}(x))\,dx + \varepsilon  |t \nabla \tilde{u}(x)|^p\md x +C_\varepsilon \right).
$$
Letting now $t \to \infty$ gives that $C_\varepsilon = 0$. Then, we may also send $\varepsilon \to 0$ to get \eqref{example}.
\end{remark}

\medskip

Since only concentration effects play a role for (ii) in Theorem~\ref{meyers4}, it is natural to expect that weak lower semicontinuity can be linked to properties of the so-called \emph{recession function} of the function ${\color{black} f}$, if it admits one. Recall, that we say that the functions 
${\color{black} f}_\infty:\overline{\O}\times\R^{m\times n}\to\R$  is a recession function for ${\color{black} f}: \overline{\O}\times\R^{m\times n}\to\R$ if for all $x\in\overline{\O}$
$$
\lim_{|A|\to\infty} \frac{{\color{black} f}(x,A)-{\color{black} f}_\infty(x,A)}{|A|^p}=0.
$$
Thus, informally speaking, the recession function describes the behavior of $f$ at ``infinitely large matrices''. Note that ${\color{black} f}_\infty$ is necessarily positively $p$-homogeneous; i.e. ${\color{black} f}_\infty(x,\lambda A)=\lambda^p {\color{black} f}_\infty(x, A)$ for all $\lambda\ge 0$, all $x\in\overline{\Omega}$, and all $A\in\R^{m\times n}$.

It follows from Remark 3.9 in \cite{Kroe10d} that if ${\color{black} f}$ admits a recession function, then quasi-subcritical growth from below is equivalent to \eqref{example} for ${\color{black} f}_\infty$.

Since weak lower semicontinuity is connected to quasiconvexity and to condition (ii) in Theorem \ref{meyers4} which is connected to effects at the boundary, it is reasonable to ask whether the two ingredients can be combined. Indeed, so-called {\it quasiconvexity at the boundary} was introduced {\color{black} by Ball and Marsden} \cite{bama} to study necessary conditions satisfied by local minimizers of variational problems -- we also refer to \cite{grabovsky-mengesha1,grabovsky-mengesha2,MiSp98,sprenger,silhavy} where this condition is analyzed, too. 
In order to define quasiconvexity at the boundary, we put for $1\le p\le+\infty$ 
\be\label{testspace}
W^{1,p}_{{\color{black} \partial D_\varrho\setminus\Gamma_\varrho}}(D_\varrho;\R^m):=\{u\in W^{1,p}(D_\varrho;\R^m);\ u=0 \mbox{ on } \partial D_\varrho\setminus\Gamma_\varrho\}\  ,\ee
where $\Gamma_\varrho$ is the planar part of $\partial D_\varrho$.

\medskip

\begin{definition}[taken from \cite{MiSp98}\footnote{The original definition in \cite{bama} considers the case $q:=0$.}]\label{qcb-def}
Let $\varrho\in\R^n$ be a unit vector. A function $f:\R^{m\times n}\to\R$ is called quasiconvex at the boundary at the point $A\in\R^{m\times n}$ with respect to $\varrho$ if there is $q\in\R^m$ such that for all $\varphi\in W^{1,\infty}_{{\color{black} \partial D_\varrho\setminus\Gamma_\varrho}}(D_\varrho;\R^m)$ it holds 
\be\label{qcbinfty}
\int_{\Gamma_\varrho} q\cdot \varphi(x)\,\md S + f(A)\mathcal{L}^n(D_\varrho)\le \int_{D_\varrho} f(A+\nabla \varphi(x))\,\md x\ . \ee
\end{definition}

\medskip

Let us remark that, analogously to quasiconvexity, we may generalize quasiconvexity at the boundary to $W^{1,p}$-quasiconvexity at the boundary (for $1 < p < \infty$) by using all $\varphi\in W^{1,p}_{{\color{black} \partial D_\varrho\setminus\Gamma_\varrho}}(D_\varrho;\R^m)$ as test functions in \eqref{qcbinfty}. For functions with $p$-growth these two notions coincide.

\medskip

\begin{remark}
Let us give an intuition on the above definition. Take a convex function $f: \R^{m \times n} \to \R$ and $\varphi\in W^{1,\infty}_{{\color{black} \partial D_\varrho\setminus\Gamma_\varrho}}(D_\varrho;\R^m)$. Then we know that
$$
f(A + \nabla \varphi(x)) \geq  f(A) + {\color{Green} g(A)}\cdot \nabla \varphi(x),
$$
{\color{Green} where $g(A)$ is a subgradient    of $f$ evaluated at $A$; see, e.g., Rockafellar-Wets \cite{rockafellar} for details about this notion}. Integrating this expression over $\Omega$ then gives
$$
\int_\Omega f(A + \nabla \varphi(x)) \md x \geq  \int_\Omega \big(f(A) + {\color{Green} g(A)}\cdot \nabla \varphi \big) \md x =  \mathcal{L}^n(\O) f(A) + \int_{\partial \Omega} \big({\color{Green} g(A)} \varrho \big) \cdot \varphi\, \md S,
$$
where $\varrho$ is the outer normal to $\partial \Omega$. Now when setting $q: = {\color{Green} g(A)\varrho}$ we obtain the definition of the quasiconvexity at the boundary. 
\end{remark}

\medskip



 
\begin{remark}
It is possible to work with more general domains than half-balls in Definition \ref{qcb-def}; namely with so-called \emph{standard boundary domains}. We say that 
$\tilde{D}_\varrho$ is a standard boundary domain with the normal $\varrho$ if there is $a\in\R^n$ such that $\tilde{D}_\varrho\subset H_{a,\varrho}:=\{x\in\R^n;\ \varrho\cdot x<a\}$ and the $(n-1)$--dimensional interior of 
$\partial \tilde{D}_\varrho\cap\partial H_{a,\varrho}$, called $\Gamma_\varrho$, is nonempty. 
Roughly speaking, this means that the boundary of $\tilde{D}_\varrho$ should contain a planar part.

As with standard quasiconvexity, if (\ref{qcbinfty}) holds for one standard boundary domain it holds for  other standard boundary domains, too.
\end{remark}

\medskip

\begin{remark}\label{remark1}
If $p>1$,  ${\color{black} f}:\R^{m\times n}\to\R$ is positively $p$-homogeneous, continuous, and $W^{1,p}$-quasiconvex at the boundary at $(0,\varrho)$ then $q=0$ in (\ref{qcbinfty}). Indeed, we have ${\color{black} f}(0)=0$ and
suppose, by contradiction, that   $\int_{D_\varrho} f(\nabla \varphi(x))\,\md x<0$ for some $\varphi\in W^{1,\infty}_{{\color{black} \partial D_\varrho\setminus\Gamma_\varrho}}(D_\varrho;\R^m)$. 
By (\ref{qcbinfty}), we  must have for all $\lambda>0$
$$ 0\le \lambda^p \int_{D_\varrho} {\color{black} f}(\nabla \varphi(x))\,\md x- \lambda \int_{\Gamma_\varrho} q\cdot \varphi(x)\,\md S\ .$$
However, this is not possible for $\lambda>0$ large enough and therefore for all $\varphi\in W^{1,\infty}_{{\color{black} \partial D_\varrho\setminus\Gamma_\varrho}}(D_\varrho;\R^m)$ it has to hold that 
 $\int_{D_\varrho} {\color{black} f}(\nabla \varphi(x))\,\md x\ge0$. Thus, we can   take  $q=0$.  
\end{remark}

\medskip

From the above remark and from \eqref{example}, we have the following lemma:

\medskip

\begin{lemma}
If a function ${\color{black} f}: \R^{m \times n} \to \R$ is $W^{1,p}$-quasiconvex at the boundary at zero and every  $\varrho\in\R^n$, a unit normal vector to $\partial\O$, then it is also of $p$-subcritical growth from below. The two notions become  equivalent if $f$ is  positively $p$-homogeneous. Here $\O$ must have a smooth boundary, so that the outer unit normal to it is defined everywhere.
\end{lemma}

\medskip

{\color{black} All the results, we presented so far concern just the case $k=1$ and  integrands $f=f(x,\nabla u)$ in \eqref{functional0}. In fact, in the general case in which $f=f(x,u,\nabla u)$ only a few results are available. One of them is, of course, Meyers' original Theorem \ref{meyers4} that applies to a general class of integrands. Another result is due to} Ball and Zhang \cite{ball-zhang} who considered  the following bound on a Carath\'{e}odory integrand $f$:
\begin{align}\label{af-growth'}
|{\color{black} f}(x,s,A)|\le a(x)+ C(|s|^p+|A|^p)\ ,
\end{align}

where $C>0$ and $a\in L^1(\O)$. Under \eqref{af-growth'}, we cannot expect weak lower semicontinuity of $I$ along generic sequences. Indeed, they proved the following weaker result.

\medskip

\begin{theorem}[Ball and Zhang \cite{ball-zhang}]\label{thm-bz}
Let $1\le p<+\infty$,  $u_k\rightharpoonup u$ in $W^{1,p}(\O;\R^m)$, ${\color{black} f}(x,s,\cdot)$ be quasiconvex for all $s\in\R^m$ and almost all $x\in\O$, and let \eqref{af-growth'} hold. 
Then there exists a sequence of sets $\{\O_j\}_{j\in\N}\subset\O$ satisfying $\O_{j+1}\subseteq\O_j$ for all $j\ge 1$, and $\lim_{j\to\infty}\mathcal{L}^n(\O_j)=0$ such that for all $j\ge 1$ 
 \begin{align}
\int_{\O\setminus\O_j} {\color{black} f}(x,u(x),\nabla u(x))\,\md x\le \liminf_{k\to\infty} \int_{\O\setminus\O_j} {\color{black} f}(x,u_k(x),\nabla u_k(x))\,\md x\ .
\end{align}
\end{theorem}

\medskip
The sets $\{\O_j\}$ that must be removed (or {\color{black} bitten off}) from $\O$ are  sets where possible concentration effects of the bounded sequence 
$\{ \vert {\color{black} f}(x,u_k,\nabla u_k)\vert \}_{k\in\N}\subset L^1(\O)$  take place. Thus, $\{\O_j\}$ depends on the sequence $\{u_k\}$ itself and  $\O_j$ are not known a-priori. \SECOND Nevertheless, in fact, $\O_j$ depends just on the sequence of gradients. Indeed, \eqref{af-growth'} and the strong convergence of $\{u_k\}_{k\in\N}$ in $L^p(\O;\R^m)$ imply that whenever $\{|\nabla u_k|^p\}_{k\in\N}$  is equiintegrable then the same holds for $\{\vert f(x,u_k(x),\nabla u_k(x))\vert \}_{k\in\N}$.  The main tool  of the proof of Theorem \ref{thm-bz} is the Biting Lemma due to Chacon \cite{brooks-chacon,ball-murat}. \EOR

\medskip

\begin{lemma}[Biting lemma]
Let $\O\subset\R^n$ be a bounded measurable set. Let $\{z_k\}_{k \in \N}\subset L^1(\O;\R^m)$ be bounded. Then there is a (non-relabeled) subsequence of $\{z_k\}_{k \in \N}$ ,  $z\in L^1(\O;\R^m)$ and a  sequence of sets $\{\Omega_j\}_{j\in\N}\subset\O$, $\O_{j+1}\subset\O_j$, $j\in\N$,  with $\mathcal{L}^n(\O_j)\to 0$ for $j\to\infty$ such that $z_k\wto z$ in $L^1(\O\setminus\O_j;\R^m)$ for $k\to\infty$ and every $j\in\N$.
\end{lemma}   

\medskip

 {\color{Green}
 
 Finally, let us remark that concentration effects do not appear if we study  lower semicontinuity of functionals with linear growth with respect to the weak $W^{k,1}(\Omega, \R^m)$-topology (see Remark \ref{rem-Wk1}). Nevertheless, this topology is too strong when it comes to the study of existence of minimizers for such functionals, cf. the discussion at the end of Section \ref{sec-BoundedFromBelow}.

\medskip
 
\begin{remark}[case $p=1$]
\label{rem-Wk1}
Let us remark that if examining weak lower semicontinuity of intergal functionals with linear growth along sequences  converging weakly in  $W^{k,1}(\Omega, \R^m)$ condition (ii) in Theorem \ref{meyers4} is also satisfied automatically. This follows from the fact that such sequences are already equiintegrable.
\end{remark}
}

\subsection{Integrands bounded from below}
\label{sec-BoundedFromBelow}
 
{\color{black} In the previous section, we saw that characterizing weak lower-semicontinuity of integral functionals with the integrand unbounded from below brings along many peculiarities if the negative part of the integrand is not equiintegrable. Naturally, all difficulties disappear if the integrand is bounded from below; notice, for example, that} condition (ii) in Theorem \ref{meyers4} is automatically satisfied. {\color{black} Thus, all the results from the previous section are readily applicable in this situation, too. Yet, as the case $f \geq 0$ for an integrand in \eqref{functional1} is the most typical one found in applications, it is worth  studying it independently. In fact, it is natural to expect that if $f$ in \eqref{functional1} has a lower bound, one can strengthen Theorem \ref{meyers4} by relaxing the }continuity assumptions stated in Definition \ref{class-F_R}. {\color{black} We review the available results in this section.}  
 
In case $k=1$ in \eqref{functional1},  the following result due to E.~Acerbi and N.~Fusco \cite{af} shows that {\color{Green} the continuity assumption on the integrand can be replaced by the Carath\'{e}odory property.}

\medskip 

\begin{theorem}[Acerbi and Fusco \cite{af}]\label{thm-af}
Let  $k=1$, $\O\subset\R^n$ be an open, bounded set, and let $f:\O\times\R^m\times\R^{m\times n}\to [0;+\infty)$ be a Carath\'{e}odory integrand, i.e., $f(\cdot,s, A)$ is measurable for all $(s,A) \in \R^m\times\R^{m\times n}$ and $f(x,\cdot,\cdot)$ is continuous for almost all $x\in\O$.  
Let further $f(x,s,\cdot)$ be quasiconvex for almost all $x\in\O$ and all $s\in\R^m$, and suppose that for some $C>0$, $1 \leq p <+\infty$, and $a\in L^1(\O)$ we have that\footnote{ This bound is often called ``natural growth conditions''.}
\begin{align}\label{af-growth}
0\le f(x,s,A)\le a(x)+ C(|s|^p+|A|^p)\ .
\end{align} 
Then $I:W^{1,p}(\O;\R^m)\to[0;+\infty)$  given in \eqref{functional1} is weakly lower semicontinuous on $W^{1,p}(\O;\R^m)$.  
\end{theorem}

\medskip

Interestingly, the paper by Acerbi and Fusco \cite{af} already implicitly contains a version of the decomposition lemma~\ref{fons}. 

Marcellini \cite{mar} proved, by a different technique of constructing a suitable non-decreasing sequence of approximations, a very similar result to Theorem \ref{thm-af} allowing also for a slightly more general growth 
\begin{equation}
-c_1|A|^r - c_2|s|^t-c_3(x) \leq f(x,s,A) \leq g(x,s)\big( 1+|A|^p\big),
\label{growth-mar}
\end{equation}
where $c_1, c_2 \geq 0$, $c_3 \in L^1(\Omega)$; $g$ is an arbitrary Carath\'{e}odory function and the exponents satisfy that $p \geq 1$, $1 \leq r < p$ (but $r=1$ if $p=1$) and $1 \leq t < np/(n-p)$ if $p < n$ and otherwise $t \geq 1$.

Note that the growth condition \eqref{growth-mar} actually allows for integrands unbounded from below but the exponent $r$ determining this growth is strictly smaller than $p$. Such integrands are of \emph{sub-critical growth} and for integrands of the class $\mathcal{F}_p(\Omega)$ weak lower semicontinuity under this growth follows also from Theorem \ref{wlsc}(ii).

Acerbi and Fusco \cite[p.~127]{af} remarked  that ``...{\it using more complicated notations as in \cite{ball-curie-olver}, \cite{meyers}, our results can be extended to the case of functionals of the type \eqref{functional1}}''.This extension has been considered by Fusco \cite{fusco} for the case $p=1$ an later by  Guidorzi and Poggilioni \cite{guidorzi} who rewrote functional \eqref{functional1} as (using the notation from Section \ref{sec-notation})
\begin{equation}
I(u)=\int_\Omega f(x, \nabla^{[k-1]} u(x), \nabla^k u(x) ) \md x
\label{gui-pogg}
\end{equation}
and proved the following.

\medskip

\begin{proposition}[{\color{black} Guidorzi and Poggilioni }\cite{guidorzi}]
Let $f:\Omega \times Y(n,m,k-1) \times X(n,m,k) \to \R$ be a Carath\'{e}odory $k$-quasiconvex function satisfying for all $H \in Y(n,m,k-1)$ and all $A \in X(n,m,k)$
\begin{align*}
0 \leq f(x,H,A) &\leq g(x, H)(1 + |A|)^p \\
|f(x, H, A)-f(x, H, B)| &\leq C(1 + |A|^{p-1} + |B|^{p-1})|A-B|
\end{align*}
where $g$ is a Carath\'{e}odory function and $C \geq 0$. Then the functional from \eqref{gui-pogg} is weakly lower semicontinuous in $W^{k,p}(\Omega; \R^n)$ for $1 \leq p < \infty$ and $k \in \N$.
\end{proposition}

\medskip

Note that in this result the continuity of the integrand in the space variable $x$ could be omitted, which is, roughly speaking, due to the fact that quasiconvexity is enough to handle the concentration effects. On the other hand, the continuity assumption from Definition \ref{class-F_R}(ii) still remains present (with $\gamma = 1$). A similar result can be drawn from the more general setting of $\mathcal{A}$-quasiconvexity (which we review in Section \ref{sec-A-quasi} below) considered in \cite{braides-fonsea-leoni}.

{\color{black} Let us end this section with some remarks on weak lower semicontinuity of integral functionals on $W^{k,1}(\Omega; \R^m)$.} While the above results handle also weak lower semicontinuity on $W^{k,1}(\Omega; \R^m)$ with respect to the standard weak convergence in this space, it is more suitable to investigate lower semicontinuity with respect to the strong convergence in $W^{k-1,1}(\Omega; \R^m)$. This is due to the fact that $W^{k,1}(\Omega; \R^m)$ is not reflexive and therefore  coercivity of \eqref{functional1} does not allow us to select a minimizing sequence that would be weakly convergent in $W^{k,1}(\Omega; \R^m)$ but the strong convergence in $W^{k-1,1}(\Omega; \R^m)$ can be assured.

The case for $k=1$ was treated by Fonseca and M\"uller \cite{fonseca-mueller} who considered continuous integrands under mild growth conditions. The result was later generalized by Fonseca, Leoni, Mal\'y, and Paroni \cite{fonseca-leoni-maly-paroni} not only with respect to the continuity of the integrand that could be partially dropped, but also to arbitrary $k$. We give the result {\color{black} for $k=1$ in Theorem \ref{thm-fonseca1} while the general case is given} in Theorem \ref{thm-fonseca}.

\medskip
{\color{black}
\begin{theorem}[due to Fonseca, Leoni, Paroni and Mal\'{y}  \cite{fonseca-leoni-maly-paroni}; case $k=1$] 
\label{thm-fonseca1}
Let $f$ in \eqref{functional0} be a Borel integrand that is moreover continuous in the following sense: For all $\varepsilon > 0$ and $(x_0, s_0) \in \Omega \times \R^m$ there exist $\delta > 0$ and a modulus of continuity $\omega$ with the property that, for some $C > 0$, $ \omega(t)  \leq C(1+t), t>0$ such that
$$
f(x_0, s_0, A) - f(x, s, A) \leq \varepsilon(1+f(x, s, A)) + \omega(|s_0-s|),
$$
for all $x \in \Omega$ satisfying $|x-x_0| \leq \delta$ and for all $s \in \R^m$ and all $A \in \R^{m \times n}.$ Suppose further that $f$ is quasiconvex and satisfies
$$
0 \leq f(x_0, s, A) \leq c (1+|A|) \qquad \forall A\in \R^{m \times n} 
$$
for some $c > 0$ or that $f$ is convex in the last variable. Then, \eqref{functional0} is lower semicontinuous with respect to the strong convergence in $L^1(\Omega; \R^m)$.
\end{theorem}

}
\medskip 

\begin{theorem}[due to Fonseca, Leoni, Paroni and Mal\'{y} \cite{fonseca-leoni-maly-paroni}] 
\label{thm-fonseca}
Let {\color{black}$f$} in \eqref{functional1} be a Borel integrand that is moreover continuous in the following sense: For all $\varepsilon > 0$ and $(x_0, H_0) \in \Omega \times Y(n,m,k-1)$ there exist $\delta > 0$ and a modulus of continuity $\omega$ with the property that, for some $C > 0$, $ \omega(s)  \leq C(1+s), s>0$ such that
$$
f(x_0, H_0, A) - f(x, H, A) \leq \varepsilon(1+f(x, H, A)) + \omega(|H_0-H|),
$$
for all $x \in \Omega$ satisfying $|x-x_0| \leq \delta$ and for all $H \in Y(n,m,k-1)$ and all $A \in X(n,m,k).$ Suppose further that $f$ is $k$-quasiconvex and satisfies
$$
\frac{1}{c}|A| - c \leq f(x_0, H_0,A) \leq c (1+|A|),
$$
for some $c > 0$ and all $A \in X(n,m,k)$. 

Then, \eqref{functional1} is lower semicontinuous with respect to the strong convergence in $W^{k-1,1}(\Omega; \R^m)$.
\end{theorem}

\medskip

For the functions $f: X(m,n,k) \to \R$, i.e. those depending only on the highest gradient, an analogous result has been obtained in \cite{cicco}. {\color{Green} We point the reader also to the suggested further reading on integrals with linear growth in Section \ref{sec:reading}.}

 \section{Null Lagrangians}
 \label{sec-NullLag}
{\color{black} After having studied weak lower semicontinuity, let us turn our attention to conditions under which the functional \eqref{functional1} is actually \emph{weakly continuous} on $W^{k,p}(\Omega; \R^m)$. As it will turn out, \eqref{functional1} is weakly continuous only for a small, special class of integrands $f$, the so-called \emph{null-Lagrangians} (cf. Theorem \ref{char-HONL} below). Null-Lagrangians are known explicitly and consist of, roughly speaking, minors of the highest order gradient; we review their characterization in this section. Null-Lagrangians play an important role in the calculus of variations, notably they are at the heart of the definition of polyconvexity that is sufficient for weak lower semicontinuity (cf. Section \ref{sec-polyconvexity} for more details).}

We start the discussion by presenting definitions of null Lagrangians of the first and higher order. 
 
 \medskip
 
\begin{definition}
\label{def-null-1st}
 We say that a continuous map $L:\R^{m\times n}\to\R$ is a  null Lagrangian of the first order, if for every $u\in C^1(\overline{\O};\R^m)$ and every  $\varphi\in C^1_0(\O;\R^m)$ it holds 
that 
\begin{align}\label{null-lagr-def}
\int_\O L(\nabla u(x)+\nabla \varphi(x))\,\md x= \int_\O L(\nabla u(x))\,\md x \ .
\end{align} 
\end{definition}

\medskip

Notice that the definition is independent of the particular Lipschitz domain $\O$. In fact, if \eqref{null-lagr-def} holds for one domain $\Omega$ it also holds for all other (Lipschitz) domains.

\medskip

\begin{remark}
The name ``null Lagrangians'' comes from the fact that,  
if $L$ is even smooth so that the variations of $J(u):=\int_\O L(\nabla u(x))\,\md x$ can be evaluated, it easily follows from \eqref{null-lagr-def}  that $J$ satisfies  $J'(u)=0$ for all  $u\in C^1(\overline{\Omega};\R^m)$. In other words, the Euler-Lagrange equations of $J$ are fulfilled identically in the sense of distributions.
\end{remark}

\medskip

\begin{remark}
Let us notice that, if $L$ is a null Lagrangian, the value of $J(u)=\int_\O L(\nabla u(x))\,\md x$ is only dependent on the boundary values of $u$. This can be seen from \eqref{null-lagr-def} as the value remains unchanged even if we add arbitrary functions vanishing on the boundary.
\end{remark}

\medskip

It is straightforward to generalize \eqref{def-null-1st} also to higher order problems.

\medskip

\begin{definition}
Let $k\ge 2$. We say that $L:X(n,m,k)\to \R$ is a (higher-order) null Lagrangian if 
\begin{align}
\int_\O L(\nabla^{k}u(x)+ \nabla^{k}\varphi(x))\,\md x= \int_\O L(\nabla^{k}u(x))\,\md x
\end{align}
for all $u\in C^k(\overline{\O};\R^m)$ and all $\varphi\in C^k_0(\O;\R^m)$. 
\end{definition}

\medskip

Similarly as in the first-order gradient case, the definition is independent of the particular (Lipschitz) domain $\O$.
In the same way as in the first order case,{\color{Green} given sufficient smoothness,} it follows that  Euler-Lagrange equations
\begin{align}\label{E-L-higher}
\sum_{|K|\le l}(-D)^K\frac{\partial L}{\partial u^i_I}(\nabla ^l u)=0
\end{align}
are satisfied identically in the sense of distributions for arbitrary $u\in C^k(\overline{\Omega};\R^m)$.

\medskip

\begin{remark}
\label{rem-genNull}
It is natural to generalize the notion of null Lagrangians to functionals of the type \eqref{functional1}, i.e., those depending also on lower order gradients, in the following way: We say that the function $L: \Omega \times Y(n,m,k) \to \R$ is a null Lagrangian for the functional \eqref{functional1} if for all $u \in C^k(\Omega; \R^m)$ and all $\varphi\in C^k_0(\O;\R^m)$ it holds that
$$
J(u+\varphi) = J (u) \qquad \text{ and } \qquad J(u) = \int_\O L(x,u(x),\nabla u(x),\ldots, \nabla^{k}u(x))\,\md x.
$$
We shall see in the end of the section that null Lagrangians for these types of functionals are actually determined by null Lagrangians at least if $k=1$.
\end{remark}

\medskip

The following result {\color{Green} highlights some of the remarkable properties of null Lagrangians $L$ of first and higher order.} In particular, it shows that null Lagrangians are the only integrands  for  which  $u\mapsto\int_\O L(\nabla^{k}u(x))\,\md x$ is  continuous in the weak topology of suitable Sobolev spaces.  It is due to {\color{black} Ball, Curie, and Olver} \cite{ball-curie-olver}.

\medskip 

\begin{theorem}[Characterization of (higher-order) null Lagrangians] 
\label{char-HONL}
Let $L:X(n,m,k)\to\R$ be continuous. Then the following statements are mutually equivalent:
\begin{itemize}
\item[(i)] $L$ is a null Lagrangian,
\item[(ii)] $\int_\O L(A+\nabla^{k}\varphi(x))\,\md x=\int_\O L(A)\,\md x$ for every $\varphi\in C^\infty_0(\O;\R^m)$ and every $A\in X(n,m,k)$ and every open subset $\O\subset\R^n$,
\item[(iii)] $L$ is continuously differentiable and \eqref{E-L-higher} holds in the sense of distributions,
\item[(iv)] The map $u\mapsto L(\nabla^k u)$ is sequentially weakly* continuous from $W^{k,\infty}(\O;\R^m)$ to $L^\infty(\O)$. This means that if $u_j\wstar u$ in $W^{k,\infty}(\O;\R^m)$ then $L(\nabla^ku_j)\wstar L(\nabla^k u)$  in $L^\infty(\O)$,
\item[(v)] $L$ is a polynomial of degree $p$ and the map $u\mapsto L(\nabla^k u)$ is sequentially weakly* continuous from $W^{k,p}(\O;\R^m)$ to $\mathcal{D}'(\O)$. This means that if $u_j\rightharpoonup u$ in $W^{k,p}(\O;\R^m)$ then $L(\nabla^ku_j)\rightharpoonup L(\nabla^k u)$ in $\mathcal{D}'(\O)$.
\end{itemize}
\end{theorem}

\medskip

While Theorem~\ref{char-HONL} provides us with very useful properties of null Lagrangians it is interesting to note that they are known \emph{explicitly} in the first as well as in the higher order. In fact, null Lagrangians are formed by minors or sub-determinants of the gradient entering the integrand in $J$. 

\subsection{Explicit characterization of null Lagrangians of the first order}

Let us start with the first order case: If $A\in\R^{m\times n}$ we denote by  $\mathbb{T}_i(A)$ the vector of all subdeterminants of $A$ of order $i$ for $1\le i\le\min(m,n)$.
Notice that the dimension of $\mathbb{T}_i(A)$  is $d(i):={m\choose i}{n\choose i}$, hence the number of all subdeterminants of $A$ is $\sigma:={m+n\choose n}-1$. Finally, we write $\mathbb{T}:=(\mathbb{T}_1,\ldots,\mathbb{T}_{\min(m,n)})$. For example, if $m=1$ or $n=1$ then $\mathbb{T}(A)$ consists only of entries of $A$, if $m=n=2$ then $\mathbb{T}(A)=(A,\det A)$  and for $m=n=3$ we obtain $\mathbb{T}(A)=(A,\cof A, \det A)$. 

Clearly, linear maps are  weakly continuous. Yet, it has been known at least since \cite{morrey-orig,reshetnyak-0,ball77} that also minors have this property (see Theorem \ref{cof} below). This result, usually called {\it (sequential) weak continuity of minors}, is unexpected because if $i>1$ then $A\mapsto\mathbb{T}_i(A)$ is a nonlinear polynomial of the $i$-th degree. As it is well-known, weak convergence generically does not commute with nonlinear mappings.

\medskip

\begin{theorem}[Weak continuity of minors (see e.g.~\cite{dacorogna})]
\label{cof}
Let $\O\subset\R^{n}$ be a bounded  Lipschitz domain. Let $1\le i\le\min(m,n)$.  Let $\{u_k\}_{k\in\N}\subset W^{1,p}(\O;\R^m)$ be such that $u_k\wto u$ in $ W^{1,p}(\O;\R^m)$ for $p>i$.
Then $\mathbb{T}_i(\nabla u_k)\wto \mathbb{T}_i(\nabla u)$  in $L^{p/i}(\O;\R^{d(i)})$.
\end{theorem}

\medskip

\SECOND
The proof of Theorem~\ref{cof} uses the structure of null Lagrangians, namely that they can be written in the divergence form. To explain briefly this idea, we just restrict ourselves to $m=n=2$.
We have for $u\in C^2(\O;\R^2)$
\begin{align}\label{det=Det}
\det\nabla u & =\frac{\partial u_1}{\partial x_1}\frac{\partial u_2}{\partial x_2}-\frac{\partial u_1}{\partial x_2}\frac{\partial u_2}{\partial x_1}  = \frac{\partial}{\partial x_1}\Big(u_1\frac{\partial u_2}{\partial x_2}\Big)-\frac{\partial}{\partial x_2}\Big(u_1\frac{\partial u_2}{\partial x_1}\Big)
\end{align}
Hence, if $\varphi\in \mathcal{D}(\O)$ is arbitrary we get  
\begin{align}\label{determinant-passage}\int_\O\det\nabla u(x)\varphi(x)\,\md x=-\int_\O\Big(u_1\frac{\partial u_2}{\partial x_2}\Big)\frac{\partial\varphi(x)}{\partial x_1}-\Big(u_1\frac{\partial u_2}{\partial x_1}\Big)\frac{\partial\varphi(x)}{\partial x_2}\ .\end{align}
If $u_k\wto u$ in $W^{1,p}(\O;\R^m)$ for $p>2$ then the right-hand side of \eqref{determinant-passage} written for $u_k$ in the place of $u$ allows us easily to  pass to the limit for $k\to\infty$ to obtain  Theorem~\ref{cof} for $m=n=i=2$.
Notice that the right-hand side of \eqref{determinant-passage} is defined in the sense of distributions even if $p\ge 4/3$, however, the integral identity \eqref{determinant-passage} fails to hold if $p<2$. Inspired by a conjecture of Ball \cite{ball77}, M\"{u}ller \cite{mueller-Det} showed that if $u\in W^{1,p}(\O;\R^2)$,  $p\ge 4/3$, then the {\it distributional determinant} 
$$\mathrm{Det}\nabla u:=\frac{\partial}{\partial x_1}\big(u_1\frac{\partial u_2}{\partial x_2}\big)-\frac{\partial}{\partial x_2}\big(u_1\frac{\partial u_2}{\partial x_1}\big)$$
belongs to $L^1(\O)$ and $\det\nabla u= \mathrm{Det}\nabla u$. Generalizations to higher dimensions are possible, defining the distributional determinant with the help of the cofactor matrix.  We refer to \cite{mueller-Det} for details.

\EOR

{\color{Green} Minors are the \emph{only} mappings depending exclusively on $\nabla u$ which are weakly continuous; and  thus in view of Theorem \ref{char-HONL1} the only null Lagrangians of the first order. We make the statement more precise in the following theorem.

\medskip

\begin{theorem}[See \cite{ball-curie-olver} or \cite{dacorogna}]\label{char-HONL1}
Let $L\in C(\R^{m \times n})$. Then $L$ is a null Lagrangian if and only if it is an affine combination of elements of $\mathbb{T}$, i.e., 
for every $A\in\R^{m\times n}$ 
\begin{align}\label{polyaff}
L(A)= c_0+c\cdot \mathbb{T}(A)\ ,\end{align}
where $c_0\in\R$ and $c\in\R^\sigma$ are arbitrary constants.
\end{theorem}}

\medskip

Let us note however, that it has been realized independently in e.g. \cite{edelen,ericksen} that minors are the only maps for which the Euler-Lagrange equation of $J(u) = \int_\Omega L(\nabla u) \md x$ is satisfied identically. 
 
As we saw in Example~\ref{det-example}, Theorem~\ref{cof} fails if $p=i$.  Nevertheless, 
the results can be much improved if we additionally assume that, for every $k\in\N$, $\mathbb{T}_i(\nabla u_k)\ge 0$ (element-wise) almost everywhere in $\O$. Indeed,  M\"{u}ller \cite{mueller-det} proved the following result.

\medskip 

\begin{proposition}[Higher integrability of determinant]
\label{high-int-det}
Assume that $\omega\subset\O\subset\R^n$ is compact, $u\in W^{1,n}(\O;\R^n)$, and that $\det\nabla u\ge 0$ almost everywhere in $\O$. Then  
\begin{align}\label{determinant-log}
\|(\det\nabla u)\ln(2+\det\nabla u)\|_{L^1(\omega)}\le C(\omega, \|u\|_{W^{1,n}(\O;\R^n)})\ 
\end{align}
for some  $C(\omega, \|u\|_{W^{1,n}(\O;\R^n)})>0$ a constant depending only on $\omega$ and the Sobolev norm of $u$ in $\O$.
\end{proposition} 

\medskip 
This proposition results in the following corollary:
\medskip
 \begin{corollary}[Uniform integrability of determinant]\label{muller}
If  $\{u_k\}_{k\in\N}\subset W^{1,n}(\O;\R^n)$ is bounded and $\det\nabla u_k\ge 0$ almost everywhere in $\O$ for all $k\in\N$ then $\det\nabla u_k\rightharpoonup \det\nabla u$ in $L^1(\omega)$ for every compact set $\omega\subset\O$.
\end{corollary}

\medskip

A related statement was achieved by Kinderlehrer and Pedregal in \cite{k-p2}.  It says that under the assumptions  of Corollary~\ref{muller}  and if 
$u_k=u$ on $\partial\O$ for all $k\in\N$ the claim of Corollary \ref{muller} holds for $\omega:=\Omega$.

\medskip

\FIRST
\begin{remark}
\label{remark-distortion-astla}
Proposition \ref{high-int-det} can be strengthened if $\det \nabla u $ of a mapping  $u\in W^{1,n}(\O;\R^n)$  is not only positive but also if the following inequality is valid for some $K \geq 1$
\begin{equation}
\label{quasiregularity}
\vert \nabla u(x)\vert^n \leq K \det \nabla u(x) \qquad \text{a.e. in $\Omega$};
\end{equation}
such mappings are called \emph{quasiregular} (and if $u$ is additionally a homeomorphism \emph{quasiconformal}) and we shall encounter them again in Section \ref{sec-quasicon-elas}. In the case of quasiregular mappings, we have even that $\det \nabla u \in L^{1+\varepsilon}(\Omega)$ with $\varepsilon > 0$ depending only on $K$ and the dimension $n$ (see e.g. \cite{Hencl} where also generalizations of this result for $K$ depending on $x$ are discussed). In the quasiconformal case in dimension 2, this observation goes back to Bojarski \cite{bojarski}; in this case even the precise value of $\varepsilon< \frac{1}{K-1}$ has been established by Astala \cite{astala-dist}.

\end{remark}
\EOR

\medskip

\subsection{Explicit characterization of null Lagrangians of higher order}

Null Lagrangians of higher order are of the same structure as those of the first order. Indeed, they also correspond to minors. In order to make the statement more precise, we assume that $K:=(k_1,\ldots,k_r)$ is such that $1\le k_i\le n$ and denote by  $\alpha:=(\nu_1,J_1;\nu_2,J_2;\ldots; J_r,\nu_r)$ with $|J_i|=k-1$ and where $1\le\nu_i\le m$. We define the $k$-th order Jacobian determinant $J^\alpha_K:X\to\R$ by the formula
$$
J^\alpha_K(\nabla u)=\det\left(\frac{\partial u^{\nu_i}_{J_i}}{\partial x^{k_j}}\right)\ .$$

Then any null Lagrangian of higher order is just an affine combination of $J^\alpha_K$, i.e., we have the following theorem

\medskip 

\begin{theorem}[See {\color{black} Ball, Currie and Olver} \cite{ball-curie-olver}]\label{HONL}
Let $L\in C(X(n,m,k))$. Then $L$ is a null Lagrangian if and only if  it is an affine combination of $k$-th order Jacobian determinants, i.e., 
$$
L= C_0+\sum_{\alpha,K}C^\alpha_K J^\alpha_K $$ 
for some constants $C_0$ and $C^\alpha_K$. 
\end{theorem}

\subsection{Null Lagrangians with lower order terms}

As pointed out in Remark \ref{rem-genNull}, the notion of null Lagrangians can also be generalized  to functionals of the type \eqref{functional1}, i.e., those containing also lower order terms. A characterization of these null Lagrangians is due to Olver and Sivaloganathan \cite{olver-si} who considered the first order case; i.e., null Lagrangians for those functionals which can also depend on $x$ and $u$. 

Based on Olver's results \cite{olver}, they showed in \cite{olver-si} that such null Lagrangians are given by the formula
$$
 L(x,u,\nabla u)=C_0(x,u)+\sum_iC_i(x,u)\cdot\mathbb{T}_i(\nabla u)\ ,$$
where $C_0$ and $C_1$ are $C^1$-functions. This means that null Lagrangians with lower order terms are determined by the already known null Lagrangians of the first order. Let us remark, that it is noted in \cite{olver-si} that the result generalizes analogously to the higher order case.

\section{Null Lagrangians at the boundary}\label{sec-weakConBD}
\label{sec-NullBd}

We have seen that null Lagrangians of the first order are exactly those functions that fulfill \eqref{quasiconvexity} in the definition of quasiconvexity with an equality. This, of course, assures that null-Lagrangians are weakly* continuous with respect to the $W^{1, \infty}(\Omega; \R^m)$ weak* topology; in addition, due to Theorem \ref{cof}, the are weakly continuous with respect to the $W^{1, p}(\Omega; \R^m)$ weak topology  if $p > \min(m,n)$ with $\Omega \subset \R^n$. 

However, in the critical case when $p = \min(m,n)$ the weak continuity fails. In fact, as we have seen in Example \ref{det-example}, for $n=m=p=2$ the functional \eqref{functional1} with $k=1$ and $f(x,u,\nabla u) = \det(\nabla u)$ is not even weakly lower semicontinuous, even though the determinant itself is definitely a null Lagrangian. Once again, the reason for the failure of weak continuity are concentrations on the boundary combined with the fact that null Lagrangians are unbounded from below. 

Nevertheless, as we have seen in Section \ref{sec-BdCon}, at least for $p$-homogeneous functions, weak lower semicontinuity can be assured for functionals with integrands that are quasiconvex at the boundary; i.e., fulfill \eqref{qcbinfty}. Thus, a proper equivalent of null Lagrangians in this case are those functions that fulfill \eqref{qcbinfty} with an equality---these functions are referred to as \emph{null Lagrangians at the boundary}. We study these functions in this section.

Clearly, null Lagrangians at the boundary form a subset of null Lagrangians of the first order. Moreover, they have exactly the sought  properties:  We know from Theorem \ref{char-HONL} that if $\cN$ is a null Lagrangian at the boundary then it is a polynomial of degree $p$ for some $p \in [1,\min(m,n)]$. 
If, additionally $\{u_k\}_{k\in\N}\subset W^{1,p}(\O;\R^m)$ converges weakly to $u\in W^{1,p}(\O;\R^m)$ then $\{\cN(\nabla u_k)\}_{k\in\N}\subset L^1(\O)$ weakly* converges to $\cN(\nabla u)$ in $\mathcal{M}(\overline{\Omega})$, i.e., in measures on the closure of the domain. This means that the $L^1$-bounded sequence  $\{\cN(\nabla u_k)\}$ converges to a Radon measure whose singular part vanishes. Thus, functionals with integrands that are null-Lagrangians at the boundary are weakly continuous even in the critical case. Null Lagrangians at the boundary can be also used to construct functions quasiconvex at the boundary; cf. Definition \ref{qcb-def}.

We first give a formal definition of null Lagrangians at the boundary.

\medskip

\begin{definition}\label{nulllagrangian}
Let $\varrho\in\R^n$ be a unit vector and let $L:\R^{m\times n}\to\R$ be a given function.
\begin{enumerate}
\item[(i)]
$L$ is called a {\it null Lagrangian at the boundary} at given $A\in \R^{m\times n}$
if  both  $L$ and $-L$ are quasiconvex at the boundary at $A$ in the sense of Definition \ref{qcb-def}; cf.~\cite{silhavy}. This means that there is $q\in\R^m$ such that for all $\varphi\in W^{1,\infty}_{{\color{black} D_\varrho \setminus \Gamma_\varrho}}(D_\varrho;\R^m)$ it holds 
\be\label{null-identity1}
\int_{\Gamma_\varrho} q\cdot \varphi(x)\,\md S + L(A)\mathcal{L}^n(D_\varrho) =  \int_{D_\varrho} L(A+\nabla \varphi(x))\,\md x\ . \ee
 \item[(ii)] If $L$ is a {\it null Lagrangian at the boundary} at every $F\in \R^{m\times n}$, we call it a null Lagrangian at the boundary.
\end{enumerate}
\end{definition}

\medskip



The following theorem explicitly characterizes all possible null Lagrangians at the boundary. It was first proved by  Sprenger in his thesis \cite[Satz 1.27]{sprenger}. Later on, the proof was slightly simplified in \cite{kkk}. Before stating the result we recall that
${\rm SO}(n):=\{R\in\R^{n\times n};\ R^\top R=RR^\top=\mathbb{I}\ ,\ \det R=1\}$ denotes the set of orientation-preserving rotations and if we write $A=(B|\varrho)$ for some $B\in\R^{n\times(n-1)}$ and $\varrho\in\R^n$ then $A\in\R^{n\times n}$, its last column is $\varrho$ and $A_{ij}=B_{ij}$ for $1\le i\le n$ and $1\le j\le n-1$. 
We remind also that $\mathbb{T}_i(A)$ denotes the vector of all subdeterminants of $A$ of order $i$.
\medskip

\begin{theorem}\label{thm:bnulllag} Let $\varrho\in\R^n$ be a unit vector
and let $L:\R^{m\times n}\to\R$ be a given continuous function.
Then the following three statements are equivalent.
\begin{itemize}
\item[(i)] $\cN$ satisfies \eqref{null-identity1} for every $F\in \R^{m\times n}$;
\item[(ii)] $\cN$ satisfies \eqref{null-identity1} for $F=0$,
\item[(iii)] There are constants $\tilde\beta_s\in\R^{{m\choose s}\times {{n-1}\choose s}}$, $1\le s\le\min(m,n-1)$, such that for all $H\in\R^{m\times n}$,
\be\label{null-form}
	\cN(H)=\cN(0)+\sum_{i=1}^{\min(m,n-1)}\tilde\beta_i\cdot \mathbb{T}_i (H \tilde R),
\ee 
where  $\tilde R\in\R^{n\times(n-1)}$ is a matrix such that $R=(\tilde R|\varrho)$ belongs to ${\rm SO}(n)$;
\item[(iv)] $\cN(F+a\otimes \varrho)=\cN(F)$ for every $F\in \R^{m\times n}$
and every $a\in \R^m$.
\end{itemize}
\end{theorem}

\medskip

If $m=n=3$ the only nonlinear null Lagrangian at the boundary with  the normal $\varrho$ is $$\cN(F)={\rm Cof }\, F\cdot(a \otimes \varrho)=a\cdot {\rm Cof }\, F\varrho$$ where $a\in\R^3$ is some fixed vector; see \v{S}ilhav\'{y} \cite{silhavy}. 

In the following theorem, we let  $\varrho$ freely move along the boundary which introduces an $x$-dependence to the problem. Then the vector $a$ may depend on $x$ as well. 

\medskip

\begin{theorem}[due to \cite{mk-cocv}]\label{cof1}
Let $\O\subset\R^3$ be a smooth bounded  domain.  Let $\{u_k\}\subset W^{1,2}(\O;\R^3)$ be such that 
$u_k \rightharpoonup u$ in $W^{1,2}(\O;\R^3)$. Let $\tilde{L}(x,F):={\rm Cof }\, F\cdot(a(x) \otimes \varrho(x))$, where $a,\varrho\in C(\overline{\Omega};\R^3)$ and $\varrho$ coincides at $\partial\O$ with the outer unit normal to $\partial\O$. Then  for all $g\in C(\overline{\Omega})$  
\begin{align}\label{weakcont}
\lim_{k\to\infty}\int_{\O} g(x)\tilde{L}(x,\nabla u_k(x))\, \md x= \int_{\O} g(x)\tilde{L}(x,\nabla u(x))\,\md x\ .\end{align}
If, moreover, for all $k\in\N$ $\tilde{L}(\cdot,\nabla u_k)\ge 0$ almost everywhere in $\O$  then $\tilde{L}(\cdot,\nabla u_k)\rightharpoonup \tilde{L}(\cdot,\nabla u)$ in $L^1(\O)$.
\end{theorem}

\medskip

Notice that even though $\{\tilde{L}(\cdot,\nabla u_k)\}_{k\in\N}$ is bounded merely in $L^1(\O)$  its weak* limit in measures is $\cN(\cdot,\nabla u)\in L^1(\O)$, i.e., a measure which is absolutely continuous with respect to the Lebesgue measure on $\O$. This holds  independently of $\{\nabla u_k\}$. Therefore, the fact that $\tilde{L}$ is a  null Lagrangian at the boundary automatically improved regularity of the limit measure, namely its singular part vanishes.
In order to understand why this happens, denote $\mathbb{P}(x):=\mathbb{I}-\varrho(x)\otimes\varrho(x)$ the orthogonal projector on the plane with the normal $\varrho(x)$, i.e., a tangent plane to $\partial\O$ at $x\in\partial\O$. Then 
$${\rm Cof}(F\mathbb{P})={\rm Cof }F {\rm Cof} \mathbb{P}=({\rm Cof}F)(\varrho\otimes\varrho)\ .$$
Consequently, 
$$
{\rm Cof}(F\mathbb{P})\varrho=({\rm Cof}F)\varrho\ ,$$
and if we plug in $\nabla u$ in the position of $F$, we see that $\tilde{L}(x,\cdot)$ only depends on the surface gradient of $u$. In other words, concentrations in the sequence of normal derivatives, $\{\nabla u_k\cdot(\varrho\otimes\varrho)\}_{k\in\N}$, are filtered out.  

The following two statements describing weak sequential continuity of null Lagrangians at the boundary can be found in \cite{kkk}. Here, a similar effect as the one on Theorem \ref{high-int-det} and Corollary \ref{muller} is observed, namely the non-negativity of the null Lagrangian allows to prove weak continuity.

\medskip

\begin{theorem}[see \cite{kkk}]\label{th:weakcontinuityup}
Let $m,n\in\N$ with $n\geq 2$, let $\Omega\subset \R^n$ be open and bounded with a  boundary of class $C^1$, and let
$L:\overline{\Omega}\times \R^{m\times n}\to \R$ be a continuous function. In addition, suppose that for every $x\in \Omega$, $L(x,\cdot)$ is a null Lagrangian and
for every $x\in \partial\Omega$, $L(x,\cdot)$ is a null Lagrangian at the boundary with respect to $\varrho(x)$, the outer normal to $\partial\Omega$ at $x$. Hence, by Theorem \ref{thm:bnulllag}, $L(x,\cdot)$ is a polynomial, the degree of which we denote by $d_{\tilde{L}}(x)$.
Finally, let $p\in (1,\infty)$ with $p\geq d_f(x)$ for every $x\in\overline{\Omega}$
and let $\{u_k\}\subset W^{1,p}(\Omega;\R^m)$ be a sequence such that $u_k\rightharpoonup u$ in $W^{1,p}$.
If
$$
	L(x,\nabla u_k(x))\geq 0~~\text{for every $k\in\N$ and a.e.~$x\in\Omega$,}
$$
then $L(\cdot,\nabla u_n)\rightharpoonup L(\cdot,\nabla u)$ weakly in $L^1(\O)$.
\end{theorem}

\medskip
The above theorem allows us to prove a weak lower semicontinuity result for convex functions of null Lagrangians at the boundary which relates to the concept of polyconvexity introduced in Section \ref{sec-polyconvexity}.
\medskip

\begin{theorem}[see \cite{kkk}] Let $h:\O\times\R\to\R\cup\{+\infty\}$ be such that $h(\cdot, s)$ is measurable for all $s\in\R$ and  $h(x,\cdot)$ is  convex for almost all $x\in\O$.  Let $L$ and $d_{L}$  be as in Theorem~\ref{th:weakcontinuityup}. Then $\int_\O h(x,L(x,\nabla u(x)))\,\md x$ is weakly lower semicontinuous on the set $\{u\in W^{1,p}(\O;\R^m); L(\cdot,\nabla u)\ge 0 \mbox{ in }\O\}$.
\end{theorem}

\medskip

Let us finally point out that $A\mapsto h(L(x,A))$ for a convex function $h$ is quasiconvex at the boundary with respect to the normal $\varrho(x)$. Therefore null Lagrangians at the boundary allow us to construct functions which are quasiconvex at the boundary.





\medskip


\section{Polyconvexity and applications to hyperelasticity}
\label{sec-polyconvexity}
We saw that, at least for integrands bounded from below and satisfying (i) in Definition \ref{class-F_R}, quasiconvexity is an equivalent condition for weak lower semicontinuity. This presents an \emph{explicit} characterization of the latter since it is not necessary to examine all weakly converging sequences. Nevertheless, in practice quasiconvexity is almost impossible to verify since, in a sense, its verification calls for solving a minimization problem itself. Therefore, it is desirable to find at least \emph{sufficient} conditions for weak lower semicontinuity that can be easily verified. Such a notion, called \emph{polyconvexity}, introduced  by J.M. Ball and can be designed by employing the null Lagrangians introduced in the last section. 

We start with the definition of polyconvexity suitable for first order functionals.

\medskip

\begin{definition}[Due to Ball \cite{ball77}]\label{def-polyconvexity}
We say that 
$f:\R^{m\times n}\to\R\cup\{+\infty\}$ is polyconvex if there exists a convex function $h:\R^\sigma\to\R\cup\{+\infty\}$ such that $f(A)=h(\mathbb{T}(A))$\footnote{Recall that $\mathbb{T}(A)$ denotes the vector of all minors of $A$.} for all $A\in\R^{m\times n}$.
\end{definition}

\medskip

\begin{remark}
Interestingly, already Morrey in \cite[Thm.~5.3]{morrey-orig} proved that one-homogeneous convex functions depending on minors are quasiconvex.
\end{remark}

\medskip 

If $h$ is affine in the above definition, we call $f$ polyaffine.
In this case, $f(A)$ is a linear combination of all minors of $A$ plus a real constant. Consequently, any polyconvex function is bounded from below by a polyaffine function. Similarly, as in the convex case, a polyconvex function is found by forming the supremum of all polyaffine  functions lying below it see e.g.~\cite[Rem.~6.7]{dacorogna}; i.e., we have the following lemma.

\medskip

\begin{lemma}
The function $f: \mathbb{R}^{m \times n} \to \mathbb{R}$ is polyconvex if and only if 
$$
f(A) = \sup \{\varphi(A); \varphi \text{ polyaffine and } \varphi \leq f\}.
$$
\end{lemma} 

\medskip

It is  straightforward to generalize  polyconvexity to higher-order variational problems, i.e., those that depend on higher-order gradients of a mapping. 
The attractiveness of such problems for applications is clear. Suitably chosen terms depending on higher-order gradients allow for compactness of a minimizing sequence in some stronger topology on $W^{1,p}(\Omega; \R^m)$, which enable us to pass to a limit in lower-order terms without restrictive assumptions on their convexity properties. Thus, for example, models of shape memory alloys (see Section \ref{sec-quasicon-elas}) can be treated by this approach; cf.~e.g.~\cite{mueller-93,mueller}.
 
We extend the notion of polyconvexity  to higher order problems \eqref{functional1} by employing the notion of null Lagrangians of higher order due to  Ball, Currie, and Olver \cite{ball-curie-olver}. 

\medskip

\begin{definition}[Higher-order polyconvexity]  Let $1 \le  r  \le n$. Let $U\subset X(n,m,k)$ be open.  A function $G:U\to\R$ is $r$-polyconvex  if
there  exists  a  convex  function  $h: {\rm Co}(J^{[r]}(U))\to\R$  such  that
$f(A)  =h(J^{[r]}(A))$    for all  $A\in U$; here ${\rm Co}(J^{[r]}(U))$ is the convex hull of $J^{[r]}(U)$.  $G$ is  polyconvex  if it is $R$-polyconvex. Here $J^r(H):=(J^{r,1}(H),\ldots, J^{r,N_r}(H))$ 
is a $N_r$-tuple with the property that any Jacobian determinant of degree $r$ can be written as a linear combination of elements of $J^r$. Consequently, $J^{[r]}:=(J^1,\ldots, J^r)$.
If $h$  is affine then we call $f$ $r$-polyaffine.
\end{definition}

\medskip

Since polyconvexity implies quasiconvexity, we may deduce by the results in Section \ref{result-meyers} that polyconvex functions in the class $\mathcal{F}_p(\Omega)$ (from Definition \ref{class-F_R}) are weakly lower semicontinuous. Yet, weak lower semicontinuity can be proved for  wider class of polyconvex functions than those in $\mathcal{F}_p(\Omega)$; in particular, the functions do not have to be of $p$-growth. This is of great importance in elasticity as explained later in this section.  
 
The proof of weak lower semicontinuity of polyconvex functions can be actually based on \emph{convexity} and weak continuity of null Lagrangians. Thus, because weak lower semicontinuity can be shown for arbitrarily growing convex functions, it generalizes to polyconvex ones, too. The following result for convex functions can be found in \cite[Thm.~5.4]{ball-curie-olver} and is based on results by Eisen \cite{eisen} who proved this theorem for $\Phi$ finite-valued.  

\medskip

\begin{theorem}[weak lower semicontinuity]\label{thm:Eisen}
Let $\Phi:\O\times\R^s\times\R^\sigma\to\R\cup\{+\infty\}$ satisfy the following properties
\begin{enumerate}[(i)]
\item $\Phi(\cdot,z,a):\O\to\R\cup\{+\infty\}$ is measurable for all $(z,a)\in \R^s\times\R^\sigma$,
\item $\Phi(x,\cdot,\cdot): \R^s\times\R^\sigma\to\R\cup\{+\infty\}$ is continuous for almost every $x\in\O$,
\item $\Phi(x,z,\cdot):\R^\sigma\to\R\cup\{+\infty\}$ is convex.
\end{enumerate}

Assume further that there is $\phi\in L^1(\O)$ such that  $\Phi(\cdot,z,a)\ge\phi$ for all $(z,a)\in \R^s\times\R^\sigma$.
Let $\{z_k\}_{k \in \N} \subset L^1(\O;\R^s)$, $\{a_k\}_{k \in \N} \subset L^1(\O;\R^\sigma)$ and let $z_k\to z$ almost everywhere in $\O$ as well as $a_k\wto a$ in $L^1(\O;\R^\sigma)$. Then 
$$\int_\O\Phi(x,z(x),a(x))\,\md x\le\liminf_{k\to\infty}\int_\O\Phi(x,z_k(x),a_k(x))\,\md x\ .$$
\end{theorem}

\medskip

Using this theorem, we may easily deduce weak lower semicontinuity of polyconvex functions. For the sake of clarity, let us start with first order problems. Then, consider $u_k\wto u$ in $W^{1,p}(\O;\R^m)$ as $k\to\infty$  where $p>\min(m,n)$. Then $u_k\to u$ in $L^p(\O;\R^m)$, so, for a (non-relabeled) subsequence, even $u_k\to u$ almost everywhere in $\O$. Hence, we can apply Theorem~\ref{thm:Eisen} with $z_k:=u_k$, $a_k:=\mathbb{T}(\nabla u_k)$ and $f(x,u,\nabla u):=\Phi(x,u,\mathbb{T}(\nabla u))$  to obtain the following corollary:

\medskip

\begin{corollary}\label{cor:polyconvexity}
Let $f:\O\times\R^m\times\R^{m\times n}\to\R\cup\{+\infty\}$ satisfy the following properties
\begin{enumerate}[(i)]
\item $f(\cdot,r,A):\O\to\R\cup\{+\infty\}$ is measurable for all $(r,A)\in R^m\times\R^{m\times n}$,
\item $f(x,\cdot,\cdot): \R^m\times\R^{m\times n}\to\R\cup\{+\infty\}$ is continuous for almost every $x\in\O$,
\item $f(x,r,A)=\Phi(x,r, \mathbb{T}(A))$ where $\Phi$ satisfies (i)--(iii) from Theorem~\ref{thm:Eisen}.
\end{enumerate}

If $u_k\wto u$ in $W^{1,p}(\O;\R^m)$ as $k\to\infty$  where $p>\min(m,n)$ then 
$$\int_\O f(x,u(x),\nabla u(x))\,\md x\le\liminf_{k\to\infty}\int_\O f(x,u_k(x),\nabla u_k(x))\,\md x\ .$$
\end{corollary}

\medskip

Similarly as in the case of first order problems, we can exploit (v) of Theorem~\ref{char-HONL} and Theorem~\ref{thm:Eisen} to show the existence of minimizers to energy functionals \eqref{functional1}. Let us present the result just for functionals \eqref{functional1} with $k =2$; generalizations for higher $k$ are straightforward and can be found in \cite{ball-curie-olver}. 

\medskip

\begin{corollary}[after \cite{ball-curie-olver}]
\label{2nd-grade}
Assume that $\O\subset\R^n$ is a bounded smooth domain and that $1\le r\le n$. Let  $f:\O\times Y(n,m,2)\to\R\cup\{+\infty\}$ and
$$
I(u) = \int_\Omega f(x,u,\nabla u, \nabla^2 u) \md x
$$
satisfy the following assumptions:
\begin{itemize}
\item[(i)] $f(x,H,A)=h(x,H,J^{[r]}(A))$, 
where $h(x,\cdot,\cdot):(\R^m \times \R^{m\times n})\times J^{[r]}(X(n,m,2))\to\R\cup\{+\infty\}$ is continuous for almost every $x\in\O$, 
\item[(ii)] $h(\cdot,H,J^{[r]}(A)):\O\to\R\cup\{+\infty\}$ is measurable for all $(H,J^{[r]}(A)) \in (\R^m \times \R^{m\times n})\times J^{[r]}(X(n,m,2))$,
\item[(iii)] $h(x,H, \cdot):J^{[r]}(X)\to\R\cup\{+\infty\}$ is convex for almost all $x\in\O$ and all $H\in(\R^m \times\R^{m\times n})$,
\item[(iv)] $f(x,H,A)\ge C(-1+|A|^p)$ for some $C>0$, $p>n$, almost all $x\in\O$ and all $A\in\R^{m\times n}$.
\end{itemize}
Then $I$ is weakly lower semicontinuous on $W^{2,p}(\O;\R^m)$.

\end{corollary}

\medskip

{\color{Green} It cannot be stressed enough} that the main strength of polyconvexity consists in the fact that convexity in subdeterminants can be advantageously combined with the Mazur lemma to show weak lower semicontinuity in a similar way like in the proof for mere convex and lower semicontinuous integrands. This contrasts with proofs  available for quasiconvex integrands where 
manipulations with boundary conditions are needed to prove the result. This is already clearly visible in 
Meyers paper \cite{meyers}. These manipulations typically destroy any pointwise constraints on the determinant of $\nabla y$, which, however, are crucial in elasticity. We shall return to this issue in Section \ref{sec-quasicon-elas}.

\SECOND
It is well known that minimizers of integral functionals with strictly convex integrands are unique. However, the same is not true for variational problems with polyconvex integrands even if $h$ in Definition~\ref{def-polyconvexity}  is strictly convex.  Examples  were provided  by Spadaro \cite{spadaro} if $m=n=2$.
\EOR

\subsection{Rank-one convexity}   

Since polyconvexity is an explicit sufficient condition for quasiconvexity, we may ask if similarly a simpler  necessary condition can be found. This is indeed so, the sought notion of convexity is the so-called \emph{rank-one convexity}:

\medskip

\begin{definition}[Due to \cite{morrey}]
\label{defin-rank1}
We say that $f: \mathbb{R}^{m \times n} \to \mathbb{R}$ is rank-one convex if
\begin{equation}
\label{rank1_def}
f(\lambda A_1 + (1-\lambda)A_2) \leq  \lambda f(A_1)+(1-\lambda)f(A_2).
\end{equation}
for all $\lambda \in [0,1]$ and all $A_1$, $A_2$ such that $\mathrm{rank}(A_1-A_2) \leq 1$.
\end{definition}

\medskip

The relations among the introduced  notions of convexity are as follows:
$$
\text{convexity} \Rightarrow \text{polyconvexity} \Rightarrow \text{quasiconvexity} \Rightarrow \text{rank-one convexity};
$$
however, none of the converse implications holds if $f: \mathbb{R}^{m \times n} \to \mathbb{R}$ and $m \geq  3$ and  $n \ge 2$.  To see that polyconvexity does not imply convexity (even for $m,m > 1$) just consider the function $f(F):= \det(F)$ which is even polyaffine but not convex. Also quasiconvexity does not imply polyconvexity even for  $m,n > 1$ as was shown in e.g.~\cite{alibert-dacorogna,terpstra}. \v{S}ver\' ak's important counter example \cite{sverak} is a construction of a function that is rank-one convex, but not quasiconvex and holds for $m\ge 3$ and $n\ge 2$. For $m=2$ and $n\ge 2$  the question of equivalence between quasiconvexity and rank-one convexity is still unsolved (see Open problem \ref{rank1open}). Notice that, if $m=1$ or $n=1$ all the generalized notions of convexity trivially coincide with standard convexity itself. 

\medskip

\begin{problem}
\label{rank1open}
Let $m=2$ and $n\ge 2$. Does rank-one convexity imply quasiconvexity for $f: \mathbb{R}^{m \times n} \to \mathbb{R}${\color{black}?}
\end{problem}

\medskip

\FIRST Many attempts can be found in literature towards the solution of Open problem \ref{rank1open} and indication both in the positive and in the negative exist. Morrey conjectured in his original paper \cite{morrey-orig} that the solution to Open problem \ref{rank1open} is negative. On the other hand, in several special cases it has been shown that rank-one convexity indeed implies quasiconvexity.
For example, if $f$ is a quadratic form, i.e., if there exists $\mathbb{C}\in\R^{m\times n\times m\times n}$  such that  $f(A)=\mathbb{C}A\cdot A:=\sum_{ijkl}\mathbb{C}_{ijkl}A_{ij}A_{kl}$ then $f$ is rank-one convex if and only if it is quasiconvex \cite{vanhove}. Additionally, if $\min(m,n)=2$ then any quasiconvex (or rank-one convex) quadratic form is even polyconvex \cite{marcellini2form,serre, terpstra}. Also, as has been realized by M\"uller \cite{muller-matrices}, rank-one convexity is equivalent to quasiconvexity on diagonal matrices. Very recently, Open problem \ref{rank1open} has been answered positive for isotropic, objective and isochoric elastic energies \cite{neff}. On the top of that, even in the case $m \geq  3$ and  $n \ge 2$ only very few examples of functions that are rank-1 convex but not quasiconvex are known. Besides the already mentioned \v{S}ver\' ak's counterexample, Grabovsky \cite{grabovsky} provided very recently a new example of a rank-one convex function which is not quasiconvex in dimensions $m=8$ and $n=2$.

On the other hand, let us also mention that if $m=n=2$ a necessary condition for the equivalence of rank-one convexity and quasiconvexity  is that every $f:\R^{2\times 2}\to\R$ quasiconvex satisfies the inequality $$f(A)\le \int_{(0,1)^2}f(A+(\nabla\varphi(x))^\top)\,\md x \quad \text{for every } \varphi\in W^{1,\infty}_0((0,1)^2;\R^2) \text{ and every } A\in\R^{2\times 2}.$$ An analogous implication for $m=n>2$ turns out to be false \cite{kruzik-transposition,mueller-transposition}, however, the two-dimensional case is open.

Before continuing to higher-dimensional equivalents of rank-one convexity, let us still point out one important function that may form a counterexample in Open problem \ref{rank1open}, the so-called Burkholder function
$$
B_p(A)=\frac{p}{2} \det(A) \vert A \vert^{p-2} + \big(1-\frac{p}{2}\big)\vert A \vert^p \qquad \text{with $p \geq 2$,}
$$
with the $\vert \cdot \vert$ denoting the operator norm. The Burkholder function emerged in the study of stochastic integrals and martingales \cite{Burkholder1,Burkholder2}, it is rank-one concave \cite{IwaniecRank1} (which means that $-B_p$ is rank-one convex). Nevertheless, it is a standing problem whether it is also quasiconcave (i.e. $-B_p$ quasiconvex).

\medskip 
\begin{problem}
\label{iwaniec-conj}
Is the Burkholder function $B_p(\cdot)$ quasiconcave? 
\end{problem}

\medskip

Naturally, Open problem  \ref{iwaniec-conj} is answered affirmatively if the answer to Open problem \ref{rank1open} is ``yes''. Yet, Open problem \ref{iwaniec-conj} is interesting for its own right due to its implications in harmonic and quasiconformal analysis. In fact, if one can prove that the Burkholder function is quasiconcave this would give a precise bound on the operator norm of the Beurling-Ahlfors transform $S: L^p(\mathbb{C}) \to L^p(\mathbb{C})$ (here $\mathbb{C}$ denotes the complex plane), which plays an important role in complex analysis since it converts the complex partial derivative $\partial_{\bar{z}}$ into $\partial_z$. Indeed, the standing conjecture due to Iwaniec \cite{Iwaniec82} is that this norm is equal to $p^\star - 1 := \max\{p-1,1/(p-1)\}$; cf. the reviews \cite{IwaniecConjrev1,IwaniecConjrev2}. It is classical that $p^\star -1$ is the lower bound for the norm of $S$ \cite{Lehto}; as for the upper bound, the attempts in literature have progressively come closer to the number $p^\star -1$ (see the review \cite{IwaniecConjrev1} and also e.g. the current improvements \cite{pstar1,pstar2}) but it is not reached as of today. Reaching this upper bound would for example imply the distortion result by Astala mentioned in Remark \ref{remark-distortion-astla}.

In the context of Open problem \ref{iwaniec-conj} and \ref{rank1open} let us also point to the recent work of Astala, Iwaniec, Prause and Saksma \cite{Astala-Burkholder} in which the authors show that the Burkholder function is quasiconcave in the identity for quasiconform perturbations\footnote{By a \emph{perturbation} we mean the function $\varphi$ in Definition \ref{def-quasiconvexity}}; in other words, in the language of Section \ref{sec-quasicon-elas}, the Burkholder function is quasiconformally quasiconcave at the identity.
\EOR

An equivalent to rank-one convexity can also be defined for higher-order problems: the corresponding notion is called $\Lambda$-convexity. Following \cite{ball-curie-olver}, we define a nonconvex cone $\Lambda\subset X(n,m,k)$ as 
$\Lambda:=\{a\otimes^k b:\, a\in\R^m,\, b\in\R^n\}$ where 
$(a\otimes^k b)^i_K=a^ib_K$. Here recall from Section \ref{sec-notation}, that $b_K =(b_1)^{k_1}(b_2)^{k_2}\ldots(b_n)^{k_n}$ and $K=(k_1, \ldots, k_n)$ with $1 \leq k_i\leq k$ for $i=1\ldots n$ is a multiindex.
\medskip
\begin{definition}
 A function $f:X\to\R$ is called  $\Lambda$-convex if
$t\mapsto f(A+tB):\R\to\R$ is convex for every $A\in X(n,m,k)$ and any $B\in\Lambda$.
\end{definition}

\medskip 

Notice that for $k=1$ $\Lambda$-convexity coincides with rank-one convexity. If $f$ is twice continuously differentiable  then $\Lambda$-convexity is equivalent to the Legendre-Hadamard condition
$$
\sum_{j,i=1}^m\sum_{|J|=|K|=k}\frac{\partial^2f(A)}{\partial A^j_J\partial A^i_K}a^ja^ib_Jb_K\ge 0
$$
for all $A\in X(n,m,k) $,  $a\in\R^m$, and $b\in\R^n$.

\medskip

\begin{proposition}[see \cite{ball-curie-olver}]
Continuous and $k$-quasiconvex functions $f:X(n,m,k)\to\R$ are $\Lambda$-convex. 
\end{proposition}

\medskip
Hence, $\Lambda$-convexity forms a necessary condition for $k$-quasiconvexity. 
This proposition was first proved by Meyers \cite[Thm.~7]{meyers} for smooth functions and then generalized by Ball, Currie, and Olver \cite{ball-curie-olver} to the continuous case.  The opposite assertion does not hold. Indeed, if $n=k=2$ and $m=3$ then we have the following example due to Ball, Currie, and Olver for $f:X(2,3,2)\to\R$
$$f(\nabla^2u)=\sum_{i,j,l=1}^3\varepsilon_{ijk}\frac{\partial^2 u^i}{\partial x_1^2}\frac{\partial^2 u^j}{\partial x_1\partial x_2}\frac{\partial^2 u^l}{\partial x_2^2}\ ,$$
where $\varepsilon_{ijk}$ is the Levi-Civita symbol. This function is even  $\Lambda$-affine (i.e., both $\pm f$ are $\Lambda$-convex)  but not a null Lagrangian and not quasiconvex. As $\Lambda$-convexity replaces rank-one convexity in the current setting we see, that this example is a reminiscent of \v{S}ver\'{a}k's example mentioned above.

\medskip

\subsection{Applications to hyperelasticity in the first order setting}

In elasticity, one is interested in modeling the response of a rubber-like material {\color{Green} exposed} to the action of applied forces. This response is obtained by solving a minimization problem; to be more specific, we are to minimize the free energy of the material. We will see that polyconvexity is perfectly fitted to the setting in elasticity and that existence of minimizers can be assured for polyconvex energies.  We give a short introduction to this matter in this section and refer the reader e.g. to the monographs \cite{gurtin,gurtin-annand,silhavy} for more details on the physical modeling.

Take a bounded Lipschitz domain $\O\subset\R^n$ which, for $n=3$, plays a role of a reference configuration of an elastic material. For given applied loads, we search for a mapping $u:\O\to\R^m$, the \emph{deformation} of the material, which describes the new ``shape'' $u(\O)$ of the body. The mapping $u$ is found by solving the following system of equations 
\begin{align}
-\text{ div } S&=b \qquad \text {in $\O$,} \label{balance-mom} \\
S\varrho&=g \qquad \text{on $\Gamma_{\rm N}$,} \label{Neumann-bc} \\
u&=u_0 \qquad \!\text{on $\Gamma_{\rm D}$.} \label{dir-bxc}
\end{align}
Here, \eqref{balance-mom} is the reduced version of Newton's law of motion for the (quasi)static case, $b$ is the applied volume force. Further, \eqref{Neumann-bc} represents the action of applied surface forces $g$ ($\varrho$ denotes the outer unit normal vector to $\Gamma_{\rm N}$) and \eqref{dir-bxc} models that the body may be clamped at some part of the boundary to a prescribed shape $u_0$. We shall require that $\Gamma_{\rm D}\subset\partial\O$  is disjoint from $\Gamma_{\rm N} \subset \partial \Omega$ and of positive $(n-1)$-dimensional Lebesgue measure. 

The material properties of the specimen are encoded in the first Piola-Kirchhoff stress tensor $S:\O\to\R^{m\times n}$ in \eqref{balance-mom} and \eqref{Neumann-bc}. The form of the Piola-Kirchhoff stress tensor cannot be deduced from first principles within continuum mechanics
but has to be prescribed phenomenologically. The prescription for $S$ is called the \emph{constitutive relation} of the given material. In the easiest case, we assume the form $S(x)=\hat S(x, \nabla u(x))$ for some given $\hat S$. Materials for which this assumption is adequate are sometimes referred to as \emph{simple} materials as opposed to non-simple materials for which $\hat S$ may depend also  on higher gradients of $u$.  Later, in Subsection \ref{sec-higherOrder}, we will consider also these sophisticated constitutive relations.     

Hyperelasticity is a part of elasticity  where  an \emph{additional assumption} is made; namely, that $S$ has a potential  $W:\R^{m\times n}\to [0;+\infty]$ such that 
$$
S_{ij}(x)=\frac{\partial W(F)}{\partial F_{ij} }|_{F=\nabla u(x)}\ .$$
This assumption emphasizes the idea that there are no energy losses in elasticity and all work, made by external forces and/or Dirichlet boundary conditions, stored in the material can be fully exploited. 

In the following, let us restrict our attention to \emph{deformations of bulks}, i.e. we do not treat plates and rods, and set thus $m=n$. In order to fulfill the basic physical requirements, $W$ has to satisfy the following relations:
\begin{align}
\label{eq:frameindiff}
W(RF)&=W(F) \text{ for all $F\in\R^{n\times n}$  and for all  $R\in {\rm SO}(n)$} \\
\label{eq:negdet}
W( F)&=+\infty \text{ if $\det F\le 0$, and }\\
\label{det2zero}
 W( F) &\to+\infty\ \text{  if $\det F\to 0_+$}.
\end{align}
Indeed, assumption \eqref{eq:frameindiff} is a consequence of the {\it axiom of frame indifference} \cite{ciarlet}; in other words the assumptions ensure that material properties are independent of the position of the observer. Conditions \eqref{eq:negdet}
and \eqref{det2zero} ensure, respectively,  that 
the material does not locally penetrate itself and that compression of a finite volume of the specimen into zero volume is not possible. These conditions, however, do not yet assure that the body does not penetrate through itself, which is also natural to assume from a physical point of view. Nevertheless, we shall see in the end of this section that with additional assumptions on the growth of the energy and, e.g., the boundary conditions even complete non-interpenetration can be assured.

The assumptions \eqref{eq:frameindiff}-\eqref{det2zero} rule out that $W$ can be convex. \SECOND Indeed, there  are matrices $A,B\in\R^{n\times n}$, both with positive determinants, such that the line segment $[A,B]$ contains a matrix $C$ with zero determinant. However, no convex function can be finite at $A,B$ and infinite at $C$.   \EOR

Moreover, due to \eqref{eq:negdet}-\eqref{det2zero} even if $W$ was quasiconvex, we could not apply theorems in Section \ref{result-meyers} since $W$ cannot be an element of the class $\mathcal{F}_p(\Omega)$.
Nevertheless, {\it polyconvexity} is fully compatible with these assumptions.

The mechanical model is that stable states of the system are found by minimizing the overall free energy
\begin{equation}
\mathcal{E}(y)=\int_\O W(\nabla u(x))\,\md x-\int_\O b(x)\cdot u(x)\,\md x-\int_{\Gamma_N} g(x)\cdot u(x)\,\md S,  \label{energy} 
\end{equation}
subject to \eqref{dir-bxc}. Smooth minimizers fulfill the balance equations \eqref{balance-mom}-\eqref{Neumann-bc}; however, even in the smooth case there might exist solutions to \eqref{balance-mom}-\eqref{Neumann-bc} which are not minimizers of \eqref{energy}. Nevertheless, such solutions are thought to be metastable and hence left after a small perturbation. Thus, minimizing \eqref{energy} is the proper way to find indeed stable states.

\medskip

\begin{remark}
Let us note that, since the minimizers of \eqref{energy} might be non-smooth, it is not guaranteed that they will satisfy the Euler-Lagrange equations either in strong or weak form. Indeed, in \cite{ball-mizel} even one-dimensional examples of smooth $W$ were given such that the minimizer does not fulfill the Euler-Lagrange equation.

One of the reasons why deducing the Euler-Lagrange equation might be difficult is that even the calculation of the variation of $\mathcal{E}$ itself can pose difficulties. Indeed, due to \eqref{eq:negdet}, the minimizer $y$ might be such that $\mathcal{E}(y+t\varphi)$ is infinite for all small enough $t>0$ and a large class of $\varphi$. Let us refer to \cite{ball-mizel} for explicit examples in which this situation occurs.
\end{remark}

\medskip 
{\color{Green}
\begin{remark}
 The way of deriving Euler-Langrange equations by taking variations of the energy functional is not completely fitting to elasticity. In particular, if two subsequent deformations are applied to an elastic body they are not added to each other but rather composed in order to obtain the final deformation of the body. 
\end{remark}
}
\medskip
\begin{remark}
Let us notice that the condition \eqref{eq:negdet} is really necessary to be stated explicitly. Namely,  from the physical point of view, the frame-indifference  \eqref{eq:frameindiff} requires that $W(F):=\tilde W(C)$ where $C:=F^\top F$ is the so-called right Cauchy-Green strain tensor. Note that $F^\top Q^\top QF=F^\top F$ for every orthogonal matrix $Q$.  Hence, pointwise minimizers of the energy density $W$  contain the set 
$\{QF_0:\, Q\in{\rm O}(n)\}$ for some given matrix $F_0$
with $\det F_0>0$ which is a pointwise minimizer itself.  Besides the physically acceptable energy wells $\{RF_0:\, R\in \text{SO}(n)\}$ other minimizers live on    wells $\{RF_0:\, R\in {\rm  O}(n)\setminus{\rm SO}(n)\} $ which are not mechanically admissible. Those wells are excluded by \eqref{eq:negdet}.
\end{remark}

\medskip 

In order to prove existence of stable states, that is minimizers of \eqref{energy}, we assume suitable {\color{Green} coercivity} of the energy density:
\begin{equation}
\label{eq:growth}
W(F)\ge c(-1+|F|^p) \text{ for all $F \in \R^{n\times n}$ and for some $c>0$,}
\end{equation}

The existence theorem follows then directly from Corollary \ref{cor:polyconvexity} where we replace $f$ by $W$.

\medskip 

\begin{theorem}
Let $\O\subset\R^3$ be a Lipschitz bounded domain, $p>3$, $u_0\in W^{1,p}(\O;\R^3)$, and $\Gamma_\mathrm{D}\subset\partial\O$ have a finite two-dimensional Lebesgue measure.  Let $W$ satisfy the assumptions (i)-(iii) from Corollary~\ref{cor:polyconvexity} posed on $f$ with $m=n=3$.  Let further 
 \eqref{eq:frameindiff}--\eqref{det2zero} and \eqref{eq:growth}  hold. If $$\mathcal{Y}:=\{u\in W^{1,p}(\O;\R^3):\, u=u_0\text{ on $\Gamma_\mathrm{D}$}, \det\nabla u>0\}\ne\emptyset$$ is such that $\inf_{\mathcal{Y}}I<+\infty$ then there is a minimizer of $\mathcal{E}$ on $\mathcal{Y}$.
\end{theorem}

\medskip

This result can be generalized for different {\color{Green} coercivity} conditions like the one considered in \eqref{eq:growth2} below. Even more general settings can be found in \cite{ciarlet} where various additional requirements on  minimizers, as e.g. conditions ensuring a friction-less contact (Signorini problem); are included, too.

Let us mention a few important examples of polyconvex stored energy densities.
Contrary to nontrivial examples of quasiconvex functions,  it is relatively easy to design a polyconvex function. To ease our notation we only define the densities for matrices of positive determinant. Otherwise, the energy density is implicitly extended by infinity. We refer to \cite{schroeder1,schroeder2,schroeder3} for examples of polyconvex functions with various special symmetries.

\medskip
\begin{example}[\it  Compressible Mooney-Rivlin material.]
This material has a stored energy of the form 
\be\label{mr}
W(F)=a|F|^2+\tilde a|\cof F|^2+\gamma(\det F)\,,
\ee
where $a,\tilde a>0$ and $\gamma(\delta)=c_1\delta^2-c_2\log\delta$, $c_1,c_2>0$.

It can be shown  that for $n=3$
$$W(F)=\frac{\lambda}{2}(\tr E)^2+\mu|E|^2 +\mathcal{O}(|E|^3)\ ,\ E=(C-\mathbb{I})/2$$
where $\lambda$ and $\mu$ are the usual Lam\'e constants, and $\mathbb{I}$ denotes the identity matrix.
Indeed, it is a matter of a tedious computation to show that, given $\lambda,\mu$,  the following equations must be fulfilled by $a,b,c_1,c_2$: $c_2:=(\lambda+2\mu)/2$
$2a+2b=\mu$, and $4b+4c_1=\lambda$.
\end{example}

\medskip
\begin{example}[Compressible neo-Hookean material.]
This material has a stored energy of the form 
\be\label{neoH}
W(F)=a|F|^2+\gamma(\det F)\,
\ee
with the same constants as for the compressible Mooney-Rivlin materials.
\end{example}
\medskip

\begin{example}[Ogden material.]
This material has a stored energy of the form (recall that $C=F^\top F$)
\be\label{ogd}
W(F)=\sum_{i=1}^Ma_i\tr \, C^{\gamma_i/2}+\sum_{i=1}^N\tilde a_i\tr (\cof \, C)^{\delta_i/2}+\gamma(\det F)
\ee
and $a_i,\tilde a_i>0$, $\lim_{\delta\to 0_+}\gamma(\delta)=+\infty$ for 
$\gamma:\R_+\to\R$ convex growing suitably at infinity.
\end{example}

\medskip

If $W$ satisfies conditions \eqref{eq:negdet}-\eqref{det2zero} then any $u \in C^1( \Omega,\R^3)$ for which $\mathcal{E}(u)$ from \eqref{energy} is finite is also locally invertible. This follows from the standard inverse function theorem. Nevertheless, what is actually desired for a physical deformation is that it is \emph{injective} \cite{ciarlet}. Indeed, non-injectivity of the deformation would mean that two material points from the reference configuration would be mapped to just one in the deformed configuration which means that the specimen penetrated through itself. Thus, additional assumptions to \eqref{eq:negdet}-\eqref{det2zero} on $W$ are needed to assure \emph{global invertibility of $u$}. Preferably, these assumptions should be compatible with polyconvexity and weak lower semicontinuity.

Take a diffeomorphism $u: \Omega \to u(\Omega)$ with $\mathrm{det} \nabla u > 0$ on $\Omega$. Then, we have by the change of variables formula for $p > 1$
$$
\int_{u(\Omega)} \!\!\!\!\!\! |\nabla u^{-1}(w)|^p \md w = \int_{\Omega} |\nabla u^{-1}(u(x))|^p \det\nabla u(x)\md x = \int_{\Omega} |(\nabla u(x))^{-1}|^p \mathrm{det}\nabla u(x)\md x =\int_{\Omega} \frac{|\cof^\top \nabla u(x)|^p}{(\det \nabla u(x))^{p-1}}\md x
$$
where we used that $\nabla u^{-1}(u(x)) = (\nabla u(x))^{-1}$ for all $x$ in $\Omega$ and that for every invertible matrix the relation $A^{-1} = \frac{\mathrm{Cof}^{\top} A}{\mathrm{det} A}$ holds. 

Therefore, for energies satisfying a stricter {\color{Green} coercivity} condition than \eqref{eq:growth} in the form of 
\begin{equation}
\label{eq:growth2}
W(F)\ge c\left(-1+|F|^p + \frac{|\cof F^\top |^p}{(\det F)^{p-1}}\right) \text{for some $c>0$,}
\end{equation}
one could rather expect that deformations on which $\mathcal{E}(u)$ is finite are invertible. This is indeed so, as theorem \ref{jmball} (below) shows.

Nevertheless, before proceeding to theorem, let us point out that the new growth condition \eqref{eq:growth2} is \emph{fully compatible with polyconvexity}. Indeed, since the function $g$ defined by  $g(x,y):= \frac{x^p}{y^{p-1}}$ is convex for $p>1$ on the set $\{(x,y) \in \R^2; y > 0\}$, $\frac{|\mathrm{cof}^\top F|^p}{(\mathrm{det}F)^{p-1}}$ is polyconvex on the set of matrices having a positive determinant.

\medskip

\begin{theorem}[Due to {\color{black} Ball} \cite{ball81}]\label{jmball}
Let $\O\subset\R^n$ be a bounded  Lipschitz domain. Let $u_0:\overline{\O}\to\R^n$ be continuous in $\overline{\O}$ and one-to-one in $\O$ such that $u_0(\O)$  is also bounded and  Lipschitz. Let $u\in W^{1,p}(\O;\R^n)$ for some $p>n$,  $u(x)=u_0(x)$ for all $x\in\partial\O$, and let $\det\nabla u>0$ a.e.~in $\O$. Finally, assume that  for some $q>n$
\begin{equation}
\int_\O|(\nabla u(x))^{-1}|^q\det\nabla u(x)\,\md x<+\infty\ .
\label{integral-Cond}
\end{equation}
Then $u(\overline{\O})=u_0(\overline{\O})$ and $u$ is a homeomorphism of $\O$ onto $u_0(\O)$. Moreover, the inverse map $u^{-1}\in W^{1,q}(u_0(\O);\R^n)$ and  $\nabla u^{-1}(w)=(\nabla u(x))^{-1}$ for $w=u(x)$ and a.a.~$x\in\O$. 
\end{theorem}

\medskip

Let us note that the Sobolev regularity needed in theorem has been weakened later in \cite{sverak-globalInv}. Indeed, in this work  it was shown that an inverse to deformation can be defined even for $p>n-1$ and $q \geq \frac{p}{p-1}$.

Theorem \ref{jmball} assures injectivity of $u$ under the growth \eqref{eq:growth2} if a up-to-the-boundary injective Dirichlet condition is prescribed. This, however, has the disadvantage that we could not model situations in which hard loads (Dirichlet boundary conditions) are prescribed only on a part on the boundary. 

One possible remedy is to minimize $\mathcal{E}$ along with the so-called \emph{Ciarlet-Ne\v{c}as} condition
\begin{equation}
\int_\Omega \mathrm{det}\nabla u(x) \md x \leq \mathcal{L}^n(u(\Omega)),
\label{Ciarlet-Necas}
\end{equation}
that was introduced in \cite{ciarlet-necas} (for $n=3$) in order to assure \emph{global injectivity} of deformations. It was shown in \cite{ciarlet-necas} that $C^1$-functions satisfying \eqref{Ciarlet-Necas} and  $\mathrm{det}\nabla u > 0$ are actually injective. The result generalizes to $W^{1,p}$-functions as well, but injectivity is obtained only almost everywhere in the deformed configuration; i.e., almost every point in the deformed configuration has only one pre-image.  

\medskip

\begin{remark}
Maps that are injective almost everywhere in the deformed configuration still include rather nonphysical situations. For example a dense, countable set of points could be mapped to one point. This can be prevented if the deformation is \emph{injective everywhere}. 

Using condition \eqref{Ciarlet-Necas}, this can be achieved for finite deformations of the energy $\mathcal{E}$ with a density $W$ satisfying \eqref{eq:growth2} for $p=m=n=2$. This setting is the most explored one due to its relations to quasiconformal maps (see Section \ref{sec-quasicon-elas}). Such deformations are open (that is they map open sets to open sets) and discrete (the set of pre-images for every point does not accumulate) and, moreover, satisfy the Lusin $N$-condition (i.e. they map sets of zero measure again to sets of zero measure); cf. e.g. \cite{Hencl}).

Then, we have by the \emph{area formula}
$$
\int_\Omega \det\nabla u\, \md x = \int_{\mathbb{R}^n} N(u,\Omega,z)\, \md z = \int_{u(\Omega)} N(u,\Omega, z)\, \md z
$$
where $N(u,\Omega, z)$ is defined as the number of pre-images of $z\in u(\O)$ in $\Omega$. So the Ciarlet-Ne\v{c}as condition is satisfied if and only if $N(u,\Omega,z) = 1$ almost everywhere on $u(\Omega)$. Also we can immediately see that the reverse inequality to (\ref{Ciarlet-Necas}) always holds.

Further, if there existed  $z \in u(\Omega)$ that had at least to two pre-images $x_1$ and $x_2$ then we could find an $\varepsilon > 0$ such that $B(x_1,\varepsilon) \cap B(x_2,\varepsilon) = \emptyset$ and $B(x_j,\varepsilon) \subset \Omega$ for $j=1,2$. On the other hand, for the images we have that $u(B(x_1,\varepsilon)) \cap u(B(x_2,\varepsilon)) \neq \emptyset$. In fact, $u(B(x_1,\varepsilon)) \cap u(B(x_2,\varepsilon))$ is of positive measure since both $u(B(x_j,\varepsilon))$  are open. Therefore, there exists a set of positive measure where $N(u,\Omega, z)$ is at least two; a contradiction to \eqref{Ciarlet-Necas}.
\end{remark}

\medskip

\subsection{Applications to hyperelasticity in the higher order setting}
\label{sec-higherOrder}
Let us now turn our attention to models of hyperelastic materials depending on higher-order gradients. Such materials are called non-simple of grade $k$, where $k$ refers to the highest derivatives appearing in the 
stored energy density.
The concept of such materials has been developing for long time, since 
the work by R.A.\,Toupin \cite{Toup62EMCS}, under various names
as non-simple materials as e.g.\ in \cite{FriGur06TBBC,Kort01FPEM, PoGiVi10HHCS,Silh88PTNB} 
or multipolar materials  (in particular 
fluids). 

Here, we  will
consider only second-grade non-simple materials, i.e., those for which second-order deformation gradients  (first-order strain gradients)  are involved. The main mathematical advantage of non-simple materials is that higher-order deformation gradients bring additional regularity of deformations and, possibly, also compactness  in a stronger topology. Moreover, {\color{Green} taking the stored energy even convex in the highest derivatives of the deformation is not in contradiction with the basic physical requirements}, which is helpful in proving existence of minimizers. The downside of this approach is that there are not many physically justified models of non-simple materials and material constants are rarely available.  

For non-simple materials of the second grade, we define an energy functional 
\begin{align}\label{non-simple-functional}
\mathcal{E}(u):=\int_\O W(\nabla u(x),\nabla^2 u(x))\,\md x-\int_\O  b(x)\cdot u(x)\,\md x -\int_{\Gamma_{\rm N}} \left( g(x)\cdot u(x)+\hat g_1(x)\cdot \frac{\partial u(x)}{\partial\varrho}\right)\,\md S\ ,
\end{align}
where $\varrho$ is the outer unit normal to $\Gamma_{\rm N}$ and $\hat g_1:\Gamma_{\rm N}\to\R^n$ is the surface density of (hypertraction) forces balancing the {\it hyperstress}
\begin{align}
x\mapsto\frac{\partial}{\partial G_{ijk}} W(F,G)|_{F=\nabla u(x),\,G=\nabla^2u(x)}\ .
\end{align} 

The corresponding first Piola-Kirchhoff stress tensor is constructed as follows. 

Denote for $i,j\in\{1,\ldots, n\}$ $$H_{ij}(F,G):= \sum_{k=1}^n\frac{\partial}{\partial G_{ijk}} W(F,G)\ .$$

Then for $x\in\O$, $F:=\nabla y(x)$,  and $G:=\nabla ^2 y(x)$ we evaluate 
the first Piola-Kirchhof stress tensor as 
$$
S_{ij}(x)=\frac{\partial W(\nabla u(x),\nabla^2u(x))}{\partial F_{ij}}-H_{ij}(\nabla u(x),\nabla^2u(x))\ . $$

We will assume that 
\begin{align}\label{forces-non-simple} u\mapsto \int_\O  b(x)\cdot u(x)\,\md x+\int_{\Gamma_{\rm N}} (g(x)\cdot u(x)+\hat g_1(x)\cdot \frac{\partial u(x)}{\partial\varrho})\,\md S\end{align}
is a linear  continuous functional evaluating the work of external forces on the specimen.
 Here we, however, assume for simplicity  that $b$, $\hat g_1$, and $ g$ depend only on $x\in\Omega$ and $x\in\Gamma_{\rm N}$, respectively.   
Notice that existence of minimizers of $\mathcal{E}$ is guaranteed by Corollary \ref{2nd-grade}.

Similarly, as in the case of simple materials, we may formally derive the Euler-Lagrange equations for minimizers of $\mathcal{E}$.  Interestingly, second-grade materials together with   suitable integrability  of $1/\det\nabla u$ imply a strictly  positive lower bound on $\det\nabla u$ on the whole closure of $\O$.  This enables us to show that minimizers of the energy functionals are weak solutions to the corresponding Euler-Lagrange equations; cf.~\cite{healey-kroemer}.  
Contrary to the simple-material situation, here the smoothness of $\partial\O$ is important because the mean curvature   of the boundary enters the equations.
Details on surface differential operators can be found, for example in \cite{GuiCaf90SIPE}.

\medskip


\section{Weak lower semicontinuity in general hyperelasticity} 
\label{sec-quasicon-elas}
We have seen in the last section that polyconvexity is relatively  easy to be verified and it  ensures weak lower semicontinuity of the corresponding energy functional. Nevertheless, there are materials  that cannot be modeled by polyconvex energy densities.

\begin{figure}[h]
\center
\begin{minipage}{0.45\textwidth}
\center
\hspace*{-2.8cm}\includegraphics[height=0.4\textheight]{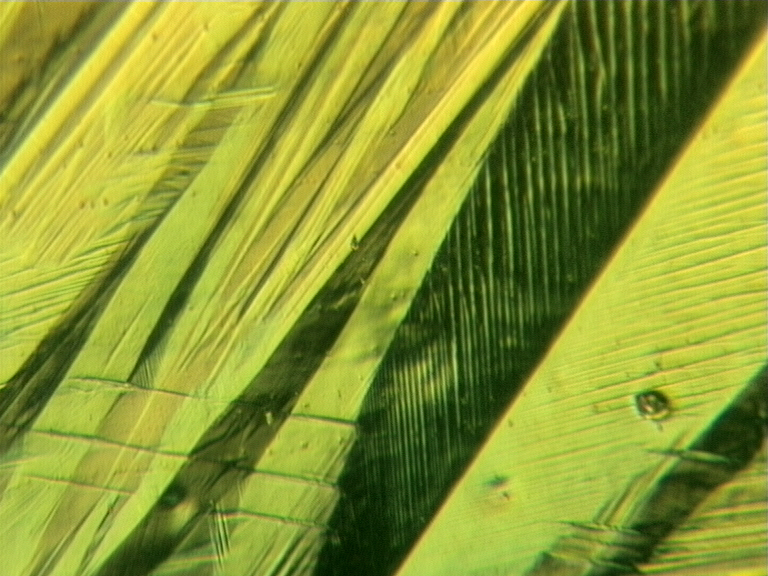}
\end{minipage}
\hspace*{-1.5cm}
\caption{Laminated microstructure in CuAlNi. Courtesy of P.~\v{S}ittner (Inst. of Physics, CAS, Prague)}
\label{figure1}
\end{figure}

A prototypical example are systems featuring phase transition with each phase characterized by some specific deformation of the underlying atomic lattice. This setup is for example found in \emph{shape-memory alloys} (see e.g. the monographs \cite{bhattacharya-1,dolzmann,fremond-1,fremond-miyazaki,Pitteri}, or a recent review \cite{jani}). Shape memory alloys are intermetallic materials which have a high-temperature  highly symmetric phase called austenite and a low temperature phase called martensite which can, however, exist in several variants. Such systems are (for a suitable temperature range) typically modeled by a multi-well stored energy of the form
\begin{align}
W(F)
\begin{cases}
 =0 &\text{if $F=QU_i$ for some $i=1 \ldots M$, and some $Q \ \in \mathrm{SO}(n)$,}\\
> 0 &\text{otherwise}
\end{cases}
\label{energy-sma}
\end{align}
where $U_1, \ldots, U_M\in\R^{n\times n}$ are  given matrices representing the phases found in the material and $\mathrm{SO}(n)$ is the set of rotations in $\R^{n \times n}$.  These materials form complicated patterns (microstructures) composed from different variants
of martensite, cf.~Figure~\ref{figure1}.

 Energy densities  satisfying \eqref{energy-sma} are generically not polyconvex nor quasiconvex and their construction is a modeling issue \cite{zimmer}.  Therefore, {\color{Green} in order to design} an appropriate model one has to find the \emph{weakly lower semicontinuous envelope} of \eqref{functional0} with an energy density given by \eqref{energy-sma}; in other words, one seeks the supremum of weakly lower semicontinuous functionals lying below the given energy. We refer also to the Subsection \ref{sect-ralaxation} for more details on how this \emph{relaxation} of the problem may be performed. Let us remark that {\color{Green} the necessity to perform a relaxation} persist even if we used  a geometrically linear description of energy wells; see  \cite{che-bha,che-sch,govindjee,kohn}, for instance. 

In order to find the weakly lower semicontinuous envelope of \eqref{energy-sma}, a precise characterization of weak lower semicontinuity in terms of convexity conditions on $W$ is needed. We have found these conditions in Section \ref{result-meyers}; however, only under the growth condition (i) in Definition \ref{class-F_R}. Yet, this is incompatible with the physical assumptions formulated in \eqref{eq:negdet}-\eqref{det2zero}. {\color{black} Indeed, notice that \eqref{eq:negdet}
requires the stored energy to be infinite if $\det F \leq 0$, while growth condition (i) in Definition \ref{class-F_R} corresponds to $p$-growth which, roughly speaking, means that the stored energy is bounded from above by a polynomial of the $p$-th order and so is finite for all its arguments.}

For energies taking infinite values it is no longer known that quasiconvexity implies weak lower semicontinuity. Indeed, this is one of the standing problems in elasticity, which was formulated by J.M.~Ball in the following way:
\medskip

\begin{problem}[Problem 1 in \cite{ballOpen}]
\label{problem1}
``Prove the existence of energy minimizers for elastostatics for quasiconvex stored-energy functions satisfying \eqref{det2zero}.''
\end{problem}

\medskip 

{\color{black} While Problem \ref{problem1} is motivated by elasticity, it leads to a study of weak lower semicontinuity that is interesting for its own right. Therefore, we shall denote the integrands by $f$ in this section as elsewhere in the paper and use the notation $W$ only if specifically the stored energy of an elastic material is meant.  

In order to see better why the methods from Section \ref{result-meyers} might fail once one wants to consider integrands taking infinite values, let us revisit Example \ref{ex-intro-quasi} given in the introduction. 

\medskip 

\begin{example}[Example \ref{ex-intro-quasi} revisited]
\label{example-infinite}
In Example \ref{ex-intro-quasi}, we studied weak lower semicontinuity of the functional 
$$
\tilde{I}(u) = \int_\Omega f(\nabla u(x)) \md x,
$$
along the sequence $\{u_k\}_{k \in \N} \subset W^{1,\infty}(\O;\R^m)$ such that $u_k \wstar Ax$ with $A$ some matrix in $\R^{m\times n}$ and saw that if $f$ is quasiconvex and finite on $\R^{n \times m}$ then $\tilde{I}$ is weakly lower semicontinuous. A key step in the proof was the \emph{cut-off} procedure; i.e. the construction of a modified sequence
$$
u_{k, \ell}(x) = \eta_\ell u_k + (1-\eta_\ell) Ax \quad \text{so that} \quad \nabla u_{k, \ell}(x) = \eta_\ell \nabla u_k + (1-\eta_\ell) A + (u_k-Ax)\otimes\nabla \eta_\ell;
$$
that fulfills $u_{k, \ell}(x) = Ax$ on $\partial \Omega$ as well as the realization that once $u_{k, \ell}(x)$ is used as a test function in \eqref{quasiconvexity}, we have that
\begin{equation}
f(A)\mathcal{L}^n(\O) \leq \int_\Omega f(\nabla u_{k(\ell), \ell}(x)) \md x = \int_{\Omega} f(\nabla u_{k(\ell)}) \md x + \int_{\Omega\setminus\Omega_\ell} f(\nabla u_{k(\ell), \ell}(x)) - f(\nabla u_{k(\ell)}(x))  \md x,
\label{eq-examplerevis}
\end{equation}
where the last integral vanishes if $f$ is continuous on $\R^{m \times n}$ because the measure of $\Omega\setminus\Omega_\ell$ converges to zero. 

Now, if $f$ is allowed to take infinite values, the last step is no longer justified. In fact, even if we assumed that the matrix $A$ as well as the sequence $\{u_k\}_{k \in \N}$ are taken such that $f(A) < \infty$ as well as $f(\nabla u_k(x)) < \infty$ for a.e. $x \in \Omega$, there is, a-priori, no guarantee that $f(\nabla u_{k, \ell}(x)) = +\infty$ on a set of positive measure in $\Omega\setminus\Omega_\ell$ can be excluded. However, in this case, $\int_{\Omega\setminus\Omega_\ell} f(\nabla u_{k(\ell), \ell}(x)) - f(\nabla u_{k(\ell)}(x))  \md x = + \infty$ so that, in particular, this term cannot vanish.
\end{example}

\medskip

Example \ref{example-infinite} presents a simplified setting that, however, contains all the difficulty in proving that quasiconvex integrands taking infinite values are lower-semicontinuous. In fact, to the best of authors' knowledge, all proofs available in literature showing that quasiconvexity implies weak-lower semicontinuity, starting from the works of Morrey \cite{morrey-orig} and Meyers \cite{meyers}, are build-up in the same principle as the one presented in Examples \ref{ex-intro-quasi} and \ref{example-infinite} (see also the monograph \cite{dacorogna}). In particular, the cut-off method plays a crucial role, which, as seen in Example \ref{example-infinite}, might result in the consequence that the method of proof becomes unusable for integrands taking values in $\R \cup \{+ \infty\}$. Here, the strength of polyconvexity should be highlighted again as for functionals with polyconvex integrands weak lower semicontinuity is proved by a different method; see Section \ref{sec-polyconvexity}.

\medskip

\begin{remark}
\label{remark-convex-cutoff}
Let us note that, if the set on which $f$ is finite is \emph{convex} the procedure in Example \ref{example-infinite} might still work. Indeed, in this case, if we assumed that $A$ as well as $\{u_k\}_{k \in \N}$ are taken such that $f(A) < \infty$ as well as $f(\nabla u_k(x)) < \infty$ for a.e. $x \in \Omega$ then $f(\eta_\ell \nabla u_k + (1-\eta_\ell) A) < \infty$ a.e. on $\Omega$ due to convexity. This is still does not mean that $f(\nabla u_{k, \ell}(x)) < \infty$ due to the small error $(u_k-Ax)\otimes\nabla \eta_\ell$ due to which $\nabla u_{k, \ell}(x)$ might no longer lie in the set where $f$ takes finite values. Nevertheless, by suitable scaling and continuity of $f$ on its domain, this difficulty can be overcome at least if the domain of $f$ is a ball; see \cite{k-p1}. On the other hand, if the domain of $f$ is \emph{not convex} we can expect that $f(\eta_\ell \nabla u_k(x) + (1-\eta_\ell) A) = \infty$ for some $x\in \Omega$ even if $f(A), f(\nabla u_k(x))$ were finite, which puts us into the situation discussed in Example \ref{example-infinite}.
\end{remark}

\medskip

The physical requirements dictated by \eqref{eq:negdet}-\eqref{det2zero} force the stored energy to be finite only on matrices with a positive determinant, which is a \emph{non-convex} set. In view of the above remark, this means that one faces exactly the situation from Example \ref{example-infinite}, in which the available techniques of proving that quasiconvex energies are weakly lower semicontinuous reach their limit. This is also what makes Problem \ref{problem1} so challenging. 
}
 
\medskip

\begin{remark}
Notice that if \eqref{energy-sma} is additively enriched by a convex term of the form $\varepsilon \int_\Omega |\nabla^2 y|^p \md x$, which is usually interpreted as some kind of interfacial energy of the microstructure, Corollary~\ref{2nd-grade} can be readily applied to show the existence of minimizers for $\mathcal{E}$.

 Let us also point out that a different approach has been proposed recently \cite{silhavy1,silhavy2}. There, a new notion  of {\it interface polyconvexity} has been introduced  which enables to prove existence of minimizers for simple materials with an additional phase field variable. 
\end{remark}

\medskip

\begin{remark}
It has been pointed out in \cite{ballOpen,ball-puzzles} that one of the reasons why Open problem \ref{problem1} is hard to solve is the fact that quasiconvexity possesses no local characterization \cite{kristensen-local}.
\end{remark}

\medskip 

Let us stress that Problem \ref{problem1} is an important attempt  towards combining quasiconvexity and elasticity but additional steps are still required. Namely, if $u:\Omega \to \R^m$ entering \eqref{functional0} ought to represent a deformation of a physical body, it should be \emph{injective} and \emph{orientation-preserving}. {\color{Green} It is even natural to require that the deformation is a \emph{homeomorphism}, at least if one restricts the attention to the case in which formation of cavities and cracks is excluded; we refer to, e.g., \cite{henao} for problems connected to the appearance of cavities.} Notice that this is not automatically satisfied for all maps on which the functional \eqref{functional0} is finite even if $f$ fulfills \eqref{eq:negdet}-\eqref{det2zero}. However, we may rely on Theorem \ref{jmball} to assure this, provided suitable coercivity of the energy.  

An alternative (and related) approach is to study directly weak lower semicontinuity along sequences found in a suitable class of mappings that are injective and orientation-preserving. As a first step, one may study classes of functions that fulfill some constraint on the Jacobian, e.g. that $\mathrm{det}\nabla u > 0$. 

\medskip
{\color{black}
\begin{remark}
Clearly, the positivity  of the Jacobian $\mathrm{det}\nabla u > 0$ does not imply that the mapping $u$ is a homeomorphism. However, surprisingly, also the converse implication does not necessarily hold, even though a homeomorphism on a domain is necessarily sense-preserving (or sense-reversing) in the topological sense. We refer to an example by Goldstein and Haj\l asz \cite{goldstein} where a sense-preserving homeomorphism with $\det \nabla u = -1$ a.e. in $\Omega$ was constructed. This phenomenon concerns functions with a low Sobolev differentiability. Indeed, it was shown by Hencl and Mal\'{y} \cite{hencl-maly} that any sense-preserving homeomorphism in $W^{1,p}(\Omega; \R^m)$ satisfies $\mathrm{det}\nabla u > 0$ a.e. on $\Omega$ if $p > \left[ \frac{n}{2} \right]$.\footnote{Here $[\cdot]$ denotes the integer part.} Recently, it was shown by Campbell, Hencl and Tengvall \cite{campbell} that the exponent is critical in the sense that if $n \geq 4$ and $p < \left[ \frac{n}{2} \right]$ then there exists a homeomorphism in $W^{1,p}(\Omega; \R^m)$ the Jacobian of which is positive and negative on a set of positive measure, respectively.
\end{remark}
}
\medskip

Even though Problem \ref{problem1} remains widely open to date, it has recently  been  approached  from different perspectives. We review the results within this section. {\color{Green} All the results presented in this section are based on changing the cut-off technique introduced in Examples \ref{example-intro-quasi}, \ref{example-infinite} in order to avoid the convex averaging.}

In \cite{krw,krw2013}, {\color{black} Koumatos, Rindler and Wiedemann} study weak lower semicontinuity along sequences in $\{u_k\} \subset W^{1,p}(\Omega; \R^m)$ with $p < n$ satisfying that $\mathrm{det}\nabla u_k > 0$. They proved that \eqref{functional0} with ${\color{black} f}=f(x,\nabla u)$ is weak lower semi-continuous along such sequences if and only if it is $W^{1,p}$-\emph{orientation preserving quasiconvex}, i.e., for almost all $x\in\O$
$$
f(x,A) \leq \frac{1}{\mathcal{L}^n(\O)} \int_\Omega f(x, \nabla \varphi(z)) \md z,
$$
for all $A$ with $\mathrm{det}(A) > 0$, all $\varphi \in W^{1,p}(\Omega; \R^m)$ satisfying that $\varphi(z) = Az$ on $\partial \Omega$ and $\det\nabla \varphi(z) > 0$ for a.a. $z \in \Omega$.

However, in \cite{krw2013} the authors also show that, in fact, for $p < n$ \emph{no} $W^{1,p}$-orientation preserving quasiconvex integrands exist that would satisfy the natural coercivity/growth condition 
$$
\frac{1}{c}\big(|A|^p+  \kappa(\mathrm{det}A)\big) \leq f(x,A) \leq {c}\big(|A|^p+  \kappa(\mathrm{det}A)\big)
$$
for almost all $x \in \Omega$. Here $c>0$ is a constant and $\kappa>0$ is a convex function satisfying that $\lim_{s \to 0} \kappa(s) = +\infty$, $\kappa(s) = +\infty$ for $s \leq 0$ and $\limsup_{s \to \infty} \frac{\kappa(s)}{s^{p/n}} < +\infty$. Notice that this growth condition is compatible with \eqref{eq:negdet}-\eqref{det2zero}.  

The proof in \cite{krw} is based on the  so-called \emph{convex integration}, a technique for solving differential inclusions. It goes back to {\color{black} Nash \cite{nash}, Kuiper \cite{kuiper1,kuiper2} and later} Gromov \cite{gromov} and it found applications in various problems including continuum mechanics ({\color{Green} interestingly in continumm mechanics of solids as well as fluid dynamics}) and regularity theory; see e.g.~{\color{black} \cite{delellis,shnirelman,sheffer,muller-sverak}}. We refer  also to the monograph \cite{dacorogna-marcellini} where solutions to partial differential inclusions by means of Baire category methods are introduced.

To the best of our knowledge, the only works in which the authors actually considered equivalent characterization of weak lower semicontinuity for \emph{injective maps} are  \cite{bbmk13} and \cite{bbkamp15} where bi-Lipschitz and quasiconformal maps \emph{in the plane}  are studied, respectively. 

Here, by bi-Lipschitz maps the following set is meant
\begin{align}
W^{1,\infty,-\infty}_+(\Omega;\mathbb{R}^2) = \Big\{&u: \Omega \mapsto u(\Omega) \text{ an orientation preserving homeomorphism}; \nonumber \\ &u\in W^{1,\infty}(\O;\R^2) \text { and } u^{-1}\in W^{1,\infty}(u(\Omega);\R^2) \text{ is Lipschitz} \Big\}\ ,
\label{deformations}
\end{align}
while quasiconformal maps are introduced as follows
\begin{align}
\mathcal{QC}(\Omega; \R^2) = \Big\{& u \in W^{1,2}(\O;\R^2): \text{ $u$ is a homeomorphism and } \exists { K\geq 1} \text{ such that } \nonumber \\ &|\nabla u|^2 \leq K \det \nabla u \text{ a.e. in $\Omega$} \Big\}. \label{qc-def}
\end{align}

It is natural to expect that weak lower semicontinuity of the functional 
$$
I(u)= \int_\Omega f(\nabla u) \md x,
$$
along sequences in $W^{1,\infty,-\infty}_+(\Omega;\mathbb{R}^2)$ or  $\mathcal{QC}(\Omega; \R^2)$ is connected with a suitable notion of quasiconvexity of $f$. One even expects a \emph{weaker} notion than the one from Definition \ref{def-quasiconvexity} since the set of possible sequences along which semicontinuity is studied is restricted. Indeed, the perfectly fitted notion to this setting seems to be an alternation of Definition \ref{def-quasiconvexity} where only function from $W^{1,\infty,-\infty}_+(\Omega;\mathbb{R}^2)$ or  $\mathcal{QC}(\Omega; \R^2)$ enter as test functions. Exactly this  has been achieved in \cite{bbmk13} and \cite{bbkamp15}; we review the result in Proposition \ref{prop-wlsc-biLip}.

First, let us introduce a notion of weak convergence on $W^{1,\infty,-\infty}_+(\Omega;\mathbb{R}^2)$ and  $\mathcal{QC}(\Omega; \R^2)$. We say that $u_k \wstar u$ in $W^{1,\infty,-\infty}_+(\O;\R^2)$ if the sequence has uniformly bounded bi-Lipschitz constants\footnote{Notice that a function $u \in W^{1,\infty,-\infty}_+(\Omega;\mathbb{R}^2)$ satisfies 
for all $x_1,x_2\in\O$ 
\begin{equation}\label{bi-li}
\frac{1}{\ell}|x_1-x_2|\le |u(x_1)-u(x_2)|\le \ell|x_1-x_2|\ .
\end{equation}
for some $\ell \geq 1$. This $\ell$ is then called the bi-Lipschitz constant of $u$.} and $u_k \wstar u$ in $W^{1,\infty}(\O;\R^2)$. Note that the weak limit is bi-Lipschitz, too.  

For a sequence $\{u_k\}_{k \in \mathbb{N}} \subset \mathcal{QC}(\Omega; \R^2)$, we say that it converges weakly to $u \in W^{1,2}(\Omega; \mathbb{R}^2)$ in $\mathcal{QC}(\Omega;\R^2)$ if $u_k \rightharpoonup u$ in $W^{1,2}(\Omega; \mathbb{R}^2)$, there exists a $K \geq 1$ such that the $u_k$ are all $K$-quasiconformal and $u$ is non-constant. Here it is important to assume that the limit function is non-constant for otherwise the limit function may not {\color{black} be }quasiconformal.\footnote{ Because a sequence of uniformly $K$-quasiconformal maps converges locally uniformly either to a $K$-quasiconformal function or a constant \cite{AstalaIwaniec}. and the locally uniform convergence is implied by the notion of weak convergence in $\mathcal{QC}(\Omega;\R^2)$ } 

Moreover, let us introduce the notions of \emph{bi-quasiconvexity} and \emph{quasiconformal quasiconvexity}.

\medskip

\begin{definition}
We say that a Borel measurable and bounded from below function $f: \R^{2 \times 2} \to \Omega$ is bi-quasiconvex if 
\be\label{bi-def}
\mathcal{L}^2(\O)f(A)\le\int_\O f(\nabla\varphi(x))\,\md x\ 
\ee 
for all $\varphi \in W^{1,\infty,-\infty}_+(\Omega;\mathbb{R}^2)$, $\varphi = Ax$ on $\partial \Omega$ and all $A$ with $\det A > 0$.

We say that $f$ is quasiconformally quasiconvex if \eqref{bi-def} holds for all $A$ with $\mathrm{det}(A) > 0$. and all $\varphi \in \mathcal{QC}(\O;\R^2)$ such that $\varphi(x) = Ax$ on $\partial \Omega$.
\end{definition}

\medskip 

Then we have the following result:

\medskip 

\begin{proposition}[{\color{black} due to Bene\v{s}ov\'{a}, Kampschulte, Kru\v{z}\'{i}k} \cite{bbmk13,bbkamp15}]
\label{prop-wlsc-biLip}
Let $\Omega \subset \R^2$ be a bounded Lipschitz domain. Let ${\color{black} f}$ be continuous on the set of matrices $2\times 2$ with a positive determinant. Then ${\color{black} f}$ is bi-quasiconvex if and only if 
$$u\mapsto I(u)=\int_\O f(\nabla u(x))\,\md x$$  is sequentially weakly* lower semicontinuous on $W^{1,\infty,-\infty}_+(\O;\R^2)$.

Moreover, let $f$ satisfy
$$
0 \leq f(A) \leq c(1+|A|^2) \qquad \text{with $c>0$}
$$
on the set of matrices with a positive determinant. Then ${\color{black} f}$ is quasiconformally quasiconvex if and only if $I$ is weakly lower semicontinuous on $\mathcal{QC}(\O;\R^2)$.
\end{proposition}

\medskip 


{\color{black} Let us shortly comment on the proof Proposition \ref{prop-wlsc-biLip} given in \cite{bbmk13} and \cite{bbkamp15}; it is based on finding a cut-off technique that can cope with the non-convexity of the set of homeomorphisms. The idea is that constructing a cut-off is very much related to understanding the trace operator. Indeed, generally speaking, we may formulate the cut-off problem as follows: given a Lipschitz domain $\Omega$ find another domain $\Omega_\delta \subset \Omega$ with $\vert \Omega\setminus \Omega_\delta \vert \leq \delta$ and a function (or \emph{deformation}) in the considered set (here $W^{1,\infty,-\infty}_+(\Omega;\mathbb{R}^2)$ or $\mathcal{QC}(\Omega; \R^2)$) such that it takes some prescribed values on $\partial \Omega$ and in $\Omega_\delta$. Reformulating this once again, we might ask to find a function from the considered set of functions on $\Omega\setminus\Omega_\delta$ (thus here a function in $W^{1,\infty,-\infty}_+(\Omega\setminus\Omega_\delta;\mathbb{R}^2)$ or $\mathcal{QC}(\Omega\setminus\Omega_\delta; \R^2)$) that has some prescribed boundary values on $\partial (\Omega\setminus\Omega_\delta)$, i.e. on $\partial \Omega$ and $\partial \Omega_\delta$. It is clear that not all boundary data will admit such an extension (even not all smooth data); to see this recall that the Jacobian is a null-Lagrangian (cf. Section \ref{sec-NullLag}) which, for example, immediately excludes all affine mappings with negative Jacobian as boundary data.\footnote{Let us remark here the contrast to the situation considered by Koumatos, Rindler and Wiedemann \cite{krw}. In fact it has been shown in \cite{krw} that every smooth function on the boundary of a Lipschitz domain admits an extension in $y \in W^{1,p}(\Omega;\R^n)$ with $\det \nabla y > 0$ a.e. in $\Omega$ if $p < n$.}   

The characterization of the trace operator on sets  $W^{1,\infty,-\infty}_+(\Omega;\mathbb{R}^2)$ as well as $\mathcal{QC}(\Omega; \R^2)$ is $\Omega$ due to \cite{dan-prat-ext,tukia-ext} and \cite{BeurlingAhlfors}, respectively. However, the works \cite{dan-prat-ext,tukia-ext,BeurlingAhlfors} consider only special geometries of $\Omega$; for example $\Omega$ can be chosen as a square but they are not suited for a doubly-connected domain like $\Omega\setminus \Omega_\delta$. This difficulty has been overcome in \cite{bbmk13,bbkamp15} by suitably meshing $\Omega \setminus \Omega_\delta$ and by defining the cut-off on the grid of the mesh.
}
Even though Proposition \ref{prop-wlsc-biLip} provides us with a weak lower semicontinuity result, this is not yet enough to prove existence of minimizers for functionals with densities from some suitable class. This is so, because bi-Lipschitz as well as quasiconformal maps include a $L^\infty$-type constraint which can be enforced by letting the stored energy density be finite only on a suitable subset of $\R^{2 \times 2}$; yet, this subset is usually left when employing cutoff methods---this happens even in the standard {\color{Green} cases where, however, the issue can be solved by scaling (see Remark \ref{remark-convex-cutoff})}. Thus letting ${\color{black} f}$ being infinite on {\color{Green} a subset of $\R^{2 \times 2}_+$} is incompatible with the proof of Proposition \ref{prop-wlsc-biLip}. 

The usual remedy for proving existence of minimizers or relaxation results is to work with $L^p$-type (with $p$ finite) constraints only. In the setting  above, this would mean to work with so-called \emph{bi-Sobolev} classes (see e.g. \cite{biSobolev}) for $1 < p < \infty$: 
\begin{align*}
W^{1,p,-p}_+(\Omega;\mathbb{R}^2) = \Big\{&u: \Omega \mapsto u(\Omega) \text{ an orientation preserving homeomorphism}; \nonumber \\ &u\in W^{1,p}(\O;\R^2) \text { and } u^{-1}\in W^{1,p}(u(\Omega);\R^2) \Big\}\ .
\end{align*}
However, for these classes of functions, the approach from \cite{bbmk13} and \cite{bbkamp15} cannot be adopted since, {\color{black} as we explained above, it relies on having} a complete characterization of the trace operator {\color{black} which}  is missing to date on these classes. In fact, we have the following

\medskip 

\begin{problem}
Characterize the class of functions $\mathcal{X}(\partial \Omega; \R^2)$ such that
$$
\mathrm{Tr}: W^{1,p,-p}_+(\Omega;\mathbb{R}^2) \stackrel{\text{onto}}{\longrightarrow} \mathcal{X}(\partial \Omega; \R^2)
$$
with $1 < p < \infty$ at least for $\Omega$ being the unit square.
\end{problem}

\medskip 

Let us note that the above problem may play a role also in constructing smooth approximationa (by diffemorphisms) of deformations in elasticity. Indeed, the standard techniques of smoothing Sobolev functions (by a mollification kernel) fail under the injectivity requirement since they essentially rely on convex averaging.

Recently, several results on smoothing even under these constraints appeared \cite{iwaniec1,daneri-pratelli,mora1,mora2,mora3} using completely different techniques and limiting their scope to planar deformations. In particular, {\color{black} Iwaniec, Koskela and Onninen} prove in \cite{iwaniec1} that a homeomorphism in $W^{1,p}(\Omega; \R^2)$ can be strongly approximated by diffeomorphisms in the $W^{1,p}$-norm for $p>1$. For $p=1$, this result has recently been extended in \cite{mora3}.

Nevertheless, in elasticity, one might rather be interested in approximating a function in $W^{1,p,-p}_+(\Omega;\mathbb{R}^2)$ \emph{together with its inverse}. To the authors knowledge, the only {\color{Green} results} in this direction are by {\color{black} Daneri and Pratelli} \cite{daneri-pratelli} who showed that bi-Lipschitz maps can be strongly approximated together with their inverse in the $W^{1,p}$-norm for every finite $p$ and {\color{Green} Pratelli \cite{pratelli-w11} who proved that diffeomorphic approximation is possible in $W^{1,1,-1}(\Omega;\R^2)$}. Yet, for functions in $W^{1,p,-p}_+(\Omega;\mathbb{R}^2)$ with {\color{Green} $1 < p<\infty$} the problem remains largely open as mentioned also in \cite{iwaniec1}.

To end this section, let us remark (by formulating several open problems) that the relation of bi-quasiconvexity and the standard notions of quasiconvexity mentioned in this paper is still unexplored. We focus here only on bi-quasiconvexity but similar problems could be formulated also for quasiconformal quasiconvexity, too.

It is clear from the definitions that any function that is quasiconvex on the set of matrices with a positive determinant is also bi-quasiconvex. Moreover, bi-quasiconvexity implies, at least in the plane, rank-one convexity on the set of matrices with a positive determinant.

\medskip 

\begin{remark}
To see why bi-quasiconvexity implies  rank-one convexity on the set of matrices with a positive determinant, we proceed as follows. First, notice that the determinant changes affinely on rank-one lines due to the formula
\begin{equation}
\det(A+\lambda a \otimes n) = \det A \big(1+\lambda n{\cdot}(A^{-1}a)\big),
\label{rank-oneLine}
\end{equation}
where $a$ and $n$ are some arbitrary vectors. Therefore, rank-one convexity on the set of matrices with a positive determinant is really meaningful, since  all matrices on a rank-one line between two matrices with a positive determinant have this property, too.

Next we mimic the proof from \cite[Lemma 3.11 and Theorem 5.3]{dacorogna} showing that quasiconvexity implies rank-one convexity. Without loss of generality, we suppose that $\Omega$ is the unit square and that we want to show rank-one convexity along the line $A + a \otimes e_1$ with $e_1$ the unit vector in the first coordinate. Then we consider the following sequence of mappings
$$
u_n(x) = u_n(x_1,x_2) = \begin{cases} Ax &\text{for $x_1 \in \left[\frac{k}{n}, \frac{k}{n}+\lambda\frac{1}{n}\right)$ for $k = 0\ldots n-1$, } \\
(A + a \otimes e_1)x &\text{for $x_1 \in \left[\frac{k}{n}+\lambda\frac{1}{n}, \frac{k+1}{n}\right)$ for $k = 0\ldots n-1$},
\end{cases}
$$
with some $\lambda \in [0,1]$. Notice that $\{u_n\}_{n \in \N}$ are Lipschitz, injective and that $(\nabla y)^{-1}$ is uniformly bounded and $\mathrm{det}\nabla u$ is bounded away from zero. Thus, $\{u_n\}_{n\in\N}$ is a sequence of uniformly bi-Lipschitz maps that converges weakly to $\lambda Ax + (1-\lambda)(A + a \otimes e_1)x$. We may therefore use the cut-off technique from \cite{bbmk13} to modify the sequence in such a way that it attains exactly the value of the weak limit at the boundary. Then, the same procedure as in \cite[Theorem 5.3]{dacorogna} gives the rank-one convexity.
\end{remark}

\medskip 

In summary, we have the following series of implications 
$$
\text{quasiconvexity on $\R^{2 \times 2}_+$} \Rightarrow \text{bi-quasiconvexity}  \Rightarrow \text{rank-one convexity on $\R^{2 \times 2}_+$},
$$
where we denoted by $\R^{2 \times 2}_+$ the two-times-two matrices with positive determinant. But it is unclear whether some of the converse implications holds, too. We have the following:

\medskip 

\begin{problem}
Does rank-one convexity on $\R^{2 \times 2}_+$ imply bi-quasiconvexity?
\end{problem}

\medskip

\begin{problem}
Does bi-quasiconvexity imply quasiconvexity on $\R^{2 \times 2}_+$?
\end{problem}

\medskip 



\subsection{Relaxation of non(quasi)convex variational problems}
\label{sect-ralaxation}
As we have already seen, mathematical (hyper)elasticity is one area of analysis  where mechanical  requirements are  above the current tools and results available  in  the calculus of variations.  Orientation preservation and injectivity for simple  non-polyconvex materials are prominent examples. Resorting to non-simple materials depending on second-order deformation gradients might seem as a way out.  What is a physically acceptable form 
of the higher-order energy density is, however, a largely open problem. We refer e.g.~to \cite{ball-corral} for a  discussion on this topic.

Another approach is to accept the fact that our minimization problem may have {\it no solution} and  to trace out the behavior of minimizing sequences driving the elastic energy functional to its infimum on a given set of deformations and to read off some effective material properties  out of their patterns. This is the idea of {\it relaxation} in the variational calculus.  {\color{Green} We explain the main ideas  on the following simplified example for which $f$ in \eqref{functional1} depends just on the first gradient.} Assume we want to 
\be\label{relaxation}
\text{ minimize } \mathcal{E}(u):=\int_\O f(\nabla u(x))\,\md x\  \qquad \text{for $u\in\mathcal{Y}$}.
\ee
Here $\mathcal{Y}$ stands for an admissible set of deformations equipped with some topology. In typical situations, $\mathcal{Y}$ is a subset of a Sobolev space and the topology is the weak one on this space.
If no minimizer exists but the infimum is finite we want to find a new  functional 
$\mathcal{E}_R $  defined on some set  $\mathcal{Y}_R$ such that the following properties hold:
\smallskip
\begin{itemize}
\item[(i)] $\min_{\mathcal{Y}_R}\mathcal{E}_R=\inf_\mathcal{Y}\mathcal{E}$,
\item[(ii)] if $\{y_k\}_{k\in\N}\subset\mathcal{Y}$ is a minimizing sequence of $\mathcal{E}$ then its convergent subsequences  converge (in the topology of $\mathcal{Y}_R$) to minimizers of $\mathcal{E}_R$ on $\mathcal{Y}_R$, and
\item[(iii)] any  minimizer of $\mathcal{E}_R$ on $\mathcal{Y}_R$ is a limit of a minimizing sequence of $\mathcal{E}$.
\end{itemize}
\smallskip

Notice that it is already implicitly assumed in  (i) that minimizers of $\mathcal{E}_R$ do exist on $\mathcal{Y}_R$.   Conditions (ii) and (iii)  state that, roughly speaking, there is a ``one-to-one'' correspondence between minimizing sequences of $\mathcal{E}$ and minimizers of $\mathcal{E}_R$. 
If (i)-(iii) hold  we say that $\mathcal{E}_R$ is the relaxation of $\mathcal{E}$  and that $\mathcal{E}_R$  is the {\it relaxed functional}.  The concept of relaxation is also very closely related to $\Gamma$-convergence and $\Gamma$-limits introduced by E.~de Giorgi. We refer to Braides \cite{braides} and Dal Maso \cite{dalmaso} for a modern exposition {\color{Green} and Section \ref{sec:reading} for further references}.  

If $\mathcal{Y}\subset W^{1,p}(\O;\R^n)$ and the continuous stored energy $f:\R^{m\times n}\to\R$ fulfills
\begin{align}\label{growth-relax}
c_0(-1+|F|^p)\le f(F)\le c_1(1+|F|^p)  
\end{align}
with $c_1>c_0>0$, and $1<p<+\infty$
then Dacorogna \cite{dacorogna-relax} showed\footnote{In fact, Dacorogna's result is stated for more general integrands, namely  $|f(F)|\le c(1+|F|^p)$ with  $Qf>-\infty$, In this case, however, fixed Dirichlet boundary conditions must be inevitably  assigned on the whole $\partial\O$. This is again strongly related to condition (ii) in Meyers ' Theorem~\ref{meyers4} and concentrations on the boundary discussed in Section~\ref{sec-BdCon}.} that $\mathcal{Y}_R=\mathcal{Y}$ equiped with the weak convergence  and 
\begin{align}\label{rel}
\mathcal{E}_R(u):=\int_\O Qf(\nabla u(x))\,\md x\ ,\end{align}
where $Qf:\R^{m\times n}\to\R$ is the {\it quasiconvex envelope} (or quasiconvexification) of $f$ being the largest quasiconvex function not exceeding $f$. It can also be evaluated at any $A\in\R^{m\times n}$ as 

\begin{align}\label{qc-envelope}
Qf(A):=\mathcal{L}^n(\O)^{-1}\inf_{\varphi\in W_0^{1,\infty}(\O;\R^m)}\int_\O f(A+\nabla \varphi(x))\,\md x\ .\end{align}

 Notice that the above formula \eqref{qc-envelope} generally holds only for $f$ locally finite (see \cite{dacorogna}).
The definition of $Qf$ does not depend on the (Lipschitz) domain $\O$ but as we see, the calculation of $Qf$ requires to solve again a minimization problem. Not surprisingly, there are only a few cases where $Qf$ is known in a closed form. We wish to point out {\color{black} the works by DeSimone and Dolzmann} \cite{desimone-dolzmann}  where the authors calculated  the quasiconvex envelope of the  stored energy density arising in modeling of nematic elastomers in three dimensions, and {\color{black} LeDret and Raoult} {\cite{ledret-raoult1,ledret-raoult2,raoult} where the quasiconvex envelope of an  isotropic homogeneous Saint-Venant Kirchoff energy density   ($m=n=3$)
$$W(F):=\frac{\mu}{4}|C-\mathbb{I}|^2+\frac{\lambda}{8}({\rm tr}\, C-3)^2\ $$
is derived. Here $\lambda,\mu$ are Lam\'{e} constants of the material, and $C=F^\top F$ is the right Cauchy-Green strain tensor. Notice that $W$ is convex in $C$ but it is not even rank-one convex in $F$. 

As to relaxation of multi-variant materials we refer to \cite{kohn} where a geometrically linear two well-problem is considered such that elasticity tensors corresponding to both wells  are equal. It was later extended in \cite{che-bha} where non-equal elastic  moduli are admitted. 

 Recently, Conti and Dolzmann \cite{conti-dolzmann} {\color{Green} proved that  the expression of the quasiconvex envelope via \eqref{qc-envelope} is valid even for functions taking infinite values and satisfying the growth condition}
$$
\begin{cases}
c_0(-1 +|F|^p+\theta(\det F))\le f(F)\le c_1(1+|F|^p+\theta(\det F)) &\text{ if $\det F>0$},\\
f(F)=+\infty &\text{ otherwise.}
\end{cases}
$$
Here, $c_1>c_0>0$,  $p\ge 1$, and $\theta:(0;+\infty]\to[0;+\infty)$ is a suitable convex function. 
They, however, require for the result to hold that $Qf$ is polyconvex. Needless to say that this assumption is extremely hard to verify.  Results applicable to a generic situation are missing, so far.

\FIRST
As we pointed out in the introduction, when solving a general minimization problem \eqref{var-problem-intro} one is interested not only in the value of the minimum but, often more importantly, in the minimizer. If the minimizer does not exist, we are still interested in the minimizing sequence. For example, if the minimization problems aims to find the stable states of an elastic  material--and these do not exist--it is still reasonable to argue that  the material will be found in a state ``near'' the actual infimum; i.e. that the minimizing sequence contains the information on physically relevant states. Actually, patterns of  minimizing sequences of functionals with a density of the type \eqref{energy-sma} can be linked to observed microstructures as shown in Figure \ref{figure1}.

The downside of the relaxation via the quasiconvex envelope is that it does not ``store'' too much information on the minimizing sequence. Therefore, a tool for relaxation would be valuable that retains some important features of the minimizing sequence by defining ``generalized functions'' for which the limits $\lim_{k\to\infty}\mathcal{E}(y_k)$  where $\{y_k\}\subset\mathcal{Y}$ are evaluated. This is the basic idea of \emph{Young measures} \cite{y}.  Instead of replacing the original integrand by its quasiconvex envelope we extend the original problem defined on vector-valued functions  to a new problem defined on parametrized measures. These measures  enable us  to describe the limit of a weakly converging sequence composed with a nonlinear function and effectively describe the ``patterns" of the minimizing sequence.

Before giving an example (see Example \ref{exa-YM}) on how Young measures can be employed and before discussing further their properties , let us start with the so-called fundamental theorem on Young measures asserting their existence. This result is originally due to L.C.~Young \cite{y} for $L^\infty$-bounded sequences; various versions of the  theorem below valid for  $L^p$ can be found in  \cite{balder,ball3, fonseca-leoni, kruzik-roubicek, schonbek}, for instance.

 
 \medskip 
 
\begin{theorem}[$L^p$-Young measures]\label{ym-thm}
If $\O\subset\R^n$ is bounded and  $\{Y_k\}_{k\in\N}\subset L^p(\O;\R^{m\times n})$, $1\le p<+\infty$ is a bounded sequence then there exists a (non-relabeled) subsequence and a family of parametrized (by $x\in\O$) probability measures $\nu=\{\nu_x\}_{x\in\O}$ supported on $\R^{m\times n}$  such that for every Carath\'{e}odory  integrand $f:\O\times\R^{m\times n}\to\R\cup\{+\infty\}$ which is bounded from below and
$\{f(\cdot,Y_k)\}_{k\in\N}$ is  relatively weakly compact in $L^1(\O)$ it holds that  
\begin{align}\label{youngmeasures1}
\lim_{k\to\infty}\int_\O f(x,Y_k(x))\,\md x=\int_\O\int_{\R^{m\times n}}f(x,A)\md\nu_x(A)\,\md x\ .
\end{align}

Conversely, if $\nu=\{\nu_x\}_{x\in\O}$  is a weakly* measurable family of probability measures supported on $\R^{m\times n}$  and either 
\begin{align}\label{Lp-condition}
&\int_\O\int_{\R^{m\times n}}|A|^p\md\nu_x(A)\,\md x<\infty \quad \text{ for some $1< p<+\infty$} \qquad \text{ or } \\
\label{Linfty-condition}
&\supp \nu_x\subset B(0,r) \text{ for almost all $x\in\O$ and some $r>0$}
\end{align} 
then there is a sequence  $\{Y_k\}_{k\in\N}\subset L^p(\O;\R^{m\times n})$ ($p=+\infty$ if \eqref{Linfty-condition} holds) such that \eqref{youngmeasures1} holds. Moreover, if \eqref{Lp-condition} holds then  $\{Y_k\}_{k\in\N}$ can be chosen such that 
$\{|Y_k|^p\}_{k\in\N}$ is relatively weakly compact in $L^1(\O)$.

\end{theorem}

\medskip

Here the adjective ``weakly* measurable'' means that $x\mapsto\int_{R^{m\times n}}f(A)\md\nu_x(A)$ is measurable  for all $f\in C_0(\R^{m\times n})$. 
The measure $\nu$ from Theorem~\ref{ym-thm} is called an $L^p$-Young measure generated by $\{Y_k\}$.  It easily follows from Theorem~\ref{ym-thm} that 
a weakly* measurable map $\nu=\{\nu_x\}:\O\to \mathcal{M}^1_+(\R^{m\times n})$ ($\mathcal{M}^1_+(\R^{m\times n})$ denotes the set of probability measures supported on $\R^{m\times n}$) is an $L^p$-Young measure for $1\le p<+\infty$, or $p=+\infty$ 
if and only if \eqref{Lp-condition}  or \eqref{Linfty-condition}, respectively holds. 

We already noted in Section \ref{sec-BdCon} that weak convergence can be essentially caused by concentrations of oscillations that can be separated from each other by Decomposition Lemma~\ref{fons}. As for the Young measure, only the oscillating part of the sequence is important. Indeed, two sequences that differ only on a set of vanishing measure generate the same Young measure \cite{pedregal}; this exactly happens for the original sequence $\{y_k\}_{k \in \N}$ and the sequence $\{z_k\}_{k \in \N}$ constructed in Lemma~\ref{fons}.

\medskip

\begin{remark}
There are finer tools than Young measures that have the ability to capture both; oscillations and concentrations in a generating sequence. These are for example Young measures and  varifolds \cite{fmp}. A detailed treatment of such generalized Young measures can be found in Roub\'{\i}\v{c}ek \cite{r}.
\end{remark}

\medskip

In summary, we see that Young measures are an effective tool to capture the asymptotic behavior of a non-linear functional along a oscillating sequence. Moreover, the Young measure carries information about the oscillations in the sequence themselves. To see this, let us return to Example \ref{ex-zigzag} from the Introduction.

\medskip
\begin{example}[Example \ref{ex-zigzag} revisited]
\label{exa-YM}
In Example \ref{ex-zigzag} we constructed a sequence of ``zig-zag'' functions $\{u_k\}_{k\in\N}$ defined by setting
$$
u(x) = \begin{cases} x & \text{if $0\le x\le 1/2$}\\ 
-x+1 & \text{if $1/2\le x\le 1$}\\
\end{cases} 
$$
extending it periodically to the whole $\R$ and letting Let $u_k(x):=k^{-1}u(kx)$ for all $k\in\N$ and all $x\in\R$. 

The sequence of derivatives $\{u_k'\}_{k \in \N}$ is converging in $L^\infty(0,1)$ weakly* to zero. This is so, because the derivatives oscillate between 1 and $-1$ with zero being the mean. We thus see that the weak limit does carry information about the ``mean'' of the oscillating sequence but not between which values the oscillations happened. This information, however, is encoded in the corresponding Young measure; in fact,  $\{u_k'\}_{k \in \N}$  generates the following sum of two Dirac measures $\frac{1}{2} \delta_1 + \frac{1}{2} \delta_{-1}$.
\end{example}
\medskip

The above example illustrates a recurring theme in Young measures: while the weak limit carries information on average values in a oscillating sequence, the Young measure contains more information on ``where'' the oscillations have taken place.

 Usually, Theorem~\ref{ym-thm} is applied to $f(x,A):=\tilde{f}(A)g(x)$ where $\tilde{f}\in C(\R^{m\times n})$, $\lim_{|A|\to\infty}\tilde{f}(A)/|A|^p=0$, and  $g\in L^\infty(\O)$. These conditions make $\{f(\cdot,Y_k)\}_{k\in\N}$  relatively weakly compact in $L^1(\O)$ so that \eqref{youngmeasures1} holds i.e., 
\begin{align}\label{youngmeasures}
\lim_{k\to\infty}\int_\O \tilde{f}(Y_k(x))g(x)\,\md x=\int_\O\int_{\R^{m\times n}}\tilde{f}(A)g(x)\md\nu_x(A)\,\md x\ .
\end{align}

\FIRST
The original  functional introduced in \eqref{relaxation} is then extended by continuity to obtain its relaxed version expressed in terms of Young measures.
 However, there are a few important issues which need to be addressed. First of all, taking a minimizing sequence $\{u_k\}_{k\in\N}\subset\mathcal{Y}\subset W^{1,p}(\O;\R^m)$ for \eqref{relaxation}  we must ensure that $\{f(\nabla u_k)\}_{k\in\N}$ is weakly relatively compact in $L^1(\O)$. This is possible, for example,  if \eqref{growth-relax} holds.
Indeed, if $\{u_k\}\subset W^{1,p}(\O;\R^m)$ is a bounded  minimizing sequence for $\mathcal{E}$ we can assume that $\{|\nabla u_k|^p\}_{k\in\N}$ is relatively weakly compact in $L^1(\O)$ due to the Decomposition Lemma~\ref{fons}. \footnote{Here, recall that gradients of the two sequences introduced in Lemma \ref{fons} generate the same Young measure.} Then  applying Theorem~\ref{ym-thm} to $Y_k:=\nabla u_k$ we get
$$\inf\mathcal{E}=\lim_{k\to\infty}\mathcal{E}(u_k)=\int_\O\int_{\R^{m\times n}}f(A)\md\nu_x(A)\,\md x\ .$$
Yet, an important issue is  that the resulting Young measure is generated by $\{\nabla u_k\}_{k\in \N}$, i.e., by gradients of Sobolev maps. Such measures are called \emph{gradient Young measures} and clearly form a subset of Young measures. Nevertheless, an explicit characterization of this subset is far more involved than with  the characterization of mere $L^p$-Young measures by means of \eqref{Lp-condition}.
Indeed,  the only  known characterization of admissible  measures, called gradient Young measures,  involves quasiconvex functions again \cite{k-p,k-p1,pedregal,mueller,r} which makes the aim of  obtaining a closed formula of relaxation  by means of parametrized  measures mostly unreachable.
The following theorem is a  characterization of gradient Young measures (also called $W^{1,p}$-Young measures) due to D.~Kinderlehrer and P.~Pedregal \cite{k-p,k-p1}. Before we state theorem we define the following set 

$$
Q(p):=
\begin{cases}
\{f:\R^{m\times n}\to\R;\  f\text{ is quasiconvex } \&\ |f|\le c(f)(1+|\cdot|^p),\ c(f)>0\}\  &\text{ if $1< p<+\infty$},\\
\{f:\R^{m\times n}\to\R;\  f\text{ is quasiconvex} \}  &\text{ if $p=+\infty$}.
\end{cases}
$$

\medskip

\begin{theorem}[$W^{1,p}$-Young measures]\label{kipe-thm}
Let $1<p\le+\infty$ and let $\O\subset\R^n$ be a bounded Lipschitz domain. An $L^p$-Young measure $\nu$ is a gradient Young measure if and only if there is 
 $u\in W^{1,p}(\O;\R^m)$ and a set  $\omega\subset \O$, $\mathcal{L}^n(\omega)=0$, such that for all $x\in\O\setminus\omega$
\begin{align}\label{Y-jensen}
 v(\nabla u(x))\le \int_{\R^{m\times n}}v(A)\md\nu_x(A)
\end{align}
for all  $f\in Q(p)$. 
\end{theorem} 

\medskip

Taking $f(A):=\pm A_{ij}$ for $1\le i\le m$ and $1\le j\le n$ in \eqref{Y-jensen} we get that $\nabla u(x)=\int_{\R^{m\times n}}A\nu_x(\md A)$, i.e., that  $\nabla u$ is the first moment (the expectation) of the Young measure $\nu$. Notice that the obtained $\nabla u$ is also the weak limit of every generating sequence of the Young measure $\nu$.

If \eqref{growth-relax} holds then the relaxed problem  looks as follows:
\begin{align}\label{Youngmeasure-relax}
\text{ minimize } \ \mathcal{E}_R(\nu):=\int_\O\int_{\R^{m\times n}}f(A)\nu_x(\md A)\,\md x \ ,
\end{align}
over  all $W^{1,p}$-Young measures generated by $\{\nabla u_k\}_{k\in\N}$ for some arbitrary bounded   $\{u_k\}\subset\mathcal{Y}$. 
These admissible Young measures then form $\mathcal{Y}_R$.

The Jensen-like inequality \eqref{Y-jensen} puts gradient Young measures into duality with quasiconvexity. As an efficient description of the set of quasiconvex functions is not available, it is practically impossible 
to decide whether a given $\nu$ is generated by gradients.

Nevertheless, one can restrict the attention to a subset of gradient Young measures consisting of so-called \emph{laminates} \cite{pedregal}. Laminates are those gradient Young measures that satisfy the Jensen-like inequality \eqref{Y-jensen} even for all functions from the following set
$$
R(p):=
\begin{cases}
\{f:\R^{m\times n}\to\R;\  f\text{ is rank-one convex } \&\ |f|\le c(f)(1+|\cdot|^p),\ c(f)>0\}\  &\text{ if $1< p<+\infty$},\\
\{f:\R^{m\times n}\to\R;\  f\text{ is rank-one convex} \}  &\text{ if $p=+\infty$}.
\end{cases}
$$
This subset can be advantageously exploited in  numerical minimization of \eqref{Youngmeasure-relax}; see, e.g,. \cite{aubry, bartels1,benesova,dolzmann, kruzik} as well as the review paper by Luskin \cite{luskin} on different finite-element approaches for such a treatment and the illustration in Figure \ref{fig-krychle} which shows computations of material microstructures by means of laminates. 

Another possibility is to minimize \eqref{Youngmeasure-relax} over a superset of gradient Young measures. For example, one can require \eqref{Y-jensen} to hold only for null Lagrangians to get a lower bound on the value of the minimizer. This approach corresponds to replacing $f$ in \eqref{relaxation} by its polyconvex envelope; i.e. the biggest polyconvex functions not exceeding $f$. We refer to \cite{bartels2,bartels-kruzik} for details on this approach. We emphasize that the two mentioned bounds are not sharp in general because rank-one convexity, quasiconvexity, and polyconvexity are different. 
\begin{figure}
\begin{center}
\hspace*{-2.3cm}\includegraphics[height=0.3\textheight]{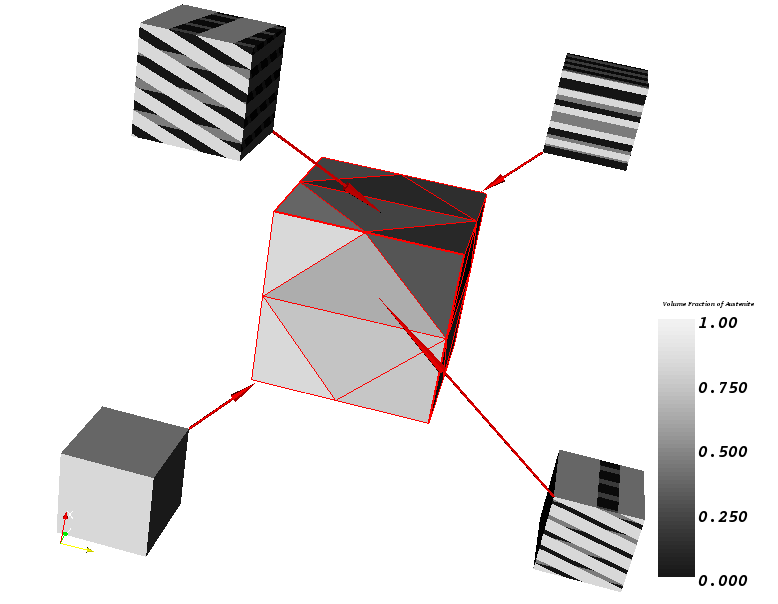}
\end{center}
\caption{An illustration of the calculation of a relaxed energy of a cube under loading. The cube in the middle is the specimen and on the sides the calculated microstructure in form of laminates  is shown in a few elements. The gray scale indicates volume fractions of the phases involved.}
\label{fig-krychle}
\end{figure}

\EOR

Nevertheless, most of the results on relaxation presented above are applicable only to energies with $p$-growth. However, as we already explained in detail in the beginning of this section, this is clearly prevents us to put the requirements \eqref{eq:negdet} and \eqref{det2zero} on $f$. Thus, we see that even current cutting-edge tools and techniques of mathematical analysis and calculus of variations need to be suitably tailored to  cope with deep  problems in elasticity and new ideas and approaches are needed to solve them.

A natural question is whether one can extend weak lower semicontinuity results known for integrands depending on gradients to more a general framework. We will 
investigate this question in the following section.

\EOR

\bigskip


\medskip

\section{$\mathcal{A}$-quasiconvexity}
\label{sec-A-quasi}

In this section, we  summarize results about weak lower semicontinuity of integral functionals  along sequence which satisfy a first-order linear  differential constraint. Clearly, gradients as curl-free fields are included in this setting but, as emphasized by L.~Tartar, besides curl-free fields there are also other  important differential constraints  on possible minimizers. Such a setting naturally arises in electromagnetism, linearized elasticity or even higher-order gradients, to name a few. Tartar's program   was  materialized by Dacorogna   in \cite{dacorogna0} and then  studied by many other authors; see~e.g. \cite{BaCheMaSa13a,braides-fonsea-leoni} and references therein.

The problem studied in this section can be formulated as follows:
Having a sequence 
$\{u_k\}_{k\in \N} \subset L^p(\O;\R^m)$, $1<p<+\infty$ such that each member  satisfies a linear differential constraint ${\mathcal A}u_k=0$ ($\mathcal{A}$-free sequence), or ${\mathcal A}u_k\to 0$ in $W^{-1,p}(\O;\R^n)$ (asymptotically $\mathcal{A}$-free sequence),   what conditions on $v$ precisely ensure  weak lower semicontinuity of integral functionals in the form 
\begin{equation}\label{eq:funct}
\mathcal{I}(u):=\int_\O f(x,u(x))\,\md x\ .
\end{equation}
 Here ${\mathcal A}$ is a first-order linear differential operator.
 
 To the best of our knowledge, the first result of this type was proved  in \cite{fomu} for nonnegative integrands.  In this case, the  crucial necessary and sufficient condition ensuring weak lower semicontinuity of $\mathcal{I}$ in \eqref{eq:funct} is the so-called $\mathcal{A}$-quasiconvexity; cf.~Def.~\ref{a-quasiconvexity} below. However, if we refrain from considering only  nonnegative integrands, this condition is not necessarily sufficient as we already observed in the case $\mathcal{A}:=$curl.

\subsection{The operator ${\mathcal A}$ and  ${\mathcal A}$-quasiconvexity}
Following \cite{fomu},  we consider the  linear operators $A^{(i)}:\R^m\to\R^d$, $i=1,\ldots, n$, and  define ${\mathcal A}:L^p(\O;\R^m)\to W^{-1,p}(\O;\R^d)$ by
\ben
{\mathcal A}u:=\sum_{i=1}^nA^{(i)}\frac{\partial u}{\partial x_i}\ ,\mbox{where }\ u:\O\to\R^m\ ,
\een
i.e.,
for all $w\in W^{1,p'}_0(\O;\R^d)$ 
\begin{eqnarray*}
\la {\mathcal A}u,w\ra=-\sum_{i=1}^n\int_\O A^{(i)}u(x)\cdot\frac{\partial w(x)}{\partial x_i}\,\md x
\ .\end{eqnarray*}
For $w\in\R^n$ we define the linear map 
\ben
{\mathbb A}(w):=\sum_{i=1}^n w_iA^{(i)}\ :\ \R^m\to\R^d\ .\een
In this review, we assume that there is $r\in\N\cup\{0\}$ such that
\be\label{rank}
{\rm rank}\ {\mathbb A}(w)=r\mbox{ for all $w\in\R^n$}\ , |w|=1\ ,\ee
i.e.,  ${\mathcal A}$ has the so-called {\it constant-rank property}. Below we use $\ker$ to denote the set of all locally integrable functions $u$ such that
$\mathcal{A} u=0$ in the sense of distributions, i.e., $\la\mathcal{A}u,w\ra=0$ for all $w\in \mathcal{D}(\O)$.
Of course, $\ker$ depends on the considered domain $\Omega$.
%
%

\medskip 

\begin{definition}[cf.~{\cite[Def.~3.1,~3.2]{fomu}}]\label{a-quasiconvexity}
 We say that  a continuous function $f:\R^m\to\R$, satisfying  $|f(s)|\le C(1+|s|^p)$ for some $C>0$, is ${\mathcal A}$-{\it quasiconvex} if for all $s_0\in\R^{m}$ and all $\varphi\in L^p_\#(Q; \R^m)\cap\ker$ with $\int_{Q}\varphi(x)\,\md x=0$ it holds
\ben
f(s_0)\le \int_Q f(s_0+\varphi(x))\,\md x\ .\een
\end{definition}

\medskip 

In the above definition, we used the space of $Q$-periodic Lebesgue integrable functions:
\[
L^p_\#(\R^n;\R^m):=\{u\in L^p_{\rm loc}(\R^n;\R^m): \mbox{ $u$ is $Q$-periodic}\},
\]
where, $Q$ denotes the unit cube $(-1/2,1/2)^n$ in $\R^n$, and we say that $u:\R^n\to\R^m$ is {\it $Q$-periodic} if for all $x\in\R^n$ and all $z\in\mathbb{Z}^n$ it holds that $u(x+z)=u(x)\ .$

Fonseca and M\"{u}ller \cite{fomu}  proved the following result linking $\mathcal{A}$-quasiconvexity and weak lower semicontinuity. Notice that the integrand is more general than that one in \eqref{eq:funct}.  

\medskip

\begin{theorem}\label{fomu-thm}
Let $\O\subset\R^n$ be open and bounded and let  $f:\O\times\R^d\times\R^m\to[0;+\infty)$ be a Carath\'{e}odory integrand.  Let
$$0\le f(x,z,u)\le a(x,z)(1+|u|^p)$$ for almost every $x\in\O$ and all $(z,u)\in\R^d\times\R^m$, $1<p<+\infty$, and some $0\le a\in L^\infty_{\rm loc}(\O;\R^d)$.  Assume that $\{z_k\}_{k\in\N}\subset L^\infty(\O;\R^m)$,  $z_k\to z$ in measure and that $u_k\wto u$ in $L^p(\O;\R^d)$,$\|\mathcal{A} u_k\|_{W^{-1,p}(\O;\R^m)}\to 0$.  

Then 
$$\liminf_{k\to\infty} \int_\O f(x,z_k,u_k)\,\md x\ge \int_\O f(x,z,u)\,\md x$$   
if and only if  $f(x,z,\cdot)$ is $\mathcal{A}$-quasiconvex for almost all $x\in\O$ and all $z\in\R^d$. 
\end{theorem}

\medskip

The following definition is motivated by our  discussion above  Theorem~\ref{wlsc}. It first appeared in \cite{ifmk}.

\medskip

\begin{definition}\label{extension}
Let $1<p<+\infty$ and $\{u_k\}_{k\in\N}\subset L^p(\O;\R^m)\cap\ker$. We say that $\{u_k\}$ has an $\mathcal{A}$-free $p$-equiintegrable extension if for every domain $\tilde\O\subset\R^n$ such that  $\O\subset\tilde\O$, there is  a sequence $\{\tilde u_k\}_{k\in\N} \subset L^p(\tilde\O;\R^m)\cap\ker$ such that\\
\noindent 
(i) $\tilde u_k= u_k$ a.e.~in $\O$ for all $k\in\N$,\\
\noindent
(ii) $\{|\tilde u_k|^p\}_{k\in\N}$ is equiintegrable on $\tilde\O\setminus\O$, and\\
\noindent 
(iii) there is $C>0$ such that $\|\tilde u_k\|_{L^p(\tilde\O;\R^m)}\le C\|u_k\|_{L^p(\O;\R^m)}$ for all $k\in\N$.  
\end{definition}

\medskip

Then we have the following result proved in \cite{ifmk}.

\medskip

\begin{theorem}\label{wlsc-A}
Let $0\le g\in C(\overline{\Omega})$, let $|f|\le C(1+|\cdot|^p)$ be ${\mathcal A}$-quasiconvex, satisfy \eqref{p-lipschitz}, have a recession function,   and let  $1<p<+\infty$.
Let  $\{u_k\}\subset L^p(\O;\R^m)\cap\ker $, $u_k\wto u$ weakly,  and assume that  $\{u_k\}$ has an $\mathcal{A}$-free  $p$-equiintegrable extension.
Then $\mathcal{I}(u)\le \liminf_{k\to\infty}\mathcal{I}(u_k)$, where
\begin{align}\label{fcional-A} \mathcal{I}(u):=\int_\O g(x)f( u(x))\,\md x. \end{align}

\end{theorem}

Surprisingly, it was shown by Fonseca and M\"{u}ller \cite[p.~1380]{fomu} that also higher-order gradients can be recast as $\mathcal{A}$-free mappings. They construct $\mathcal{A}$ such that $\mathcal{A}u=0$ if and only if $u=\nabla^k w$ for some $w\in W^{k,p}(\O;\R^m)$. In this situation, $\mathcal{A}$-quasiconvexity coincides with Meyers' $k$-quasiconvexity. This allows us to study the weak lower semicontinuity of 
\begin{equation}
\mathcal{I}(w):=\int_\O g(x)f(\nabla^k w(x))\,\md x
\label{A-higherDer}
\end{equation}
with $g$, and $f$ as in Theorem~\ref{wlsc-A}, on the set
$$\{ w\in W^{k,p}(\O;\R^m):\, w=w_0\text{ on $\partial\O$}\}\ $$
with $w_0\in W^{k,p}(\O;\R^m)$, in the context of $\mathcal{A}$-quasiconvexity. In particular, in this case, we may construct the extension required in Definition \ref{extension} by extending the Dirichlet boundary condition and thus Theorem \ref{wlsc-A} is applicable. It follows that $k$-quasiconvexity if $f$ is sufficient to make \eqref{A-higherDer} weakly lower semicontinuous which affirmatively answers Problem \ref{Meyers-conj} in this particular setting.\footnote{In fact, one can consider integrands of the form $f(x,\nabla^k w(x))$ whenever $f(x,\cdot)$ is $k$-quasiconvex for all $x\in\overline{\Omega}$, $|f(x,A)|\le C(1+|A|^p)$,  $f(\cdot; A)$ is continuous in $\overline{\Omega}$ for all $A\in X(n,m,k)$, and $f(x,\cdot)$ possesses a recession function.}

Let us finally point out that treatment of $\mathcal{A}$-quasiconvexity for integrands whose negative part growth with the $p$-th power is a very subtle issue which has recently been treated in \cite{kraemer}. There is a new condition called $\mathcal{A}$-quasiconvexity at the boundary which is introduced in two forms depending whether $u$ can be extended to a larger domain preserving the $\mathcal{A}$-free property or not. This allows us to remove the assumption on the existence of an $\mathcal{A}$-free $p$-equiintegrable extension from Theorem~\ref{wlsc-A}. 

\medskip

\section{Suggestions for further reading}\label{reading}
\label{sec:reading}
\FIRST

The above exposition aims at reflecting developments in  weak lower semicontinuity  related to Morrey's \cite{morrey-orig} and  Meyers' \cite{meyers} papers  with the emphasis on applications to static problems in  continuum mechanics of solids. We dare to  hope that it provides a fairly complete picture of theory starting in 1952  to current trends.  On the other hand, it is clearly influenced by our personal point of view.   We  believe to have  convinced the reader that the calculus of variations has been a very active research area   for last fifty years with many important results and with many persisting challenging open problems. Below we mention some additional research directions and we invite the interested reader to find more details in the listed references.

\subsection{Applications to continuum mechanics and beyond} We saw that weak lower semicontinuity serves as a main ingredient of proofs of existence of minimizers to variational integrals and we outlined applications in elastostatics. Even in the static case, models of elasticity can be combined with other phenomena, as magnetism, for instance. This leads to magnetoelasticity (magnetostriction), a property of NiMnGa, for instance. We refer e.g.~to the monograph \cite{ogden} or \cite{eringen-maugin} for a physical background. 

The idea to draw macroscopic properties of composite  materials from their microscopic ones is in the core of {\it homogenization theory}. We cite the classical books by Jikov, Kozlov, and Olejnik \cite{jikov}, by Braides and Defranceschi \cite{braides-de} or by Bensoussan, Lions, and Papanicolaou \cite{papa} for a thorough overview. $\Gamma$-convergence of integral functionals (see  monographs by  Braides or Dal Maso \cite{braides,dalmaso}) plays a key role in this research.   The main objective is to study 
properties for $\varepsilon\to 0$ of 
$I_\varepsilon(u_\varepsilon):=\int_\O f(x,\frac{x}{\varepsilon},\nabla u_\varepsilon(x))\,\md x$, where 
$f$ is $(0,1)^n$-periodic in the middle variable.
Homogenization problems are not only restricted to gradients but new results in the context of $\mathcal{A}$-quasiconvexity (even with non-constant coefficients) recently appeared, e.g., in \cite{davoli-fonseca}. 
Various other generalizations including stochastic features to describe  randomness in the distribution of inhomogeneities can be found in the papers by Dal Maso and Modica \cite{dalMaso-stoch},  Messaoudi and Michaille \cite{gerard} or the recent work of Gloria, Neukamm, and Otto \cite{gloria-neukamm-otto}.
$\Gamma$-convergence is also one of the main tools to dimension-reduction problems in mechanics to obtain  various plate models, see for instance  LeDret and   Raoult \cite{ledret-raoult3},  Friesecke, James, and M\"{u}ller \cite{fjm} or Hornung, Neukamm and Vel\v{c}i\'{c} \cite{hornung} and references therein.

 Quasistatic evolutionary problems and their treatment by means of the so-called ``energetic solution'' has been a lively area of mathematical continuum mechanics for many years .  This theory is     thoroughly   discussed  and summarized in the book by Mielke and Roub\'{\i}\v{c}ek \cite{mielke-roubicek} where one can find  applications to   continuum mechanics of solids including plasticity, damage, or  mechanics of shape memory materials including numerical approximations. Various combinations of inelastic processes are also treated e.g. in \cite{bonetti, francfort}.   Weak lower semicontinuity of corresponding energy functionals  is a key part of theory. It also includes relaxation and  Young measures.  We also point out the paper  \cite{mielke-roubicek-stefanelli} for (evolutionary) $\Gamma$-convergence treatment of quasistatic problems or  \cite{bonetti} for combinations of damage and plasticity.  This opens many possibilities to apply this notion of solution to dimension-reduction problems, linearization, or brittle damage, for instance.    
Weak lower semicontinuity finds its application in dynamical problems, too. For example, it is the main tool to prove existence of solutions in time-discrete approximations of evolution in various models. We refer e.g.~to the book of Braides \cite{braides-book} or \cite{mielke-roubicek} for  many such instances. We also refer to \cite{gms} and references therein for further results concerning 
mathematical treatment of nonlinear elasticity.  

Derivation of nonlinear continuum models from discrete atomic ones  and corresponding  numerical computations has been a lively field of research in recent years. We refer to Blanc, Le Bris, and Lions \cite{blanc} for a survey article and to work of Braides and Gelli \cite{DC2}, Allicandro and Cicalese \cite{DC1} or Braun, Friedrich and Schmidt \cite{schmidt-1,schmidt-2} and references therein. Various aspects of  continuum/atomistic coupling  which combines accurate atomistic models with efficient continuum elasticity are thoroughly discussed in Luskin and Ortner \cite{luskin-ortner}, see also Lazzaroni, Palombaro Sch\"{a}ffner and Schl\"{o}merkemper \cite{Anja,schaeffner} for further analytical results  in this direction.   

Nonlocal theories of elastodynamics, as e.g. peridynamics, has  initiated research activities in nonlocal  variational problems. Here we wish to mention recent work of Bellido, Mora-Corral \cite{bellido-1} and Bellido Mora-Corral and Pedregal \cite{bellido-2}, for instance.

\subsection{Functionals with linear growth} Contrary to weak lower semicontinuity in $W^{1,1}(\O;\R^m)$ where the definition excludes concentrations in gradient, a serious  analytical problem appears  if one wants  to minimize  functionals defined  in  \eqref{functional0} for $f(x,r,\cdot)$ having linear growth at infinity. In this case, we would naturally work in $W^{1,1}(\O;\R^m)$ which is, however,  not reflexive and therefore the weak limit of a minimizing (sub)sequence does not necessarily exist. 
This leads to various extensions of $W^{1,1}(\O;\R^m)$ as well as of the functional $I$. Usually, we embed $W^{1,1}(\O;\R^m)$ into the space of functions with bounded variations $\mathrm{BV}(\O;\R^m)$ which contains integrable mappings whose gradient is a Radon measure in $\O$. Detailed  description and properties can be found in monographs by Ambrosio, Fusco, and Pallara \cite{ambrosio} or by Attouch, Buttazzo and Michaille \cite{attouch}. Besides the function space one must also suitably extend
the functional $I$ to allow for measure-valued gradients.  This  uses the notion of 
the recession function to $f(x,r,\cdot)$. We refer to Fonseca and M\"{u}ller \cite{fonseca-mueller} and \cite{fonseca-mueller2} and references therein  for the case of nonnegative integrands.  Recently, Kristensensen and Rindler \cite{KriRi10a} resolved the weak* lower semicontinuity  relaxation of the functional  if $f=f(x,\nabla u)$  with $|f(x,s)|\le C(1+|s|)$ along sequences with prescribed Dirichlet boundary conditions.  This result was then generalized by Bene\v{s}ov\'{a}, Kr\"{o}mer and Kru\v{z}\'{\i}k \cite{BKK2013} to avoid restrictions on the boundary.   
A closely related topic is weak* lower semicontinuity in the space of functions with bounded deformations \cite{temam-strang}, i.e., in the subspace of $L^1(\O;\R^m)$ of mappings whose symmetrized gradient is a measure on $\O$. This set   naturally appears in problems of linearized  perfect plasticity \cite{suquet}, for instance. We refer to Rindler  \cite{rindler}  for a general lower semicontinuity statement.

Functionals with linear growth and their setting in the space of maps with bounded variations many times arise as $\Gamma$-limits of various phase-field problems.  Here one tries to accommodate two or more phases/materials of a given volume in such a way that the interfacial energy assigned to two mutual interfaces is minimal. Besides classical works of Modica and Mortola \cite{modica,modica-mortola} on phase transition nowdays  exist versatile applications in various areas of mechanics like topology optimization \cite{garcke1} and many others.

\EOR

\medskip

\noindent
{\bf Acknowledgment.}  We thank Anja Schl\"{o}merkemper, Tobias Geiger, and Luisa Wunner  for their helpful remarks on an earlier version of the manuscript. We are also indebted to  the anonymous referees for their comments and suggestions which thoroughly improved the exposition of the article. This work reflects  long-term research of the authors which   was supported by many grants, in particular by  GA\v{C}R projects 14-15264S, 14-00420S, P107/12/0121, 16-34894L, 13-18652S, M\v{S}MT project 7AMB16AT015, and the DAAD-CAS project DAAD-15-14. A revised version of this paper has been partly prepared when MK hold a visiting Giovanni-Prodi Chair in the Institute of Mathematics, University of W\"{u}rzburg. The support and hospitality of the Institute is gratefully acknowledged.

\end{document}